\documentclass[invmat,
final,
]{bcr-svjour}

\usepackage{amsmath}
\usepackage{amssymb}

\usepackage{
epsfig,
graphicx,
psfrag,
subfigure
}

\newcommand{\dD}{\ensuremath{\mathbb{D}}} 
\newcommand{\dE}{\ensuremath{\mathbb{E}}} 
\newcommand{\dL}{\ensuremath{\mathbb{L}}} 
\newcommand{\dN}{\ensuremath{\mathbb{N}}} 
\newcommand{\dP}{\ensuremath{\mathbb{P}}} 
\newcommand{\dR}{\ensuremath{\mathbb{R}}} 
\newcommand{\dZ}{\ensuremath{\mathbb{Z}}} 

\newcommand{\C}[1]{\ensuremath{{\mathcal C}^{#1}}} 
\newcommand{\PT}[1]{\mathbf{P}_{\!#1}} 
\newcommand{\GI}{\mathbf{L}} 
\newcommand{\Ga}{\boldsymbol{\Gamma}} 
\newcommand{\ind}{\mathrm{1}\hskip -3.2pt \mathrm{I}} 
\newcommand{\NRM}[1]{\ensuremath{{\left\Vert #1\right\Vert}}} 
\newcommand{\ent}{\mathbf{Ent}} 

\spnewtheorem*{assumption}{Assumption}{\bf}{\rm}

\begin{document}   

\makeatletter
 \renewcommand{\@evenfoot}{%
   \hfil{\small%
         \today{}.
   \hfill}}
 \renewcommand{\@oddfoot}{\@evenfoot}
\makeatother


\title{Interpolated inequalities between exponential and Gaussian,
Orlicz hypercontractivity and  isoperimetry.}

\titlerunning{Interpolated inequalities between exponential and Gaussian.}

\author{F. Barthe\inst{1} \and P. Cattiaux\inst{2} \and C. Roberto\inst{3}}

\institute{
Universit\'e Toulouse III, \\
\email{barthe@math.ups-tlse.fr}
\and
Ecole Polytechnique et Universit\'e Paris X,\\
\email{cattiaux@cmapx.polytechnique.fr}
\and
Universit\'es de Marne la Vall\'ee et de Paris XII Val de Marne,\\
\email{roberto@math.univ-mlv.fr}
}
\maketitle

\begin{abstract}
   We introduce and study a notion of Orlicz hypercontractive
   semigroups. We analyze their relations with general $F$-Sobolev
   inequalities, thus extending Gross hypercontractivity theory.  We
   provide criteria for these Sobolev type inequalities and for
   related properties. In particular, we implement in the context of
   probability measures the ideas of Maz'ja's capacity theory, and
   present equivalent forms relating the capacity of sets to their
   measure.  Orlicz hypercontractivity efficiently describes the
   integrability improving properties of the Heat semigroup associated
   to the Boltzmann measures $\mu_{\alpha}(dx) = (Z_{\alpha})^{-1}
   e^{-2|x|^{\alpha}} dx$, when $\alpha\in (1,2)$.  As an application
   we derive accurate isoperimetric inequalities for their products.
   This completes earlier works by Bobkov-Houdr\'e and Talagrand, and
   provides a scale of dimension free isoperimetric inequalities as
   well as comparison theorems.  
\end{abstract}


\section{Introduction.}\label{I}

Sobolev type inequalities play an essential role in the study of the
concentration phenomenon for probability measures.  They are also a
powerful tool to analyze the regularizing effects and the convergence
to equilibrium of their associated symmetric semigroups. In
particular, several surveys deal with the celebrated Poincar\'e (or
spectral gap) inequality and the stronger logarithmic-Sobolev
inequality and provide striking applications \cite{Gross},
\cite{Bak94}, \cite{Ane}, \cite{Led99}, \cite{GZ99}, \cite{Ro99}.

\medskip A concrete illustration can be given for the family of
probability measures on the real line
$$
\mu_{\alpha}(dx) = (Z_{\alpha})^{-1} e^{-2|x|^{\alpha}} dx,\qquad
\alpha>0.
$$
These measures and their products $\mu_{\alpha}^{\otimes n}$
deserved particular attention in recent years, where the focus was on
dimension free properties. They enter Talagrand's work on the
concentration phenomenon for product measures. His study was continued
by Ledoux \cite{Led99}, who strongly put forward the use of the logarithmic Sobolev
inequality for concentration, and more recently by Bobkov and Houdr\'e 
\cite{bobkh97cbis}  who
introduced $\dL_1$-Sobolev type inequalities in order to study the more
delicate isoperimetric problem. We review the main results concerning
these measures and the associated semigroup $(\PT{t}^{\alpha})_{t\ge 0}$
 generated by the operator
$\GI_\alpha$ such that
$$
\GI_{\alpha}f = \frac 12 f'' - \alpha |x|^{\alpha - 1} {\rm
  sign}(x) f'.
$$

For $\alpha >0$ the measures $\mu_{\alpha}$ verify a Weak Spectral
Gap property introduced by Aida and Kusuoka as shown in \cite{RW01}.
They satisfy the Spectral Gap inequality exactly when $\alpha \ge 1$,
and the logarithmic Sobolev inequality if and only if $\alpha\ge 2$.

When $\alpha <1$ there is no dimension free concentration property
\cite{Tal91}.
Oppositely, the measures enjoy very strong properties when $\alpha
>2$.  The corresponding semigroup is ultracontractive \cite{KKR},
meaning that for positive time it is continuous from $\dL^2(\mu_\alpha)$
to $\dL^\infty$.  The measures satisfy a dimension free Gaussian
isoperimetric inequality \cite[Theorem~9]{bart01epch}, and this is as
bad as it gets by the Central Limit Theorem. Recently, Bobkov and
Zegarlinski \cite{BZ02} obtained concentration inequalities for these
measures but for the $\ell_\alpha^{n}$-distance on $\dR^n$.  Their
results are based on appropriate modification of the logarithmic
Sobolev inequality, and show different behaviors for different values
of $\alpha$. This was not the case when considering the Euclidean
distance.

The range $\alpha\in[1,2]$ presents very interesting properties. We
start with the Gaussian case, $\alpha=2$, which is best understood.
Concentration of measure and  isoperimetry in Gauss space are now classical
 (see {\it e.g.} \cite{Led99,BL96}). It is remarkable  that they are both
dimension free. Recall that the isoperimetric inequality asserts in particular
that for $A\subset \dR^{n}$ with $\mu_2^{\otimes n}(A)=\mu_2((-\infty,t])$
one has for all $h>0$,
\begin{equation}\label{isop-mu2}
 \mu_2^{\otimes n}\Big(A+hB_2^{n}\Big)\ge\mu_2\Big((-\infty,t+h]\Big).
\end{equation}
Here $B_2^{n}$ is the $n$-dimensional Euclidean ball. Taking limits one obtains
that among sets of given Gaussian measure, half-spaces have minimal Gaussian 
boundary measure.

On the other hand, the Gaussian measure has remarkable analytic properties:
the corresponding Ornstein-Uhlenbeck semigroup is hypercontractive, as discovered
by Nelson \cite{nels73fmf}. Gross proved that this fact is equivalent to the logarithmic
Sobolev inequality \cite{gross75}. Let us also mention that the Gaussian measure
is the prototype of strictly log-concave measures. It was a  success of the 
Bakry-Emery formalism to allow the extension of most of the previous results
to abstract semigroups with positive curvature (see \cite{bakry-emery}
for logarithmic Sobolev inequalities and  \cite{BL96}
 for  Gaussian isoperimetry and an abstract version of the Levy-Gromov
theorem).

\medskip

The two sided exponential measure $\alpha=1$ is also well understood.
Talagrand's paper \cite{Tal91} provides the following very precise estimate:
if $A\subset \dR^{n}$ verifies $\mu_1^{\otimes n}(A)=\mu_1((-\infty,t])$ then
for all $h\ge 0$ it holds
$$ \mu_1^{\otimes n}(A+h B_1^{n}+\sqrt h B_2^{n})\ge \mu_1((-\infty,t+h/K]),$$
where $K$ is a universal constant and $B_p^{n}=\{x\in \dR^{n};\; \sum_{i=1}^{n}|x_i|^{p}
\le 1\}$. See also \cite{maur91sdi,tala96tcgo}.
 In a slightly weaker form, such a statement was recovered by 
 Bobkov and Ledoux \cite{BL97}, via a modified logarithmic Sobolev inequality which 
is equivalent to Poincar\'e inequality. Thus products of measures on $\dR^{d}$
with a spectral gap satisfy a concentration inequality on the exponential model.
Moreover, Bobkov and Houdr\'e \cite{BH97} proved that $\mu_1^{\otimes n}$ 
satisfies a dimension free isoperimetric inequality of Cheeger. The proof relies
on an $\dL_1$ version of the Poincar\'e inequality, and the statement can be rephrased
as follows: let  $A\subset \dR^{n}$ with $\mu_1^{\otimes n}(A)=\mu_1((-\infty,t])$ then
for all $h\ge 0$ 
\begin{equation}\label{isop-mu1} \mu_1^{\otimes n}\Big(A+h B_2^{n}\Big)\ge
 \mu_1\Big(\big(-\infty,t+\frac{h}{2\sqrt 6}\big]\Big).
\end{equation}
This result completes the one of Talagrand. It is weaker for large values of $h$
but gives isoperimetric information as $h$ goes to zero.

\medskip

Understanding the case  $\alpha \in (1,2)$ is the task of the present paper.
The concentration phenomenon is already well described. Indeed, Talagrand's exponential
inequality  transfers to $\mu_\alpha$ for any $\alpha>1$ \cite{tala94scp} and
ensures that for every $A\subset \dR^{n}$ and $h>0$
$$ \mu_\alpha^{\otimes n}\Big(A+h^{1/\alpha}B_\alpha^{n}+\sqrt h B_2^{n}\Big)\ge 1-
\frac{1}{\mu_\alpha^{\otimes n}(A)} e^{-h/K}.$$
In particular if $\alpha\in (1,2)$ and $\mu_\alpha^{n}(A)\ge 1/2$ one gets that
for $h\ge 1$,
 $\mu_\alpha^{n}(A+hB_2^{n})\ge 1-2e^{-h^{\alpha}/K}$.  
A functional approach to this fact was recently discovered by 
 Lata{\l}a and Oleszkiewicz \cite{LO00}. 
These authors established the following family of Sobolev inequalities:
 there exists a universal constant $C$ such that for all $1 < p < 2$ it holds
\begin{equation}\label{1.1}
  \int f^2 d\mu_{\alpha} - \Big(\int |f|^p  d\mu_{\alpha}\Big)^{\frac 2p}
  \leq
  C (2-p)^{2(1-\frac{1}{\alpha})}  \int (f')^2 d\mu_{\alpha}
\end{equation}
for smooth enough $f$. For $\alpha = 2$ these inequalities are due to
 Beckner \cite{Bec89}.
Inequalities~\eqref{1.1} interpolate between Poincar\'e and log-Sobolev.
They enjoy the tensorisation property and
imply  dimension free concentration  with decay $e^{-K t^{\alpha}}$ as expected.
Obviously  \cite{LO00} was the starting point of an extension of the log-Sobolev
approach to concentration, encompassing  more general behaviors. Recently two of us
simplified the proof of \eqref{1.1} and characterized all measures on $\dR$
satisfying the same property \cite{BR03} (such a criterion for log-Sobolev
already existed, thanks to  Bobkov and G\"{o}tze \cite{BG99}). 
See \cite{boucblm04mifi} for other developements.
Inequalities \eqref{1.1} above are part of a more general family
denoted $\Phi$-Sobolev inequalities.  A study of this family in connection with some
aspects of semi-group theory is done in \cite{Ch02}.

\medskip

The initial goal of this work is to obtain a precise dimension free
 isoperimetric inequality for $\mu_\alpha^{\otimes n}$ when $\alpha\in (1,2)$.
  Namely we
want to prove that  there exists a constant $C$
such that for all $n \in \dN$
\begin{equation}\label{1.2}
  (\mu_{\alpha}^{\otimes n})_s(\partial A)
  \geq
  C  \mu_{\alpha}^{\otimes n}(A)
  \Big(\log \big(
  \frac{1}{\mu_{\alpha}^{\otimes n}(A)}
  \big)\Big)^{1-\frac{1}{\alpha}},
\end{equation}
for all $A$ such that $\mu_{\alpha}^{\otimes n}(A) \leq \frac12$,
where $\mu_s(\partial A)$ denotes the surface measure of $A$ (see
section \ref{sec:isoperimetrie}). This  bound is known
for $\alpha =1$ \cite{BH97} and $\alpha=2$ \cite{BL96} and can be
deduced from \cite{bobk96ephs} in dimension 1. So  \eqref{1.2} is
exactly what is expected. This result is stronger than the
concentration result. Indeed it implies that for all $n$ and $A\subset \dR^{n}$
with $\mu_\alpha^{\otimes n}(A)=\mu_\alpha((-\infty,t])$ one has
$$\mu_\alpha^{\otimes n}\Big(A+hB_2^{n}\Big)\ge \mu_\alpha\Big((-\infty,t+h/K]\Big).$$
This interpolates between \eqref{isop-mu1} and \eqref{isop-mu2}.

Inequality \eqref{1.2} will be shown in Theorem \ref{th:iso} as the
achievement of a somewhat intricate story. Actually, we  prove much 
more and develop several useful methods on the way. They should find  a field of
applications in the study of empirical processes or in statistical physics.

Before describing the organisation of the paper, let us explain that
our proof relies on a method initiated by Ledoux \cite{ledo94sapi} and
improved in \cite{BL96}. It can be summarized as follows: any integrability
improving property of a semigroup with curvature bounded from below provides 
isoperimetric information for the invariant measure. Hence our problem 
translates to a question on the semigroup $(\PT{t}^{\alpha})_{t\ge 0}$ for 
$\alpha\in (1,2)$. However, a theorem of H{\o}egh-Krohn and Simon \cite{hoegs72hstd}
shows that $\PT{t}^{\alpha}$ is never continuous from $\mathbb L^2(\mu_\alpha)$
to  $\mathbb L^{2+\varepsilon}(\mu_\alpha)$. Since the $\dL_p$ scale is too rough for 
our problem, we analyze the regularizing properties in appropriate scales
of Orlicz spaces and ask wether the semigroup maps $\dL^{2}(\mu_\alpha)$ into a
smaller  Orlicz space. 

\smallskip
Section \ref{II} contains the required elements on Orlicz spaces. 

\smallskip
Section \ref{III} presents  a sufficient condition on the Young function
$\tau$ for $\mathbf{Q}_t^{\alpha}$ (a slightly modified
$\PT{t}^{\alpha}$) to map continuously $\dL^2$ into $\dL_{\tau}$, for a
fixed $t$.  This condition relies on the probabilistic
representation of $\PT{t}^{\alpha}$ (Girsanov transformation) and on
martingale methods inspired by \cite{KKR,Cat03}.
Unfortunately the method cannot reach  the contraction property (only
boundedness, simply called $\tau$-Orlicz hyperboundedness) and does
not easily yield explicit bounds. The criterion readily extends to
certain perturbations of an ultracontractive semi-group.


\smallskip
In order to get  contraction results and explicit bounds, we 
build in Section~\ref{IV}  the full analogue of Gross theory.
Following  \cite{Cat02} we start with the analogue of a result  by
H{\o}egh-Krohn and Simon (Theorem \ref{4.3}):
if $\PT{t_0}$ is continuous from $\dL^{2}(\mu)$ into $\dL_\tau(\mu)$ 
then $\mu$ satisfies a
defective logarithmic Orlicz-Sobolev inequality (DLOSI).  This
 is actually a particular $F$-Sobolev inequality as 
studied in \cite{Wa00,GWa02} (the notion apparently goes back to
Concordet). For the $|x|^{\alpha}$ Boltzmann measure $\mu_{\alpha}$,
it is equivalent (see \eqref{Rosen}) to the following result of
Rosen \cite{Ro76}: there exist $A$ and $B$ such that for $\int f^2
e^{-2|y|^{\alpha}} dy = 1$,
\begin{equation}\label{1.3}
  \int f^2(y) \Big(\log_+(|f(y)|)\Big)^{2(1-\frac{1}{\alpha})}
  e^{-2|y|^{\alpha}}  dy
  \leq
  A \int |\nabla f|^2 e^{-2|y|^{\alpha}}  dy  +  B  .
\end{equation}
See Adams \cite{Ad79} for  extensions and Zegarlinski \cite{Zeg01a}
 for an application of  Rosen type  inequalities to the study for Gibbs measures
with non-Gaussian tails.

Next we  consider   homogeneous $F$-Sobolev inequalities. 
 One of our main results is Theorem~\ref{th:go}
 where we obtain the equivalence between a $F$-Sobolev
inequality and the $\tau_q$ Orlicz-hyperboundedness (or
hypercontractivity) of the whole semigroup for $\tau_q(x):=x^p e^{qF(x^p)}$.
Under a few assumptions of $F$, 
the time evolution of the regularizing effect is quantified. 
A weak  form of part of these results
appeared in \cite[Theorem 1.2 and Theorem 2.4]{GWa02}.
 These authors proved that  a particular tight $F$-Sobolev
inequality is equivalent to Orlicz-hyperboundedness for some time.
Their motivation was  a criterion for the generator to have a non-empty 
 essential spectrum  (see \cite{GWa02,Wa00} for connections with 
 super-Poincar\'e inequalities).
By  Theorem~\ref{th:go}, a tight $F$-Sobolev inequality for
 a nonnegative $F$ garantees that  the semigroup is Orlicz-hypercontractive.
 We conclude this section by extending the well known inequality of Rothaus
\cite{Rot85}:
 under spectral gap assumption this allows to turn certain defective
$F$-Sobolev inequalities into tight ones.  

\smallskip 
Section~\ref{sec:fsobolev}  provides a thorough study of Sobolev type 
inequalities. In the Gaussian context the log-Sobolev inequality is
canonical and has plenty of remarkable properties: it tensorizes,
provides concentration via Herbst argument, hypercontractivity
and entropy decay along the semigroup. In our more general setting,
 in particular for $\mu_\alpha$, $\alpha \in (1,2)$, no such miracle
happens. Several Sobolev inequalities are available. However none of 
them concentrates all good properties. This is why we undertake a precise
study of Beckner type inequalities, of homogeneous $F$-Sobolev inequalities
and additive $\varphi$-Sobolev inequalities also called $\Phi$-Entropy inequalities
(wee shall not discuss the latter in terms of exponential decay of $\Phi$-entropy.
 See \cite{Ch02,Wa03}). Our strategy is to provide each inequality with 
a simpler reduced form relating the measure of sets to their $\mu$-capacity.
This notion was alluded to by the first and last-named authors in   \cite{BR03}.
Here we use it systematically in the spirit of Maz'ja \cite{mazy85ss}. Note that
the probabilistic setting is delicate since constant functions are equality 
cases in all our inequalities. Our approach is an extension to any dimension
of the criteria on the real line recently obtained through Hardy inequalities
\cite{BG99,BR03}. It provides new criteria and equivalences between
several Sobolev inequalities. A final figure summarizes the situation.

\smallskip

Section \ref{concentration} deals with the consequences of generalized
Beckner inequalities for the concentration of measure. They are immediate
from the method of Lata{\l}a-Oleszkiewicz, and where discussed independently 
 by Wang \cite{Wa03b}. Our contribution here comes from our sharp criteria 
for these inequalities. In particular we give general neat conditions for products
of measures on $\dR$ to enjoy  dimension free concentration with rate
$e^{-\Phi(t)}$ where $\Phi(t)$ is convex, but less than $t^{2}$.
Under reasonable assumptions the criterion is satisfied by the measure
$e^{-\Phi(t)}/Z$ itself, so the concentration is sharp.
For other results in connection with mass transportation,  see also
\cite{Wa03,CG03,GGM}.

\smallskip

Section~\ref{examples} illustrates all the previous results in the
case of $|x|^{\alpha}$ Boltzmann measures. In this concrete situation
we explain how to deal with the technical conditions involved. We
also develop a perturbation argument similar to the one of 
\cite[section 4]{Cat03}. 

\smallskip

The final section deduces isoperimetric inequalities from 
semigroup hyperboundedness properties. The claimed infinite
dimensional isoperimetric bound \eqref{1.2} is derived.
As a consequence a family of comparison theorems is provided.

\smallskip
For sake of clarity we decided not to develop our argument
in its full generality. However,
most  of our results easily extend to more general situations,
encompassing diffusion operators on Riemannian manifolds. 
This is the case of the Gross-Orlicz theory, of the reductions to 
inequalities between capacity and measure. The final
isoperimetric lower bounds would work  when
the curvatures of the generators is bounded from below.


\section{Orlicz hypercontractivity.}\label{II}

In this section we shall discuss a weakened form of
\emph{hypercontractivity} and \emph{hyperboundedness},
replacing $\dL^p$ spaces by \emph{Orlicz spaces}.
Before to state the natural definition we have in mind,
and because Orlicz spaces are somewhat intricate to use,
we shall first introduce in this section
the material we need. Some definitions are not the usual ones used
{\it e.g.} in the book by Rao and Ren (\cite{RR}).

In the sequel we shall consider a \emph{complementary pair}
$(\tau^*,\tau)$ of continuous and even Young's
functions ({\it i.e.} $\tau^*$ is the Fenchel-Legendre dual function of
$\tau$, both being convex functions vanishing at the origin) satisfying
\begin{equation}\label{2.2}
  \lim_{y \rightarrow +\infty} \frac{\tau(y)}{y^2} = +\infty
  \quad \textrm{ and } \quad
  \lim_{y \rightarrow +\infty} \frac{\tau(y)}{y^p} =  0 \quad
  \textrm{ for }  p>2 .
\end{equation}
It follows that
\begin{equation}\label{2.3}
  \lim_{y \rightarrow +\infty} \frac{\tau^*(y)}{y^2}  =  0
  \quad \textrm{ and } \quad
  \lim_{y \rightarrow +\infty} \frac{\tau^*(y)}{y^p}  =  +\infty
  \quad \textrm{ for }  p<2 .
\end{equation}
We shall also assume that $\tau$ and $\tau^*$ both satisfy the
\emph{$\Delta_2$ condition} ({\it i.e.} $\tau(2y) \leq K \, \tau(y)$
for some $K>1$ and $y\geq y_1 \geq 0$, and a similar result for $\tau^*$ with
possibly different $K^*$ and $y_1^*$).
It follows that they both satisfy the \emph{$\nabla_2$ condition}
too ({\it i.e.} $2l \, \tau(y) \leq \tau(ly)$ for some
$l>1$ and $y\geq \bar{y_1}\geq 0$ and similarly for
$\tau^*$ with $l^*$ and ${\bar y_1}^*$), see \cite[p. 23]{RR}.

We also assume that the pair $(\tau, \tau^*)$ is normalized, {\it i.e.}
\begin{equation}\label{2.4}
  \tau(0)  =  \tau^*(0)  =  0
  \quad \textrm{ and } \quad
  \tau(1)  +  \tau^*(1)  =  1.
\end{equation}

The space $\dL_{\tau}(\mu)$ is the space of measurable functions $f$ such that
\begin{equation}\label{2.5}
  I_{\tau}(f) \stackrel{\rm def}{=}  \int \tau(f)  d\mu  <  +\infty .
\end{equation}

Thanks to the $\Delta_2$ property, $\dL_{\tau}$ and $\dL_{\tau^*}$
are linear spaces.

We shall use two norms on each space,
\begin{eqnarray}\label{2.6}
  N_{\tau}(f) \stackrel{\rm def}{=}
  \inf \{k>0 ; I_{\tau}\big(\frac fk\big) \leq \tau(1) \}, \\
  \NRM{f}_{\tau} \stackrel{\rm def}{=}
  \sup \{\int |fg| d\mu  ;  N_{\tau^*}(g) \leq 1  \} , \nonumber
\end{eqnarray}
with similar definitions for $\tau^*$.
The first one is called the \emph{gauge} (or Luxemburg) \emph{norm} while
the second definition does not coincide with the usual
Orlicz norm (\cite[Definition 2 page 58]{RR}).

Indeed in order to compare Orlicz norms and usual
$\dL^p(\mu)$ norms in the framework of Markov
semi-groups (for a given Probability measure $\mu$),
we need the norms of constant functions to be
equal to the same constant. With our definitions,
and thanks to the regularity properties of $\tau$
and $\tau^*$ it is known that for $f\in \dL_{\tau}$ and
$g\in \dL_{\tau^*}$ it holds (see \cite[Proposition 1 p.58]{RR})
\begin{equation}\label{2.7}
  I_{\tau}\big(\frac{f}{N_{\tau}(f)}\big) = \tau(1)
  \quad \textrm{ and } \quad
  \int  |fg| d\mu  \leq  N_{\tau}(f)  N_{\tau^*}(g) .
  \end{equation}
  Accordingly
  \begin{equation}\label{2.8}
  \NRM{f}_{\tau} \leq  N_{\tau}(f)
  \quad \textrm{ and } \quad
  \NRM{\ind}_{\tau} = N_{\tau}(\ind) = 1,
\end{equation}
similar results being true for $\dL_{\tau^*}$.

Note that if we replace $\tau(1)$ by $1$ in the definition of
$N_{\tau}$ we get another gauge norm $N^1_{\tau}$ (the usual one)
which is actually equivalent to $N_{\tau}$, more precisely
$$
N_{\tau}^1(f) \leq N_{\tau}(f) \leq \frac{1}{\tau(1)} N_{\tau}^1(f).
$$
It follows thanks to $\Delta_2$ (see \cite[chapter IV]{RR}) that
$(\dL_{\tau},N_{\tau})$ is a reflexive Banach space with dual space
$(\dL_{\tau^*}, \NRM{\cdot}_{\tau^*})$. Also note that
$N_{\tau}$ and $\NRM{\cdot}_{\tau}$ are equivalent
(see (18) p.62 in \cite{RR}) and that
the subset of bounded functions is everywhere dense in $\dL_{\tau}$.
The same holds when we replace $\tau$ by $\tau^*$.

Finally remark that if $N_{\tau}(f) \geq 1$,
$$
  \tau(1)
  =
  I_{\tau}\Big(\frac{f}{N_{\tau}(f)}\Big)
  \leq
  \frac{1}{N_{\tau}(f)} I_{\tau}(f),
$$
so that
\begin{equation}\label{2.9}
  N_{\tau}(f) \leq \max \Big(1, \frac{I_{\tau}(f)}{\tau(1)}\Big).
\end{equation}
Conversely if $f(x) \geq N_{\tau}(f) y_1$
(recall the definition of $\Delta_2$) then
$$
  \tau(f(x))
  =
  \tau\Big(N_{\tau}(f) \frac{f(x)}{N_{\tau}(f)}\Big)
  \leq
  K^{\frac{\log(K)}{\log(2)} + 1} \tau\Big(\frac{f(x)}{N_{\tau}(f)}\Big).
$$
It follows that
\begin{equation}\label{2.10}
  I_{\tau}(f)
  \leq
  \tau( N_{\tau}(f) y_1) + K^{\frac{\log(K)}{\log(2)}+1}\tau(1).
\end{equation}

\begin{definition}[Orlicz-hyperboundedness]\label{OH}
 We shall say that a $\mu$-symmetric semi-group $(\PT{t})_{t \geq 0}$ is
 \emph{Orlicz-hyperbounded} if for some $t>0$, $\PT{t}$ is
 a continuous mapping from $\dL^2(\mu)$ into some Orlicz space
 $\dL_{\tau}(\mu)$ for some Young function $\tau$ satisfying
 $$
   \lim_{y \rightarrow +\infty} \frac{\tau(y)}{y^2} = +\infty
   \quad \textrm{ and } \quad
   \lim_{y \rightarrow +\infty} \frac{\tau(y)}{y^p} = 0
   \quad \textrm{ for }  p>2.
 $$
\end{definition}

We can now give the definition of Orlicz-hypercontractivity

\begin{definition}[Orlicz-hypercontractivity]\label{2.11}
 We shall say that a $\mu$-symmetric semi-group $(\PT{t})_{t \geq 0}$ is
 \emph{Orlicz-hypercontractive} if for some $t>0$, $\PT{t}$
 is a \emph{contraction} from $\dL^2(\mu)$ into
 $(\dL_{\tau}(\mu),N_{\tau})$ for some Young function $\tau$ as in Definition
 \ref{OH}.
 Equivalently $\PT{t}$ is a contraction from
 $(\dL_{\tau^*}(\mu),\NRM{\cdot}_{\tau^*})$ into $\dL^2(\mu)$.
\end{definition}

Note that this definition is coherent.
With our definitions and thanks to Jensen's inequality, $\PT{s}$
is for all $s>0$ a contraction in both $(\dL_{\tau}(\mu),N_{\tau})$ and
$(\dL_{\tau^*}(\mu),\NRM{\cdot}_{\tau^*})$
(as well as in both $(\dL_{\tau^*}(\mu),N_{\tau^*})$ and
$(\dL_{\tau}(\mu),\NRM{\cdot}_{\tau}$)
whose norm (equal to $1$) is attained for the constant functions.
In particular if the contraction property in Definition \ref{2.11}
holds for some $t$, it holds for all $s>t$.

The next section will give some criterion for semi-group to be
Orlicz-hyperbounded.


\section{Orlicz hyperboundedness for $|x|^{\alpha}$ and general
Boltzmann measures.}\label{III}

For $\alpha \in ]1,2[$ we may consider the $\C{2}$ function
$u_{\alpha}$ defined on $\dR$ by
\begin{equation}\label{3.1}
  u_{\alpha}(x) = \left\{
    \begin{array}{ll}
    |x|^{\alpha} & \textrm{ for } |x|>1 \\
    \frac{\alpha (\alpha-2)}{8}  x^4  +  \frac{\alpha(4-\alpha)}{4}  x^2
    + (1-\frac 34 \alpha + \frac 18 \alpha^2) & \textrm{ for } |x|\leq 1 .
    \end{array}
  \right.
\end{equation}
With this choice it is easy to see that $u_{\alpha}$
is convex and bounded below by $1-\frac 34 \alpha  + \frac 18 \alpha^2$
which is positive.

The associated \emph{Boltzmann measure} on $\dR^n$ is defined as
\begin{equation}\label{def:nu}
  \nu_{\alpha}(dx)
  = Z_{\alpha}^{-1} e^{- 2 \sum_{i=1}^n u_{\alpha}(x_i)} dx
  = Z_{\alpha}^{-1} e^{- 2 U_\alpha (x)} dx
\end{equation}
where $Z_{\alpha}$ is the proper normalizing constant such that
$\nu_{\alpha}$ is a Probability measure and
$U_\alpha (x) \stackrel{\rm def}{=} \sum_{i=1}^n u_\alpha(x_i)$.

To the Boltzmann measure is associated a symmetric semi-group
$(\PT{t}^{\alpha})_{t \geq 0}$ generated by the
operator
$$
  A_{\alpha} = \frac 12 \Delta - \nabla  U_{\alpha} \cdot \nabla .
$$
One can show (see \cite{Cat02} or \cite{Cat03}) that the semi-group is given by
\begin{equation}\label{3.2}
  \big(\PT{t}^{\alpha} h\big)(x)
  = e^{U_{\alpha}(x)} \dE^{\dP_x} [h(X_t) e^{ - U_{\alpha}(X_t)} M_t],
\end{equation}
where $\dP_x$ is the Wiener measure such that $\dP_x\big(X_0=x\big)=1$
({\it i.e.} under $\dP_x$, $X_{.}$ is a $n$-dimensional Brownian motion
starting from $x$) and $M_t$ is defined as
\begin{equation}\label{3.3}
  M_t = \exp \Big(
  \frac 12 \int_0^t (\Delta U_{\alpha} - |\nabla U_{\alpha}|^2)(X_s) ds \Big).
\end{equation}
(see \cite[section 7]{Cat02} for a proof). We have chosen this form of $A_{\alpha}$ (with 1/2 in
front of $\Delta$) in this section not to introduce extra variance on the Brownian motion.

Since $e^{U_{\alpha}}$ belongs to all $\dL^p(\nu_{\alpha})$ for
$p<2$, an almost necessary condition for Orlicz hypercontractivity
is that $\PT{t}^{\alpha}(e^{U_{\alpha}})$ belongs to some
$\dL_{\tau}(\nu_{\alpha})$.

In \cite[section 3]{Cat03}, the ``Well Method'' originally due to Kavian, Kerkyacharian and Roynette
\cite{KKR} is further developed and allows to get estimates for 
$\PT{t}^{\alpha}(e^{U_{\alpha}})$,
but for $\alpha \geq 2$.

We shall below extend the ``Well Method'' to the case $1<\alpha<2$.

\begin{theorem}\label{3.4}
 Let $\tau$ be a Young function satisfying
 $\tau(y)=y^2 \, \psi(y)$ for some positive and non
 decreasing function $\psi$ going to $+ \infty$ at infinity.
 We also assume that there exists a constant $k(\psi)$  such that
 $\psi(2 y) \leq k(\psi) \psi(y)$ ({\it i.e.} $\tau$ satisfies
 $\Delta_2$). Let $\nu_\alpha$ be the Boltzmann measure
 defined on $\dR^n$ in \eqref{def:nu}.

 Then for any integer $n$, $\PT{t}^{\alpha}(e^{U_{\alpha}})$ belongs to
 $\dL_{\tau}(\nu_{\alpha})$ if there exists some constant
 $d_{\alpha}<\alpha^2$ such that
 $$
   \int^{+\infty} \psi(e^{|x|^{\alpha}})
   e^{- d_{\alpha} t|x|^{2(\alpha - 1)}} dx < + \infty .
 $$
\end{theorem}

The proof below can be used (or improved) to get explicit bounds.

\begin{proof}
 First remark that $U_{\alpha}$ satisfies
 \begin{equation}\label{3.5}
   \frac 12 \Big(|\nabla U_{\alpha}|^2(x) - \Delta U_{\alpha}(x)\Big)
   \geq
   G_{\alpha}(U_{\alpha}(x)) - c_{\alpha} = H_{\alpha}(U_{\alpha}(x)),
 \end{equation}
 with
 $G_{\alpha}(y)= \frac{\alpha^2}{2} |y|^{2 (1-\frac{1}{\alpha})}$ and
 $c_{\alpha}=n(1 + \frac 12 \alpha (\alpha - 1))$,
 with our choice of $u_{\alpha}$ for $|x|\leq 1$.
 Note that $H_{\alpha}$ admits an inverse $H_{\alpha}^{-1}$ defined on
 $[-c_{\alpha};+\infty)$ with values in $\dR^+$.

 For $0<\varepsilon $ define the stopping time $T_x$ as
 \begin{equation}\label{3.6}
   T_x =
   \inf\{ s>0 ; \frac 12 (|\nabla  U_{\alpha}|^2 - \Delta U_{\alpha})(X_s)
   \leq H_{\alpha} \big( U_{\alpha}(x) - \varepsilon \big)\}.
  \end{equation}

 Note that for all $x \in \dR^n$, $T_x > 0$ $\dP_x \textrm{ a.s.}$
 provided $U_{\alpha}(x) - \varepsilon \geq 0$ and that on $T_x < +\infty$,
 \begin{equation}\label{3.7}
   U_{\alpha}(X_{T_x})
   \leq
   H_{\alpha}^{-1} \Big(\frac 12 (|\nabla U_{\alpha}|^2
   - \Delta U_{\alpha})(X_{T_x})\Big)
   \leq
   U_{\alpha}(x) - \varepsilon .
 \end{equation}

 Introducing the previous stopping time we get
 \begin{equation*}
   \dE^{\dP_x}[M_t]
   =
   \dE^{\dP_x}[M_t \ind_{t<T_x}]+\dE^{\dP_x}[M_t \ind_{T_x \leq t}]
   =
   A + B,
 \end{equation*}
 with
 \begin{equation}\label{3.8}
   A  = \dE^{\dP_x}[M_t \ind_{t<T_x}]
   \leq
   \exp \big( - t H_{\alpha}(U_{\alpha}(x) - \varepsilon) \big),
 \end{equation}
 and
 \begin{eqnarray}\label{3.9}
   B
   & = &
   \dE^{\dP_x}[M_t \ind_{T_x \leq t}]
   \nonumber\\
   & \leq &
   e^{c_{\alpha} t} \dE^{\dP_x} \big[ \exp
   \Big(\int_0^t \big( \frac 12 (\Delta U_{\alpha} - |\nabla U_{\alpha}|^2) -
   c_{\alpha} \big)(X_s) ds \Big) \ind_{T_x \leq t} \big]
   \nonumber \\
   & \leq &
   e^{c_{\alpha} t} \dE^{\dP_x} \big[ \exp
   \Big(\int_0^{T_x} \big( \frac 12 (\Delta U_{\alpha} - |\nabla U_{\alpha}|^2)
   - c_{\alpha} \big)(X_s) ds \Big) \ind_{T_x \leq t} \big]
   \nonumber \\
   & \leq &
   e^{c_{\alpha}  t} \dE^{\dP_x}\big[ \exp
   \Big(\int_0^{T_x} \big( \frac 12 (\Delta U_{\alpha} - |\nabla U_{\alpha}|^2)
   \big)(X_s) ds \Big) \ind_{T_x \leq t} \big]
   \nonumber\\
   & = & e^{c_{\alpha} t} \dE^{\dP_x}[M_{T_x} \ind_{T_x \leq t}].
 \end{eqnarray}
 But $e^{-U_{\alpha}(X_s)} M_s$ is a  bounded  $\dP_x$ martingale.
 Hence, according to Doob optional sampling Theorem
 \begin{equation}\label{3.10}
   \dE^{\dP_x}[e^{-U_{\alpha}(X_{T_x})} M_{T_x} \ind_{T_x \leq t}]
   \leq
   \dE^{\dP_x} \big[ e^{-U_{\alpha}(X_{t\wedge T_x})} M_{t\wedge T_x} \big]
   = e^{-U_{\alpha}(x)}.
 \end{equation}
 According to \eqref{3.7}, $\displaystyle
 e^{- U_\alpha(X_{T_x})} \geq e^{\varepsilon} e^{- U_{\alpha}(x)}$,
 so that thanks to \eqref{3.10},
 $$
   \dE^{\dP_x}[M_{T_x} \ind_{T_x \leq t}] \leq e^{- \varepsilon}.
 $$

 Using this estimate in \eqref{3.9} and using \eqref{3.8} we finally obtain
 \begin{equation}\label{3.11}
   \dE^{\dP_x}[M_t] \leq e^{- t H_{\alpha}(U_{\alpha}(x) -\varepsilon)}
   + e^{- \varepsilon} e^{c_{\alpha} t}.
 \end{equation}

 It remains to choose $\varepsilon = \beta U_{\alpha}(x)$ for some $\beta < 1$.
 Inequality \eqref{3.11} yields for $|x|>1$
 \begin{equation}\label{3.12}
   \dE^{\dP_x}[M_t]
   \leq
   e^{c_{\alpha} t} \Big(e^{-\big( t (1-\beta)^{2(1-\frac{1}{\alpha})}
   \frac 12 \alpha^2 |x|^{2(\alpha - 1)} \big)} +
   e^{- \beta |x|^{\alpha}}\Big),
 \end{equation}
 while a rough bound for $|x|\leq 1$ is
 $\dE^{\dP_x}[M_t] \leq e^{c_{\alpha} t}$,
 since $M_t \leq e^{c_{\alpha} t}$ according to \eqref{3.5}.

 Finally recall that
 $\PT{t}^{\alpha}(e^{U_{\alpha}})=e^{U_{\alpha}} \dE^{\dP_x}[M_t]$
 and remark that since $\alpha < 2$ the dominating term in \eqref{3.12}
 is the first one in the sum, at least for large $|x|$.
 Together with the property of $\psi$ all this yields the statement
 in the Theorem.
 \qed
\end{proof}

As in \cite[Theorem 2.8.]{Cat03} we shall see below that the condition
$\PT{t}^{\alpha}(e^{U_{\alpha}}) \in \dL_{\tau}(\nu_{\alpha})$ 
is also a sufficient condition for
$\tau$-Orlicz hyperboundedness.

\begin{theorem}\label{3.14}
Let $\tau$ be as in Theorem \ref{3.4}.
A sufficient condition for $(\PT{t}^{\alpha})_{t \geq 0}$ to be $\tau$-Orlicz
hyperbounded is that
$\PT{t}^{\alpha}(e^{U_{\alpha}}) \in \dL_{\tau}(\nu_{\alpha})$
for some $t>0$.
\end{theorem}

\begin{proof}
First recall that thanks to \eqref{3.5},
$M_t \leq e^{c_{\alpha} t}$. On the other hand, the Brownian
semi-group $(\PT{s})_{s \geq 0}$ on $\dR^n$ is ultracontractive and
$\NRM{\PT{s}}_{\dL^2(dx) \rightarrow \dL^{\infty}(dx)} =(4\pi s)^{-\frac n4}$.

Pick some smooth function $f$ on $\dR^n$ with compact support.
Since $|f|e^{-U_{\alpha}}\in \dL^2(dx)$ and using the
Markov property, for $s>0$ and $t>0$, it holds
\begin{eqnarray*}
\dE^{\dP_x}[M_{t+s} \big(e^{-U_{\alpha}} |f|\big)(X_{t+s})]
& = &
\dE^{\dP_x}[M_t \dE^{\dP_{X_t}}[M_s \big(e^{-U_{\alpha}} |f|\big)(X'_s)]] \\
& \leq &
e^{c_{\alpha} s} \dE^{\dP_x}[M_t \big(P_s(|f| e^{-U_{\alpha}})\big)(X_t)] \\
& \leq &
e^{c_{\alpha} s} (4\pi s)^{-\frac n4}
\NRM{f}_{\dL^2(\nu_{\alpha})} \dE^{\dP_x}[M_t] .
\end{eqnarray*}
Hence
\begin{eqnarray*}
&&
\qquad
\int \tau(P_{t+s}^{\alpha}(|f|)) d\nu_{\alpha}
=
\int \tau\Big(e^{U_{\alpha}} \dE^{\dP_x}\big[
M_{t+s} \big(e^{-U_{\alpha}}|f|\big)(X_{t+s}) \big] \Big) d\nu_{\alpha} \\
&&
\leq  C^2(\alpha, s) \NRM{f}_{\dL^2(\nu_{\alpha})}^2
k(\psi)^{1+\log_2\big(C(\alpha,s)\NRM{f}_{\dL^2(\nu_{\alpha})}\big)}
\! \int \! \tau \big( e^{U_{\alpha}} \dE^{\dP_x}[M_t]\big) d\nu_{\alpha},
\end{eqnarray*}
with $C(\alpha,s)=e^{c_\alpha s}(4\pi s)^{-\frac n4}$.
In particular, if $\NRM{f}_{\dL^2(\nu_{\alpha})}=1$,
$$
\int \tau(\PT{t+s}^{\alpha}(|f|)) d\nu_{\alpha} \leq K(t,s,\alpha).
$$
According to \eqref{2.9}, $\PT{t+s}^{\alpha}$ is thus continuous.
\qed
\end{proof}

\begin{example}\label{3.15}
The best possible choice of $\psi$ in Theorem \ref{3.4} is given by
$$
\psi(y) = \exp \Big((\log(|y|)^{2 (1 - \frac{1}{\alpha})} \Big)
\quad \textrm{ for } |y| \textrm{ large enough.}
$$
According to Theorem \ref{3.4} and \ref{3.14},
$(\PT{t}^{\alpha})_{t \geq 0}$ is then $\tau$-Orlicz hyperbounded for
$t>\frac{1}{\alpha^2}$.
\end{example}

The previous scheme of proof, without any change,
obviously extends to the general framework we have
introduced. Let us describe the situation.

Let $\PT{t}$ be a $\mu$-symmetric semi-group on a space $E$ as
describe in section \ref{I}, with generator $\GI$. For $V$ in the
domain $\dD(\GI)$ of $\GI$, we introduce the general Boltzmann measure
$d\nu_V = e^{-2V} d\mu$ and assume that $\nu_V$ is a probability
measure. Under some assumptions it is known that one can build a
$\nu_V$-symmetric semi-group $(\PT{t}^V)_{t \geq 0}$, via
\begin{equation}\label{GBSG}
(\PT{t}^V h)(x) = e^V(x) \dE^{\dP_x}\big[h(X_t) e^{-V(X_t)} M_t\big],
\end{equation}
with
$$
M_t= \exp \Big(\int_0^t \big(\GI V(X_s) - \Ga(V,V)(X_s)\big) ds
\Big).
$$
In the general case these assumptions are denoted by (H.F) in
\cite{Cat03}. Here we have chosen the usual definition
$$\Ga(V,V)=\frac 12 \, \left(\GI V^2 - 2 V \GI V\right) .$$

When $E=\dR^n$, $L=1/2 \Delta$ and $\mu=dx$ each of the following
conditions (among others) is
sufficient for \eqref{GBSG} to hold:\\
$(i)$ there exists some $\psi$ such that $\psi(x) \rightarrow +\infty$
as $|x| \rightarrow +\infty$
and $\nabla V \cdot \nabla \psi - \Delta \psi$ is bounded from below, \\
$(ii)$ $\int |\nabla V|^2 d\nu_V < +\infty$.

For the first one see {\it e.g.} \cite[p.26]{Ro99} and
\cite[(5.1)]{Cat03}, for the second one see {\it e.g.} \cite{CL96}.

Introduce the analogue of \eqref{3.5}:

\begin{assumption}[OB(V)]\label{OBF}
We shall say that $V$ satisfies assumption OB(V), if\\
$(i)$
$V$ is bounded from below by some (possibly negative) constant $d$.\\
$(ii)$
There exist some $c \in \dR$, $u_0>0$ and a function
$G: \dR^+ \to \dR^+$ such that $G(u) \to +\infty$ as $u \to +\infty$
and $G(u)/(u+1)$ is bounded for $u\geq u_0$, and such that for all
$x\in E$,
$$
 \Ga(V,V)(x) - \GI V(x) \geq  G(|V(x)|) - c .
$$
\end{assumption}

Assumption OB(V) ensures that the dominating term in
the analogue of \eqref{3.12} is the former for $x$ large enough.
Then

\begin{theorem}\label{OHGB}
Let $\tau$ be as in Theorem \ref{3.4}.
If $(\PT{t})_{t \geq 0}$ is \emph{ultracontractive} and $V$ satisfies
assumption OB(V),
then the perturbed semi-group $(\PT{t}^V)_{t \geq 0}$ is
$\tau$-Orlicz hyperbounded as soon as for some $C>0$
$$
\int \psi(e^{V}) e^{- C G(|V|)} d\mu < + \infty .
$$
\end{theorem}

\begin{remark}
   The assumption OB(V) appeared first in Rosen's work with
   $G(u)=u^{2(1-\frac{1}{\alpha})}$ (\cite[condition (5) in Theorem
   1]{Ro76}) for $|x|^{\alpha}$ Boltzmann measures on $\dR^n$, for
   which a modified version of (defective) logarithmic Sobolev
   inequalities is obtained. Though Rosen proved that this condition
   is in a sense optimal (see his Theorem 5) for his log-Sobolev like
   inequality for the $|x|^{\alpha}$ Boltzmann measure, he did not
   rely this inequality to the Orlicz hyperboundedness of the
   associated semi-group. Furthermore, we think that the meaning of
   assumption OB(V) is enlightened by our probabilistic approach. We
   shall discuss later Rosen's results in relationship with
   $F$-Sobolev inequalities (see section \ref{examples}).
\end{remark}

\begin{remark}\label{Kunz}
   If $E=\dR^n$, $L=1/2 \Delta$ and $\mu=dx$, the second condition in
   assumption OB(V) implies the existence of a spectral gap (see {\it
     e.g.} \cite[Proposition 5.3.(2)]{Cat03}), provided $V$ goes to
   infinity at infinity.
\end{remark}


\section{Gross theory for Orlicz hypercontractivity.}\label{IV}

In this section we assume that $\tau$ and $\tau^*$ are smooth and
increasing on $\dR^+$ (hence one to one on $\dR^+$).
We also simply denote by $\NRM{f}_p$ the $\dL^p(\mu)$ norm of $f$ (when
no confusion on the underlying measure $\mu$ is possible).
In this section we shall assume for simplicity that $\mu$
\emph{is a probability measure}. The framework is the one
described in the introduction (see the notations therein).

\subsection{An Orlicz version of H{\o}egh-Krohn and Simon Theorem.}\label{IV1}

Since we do not {\it a priori} consider a parametrized family of
Orlicz functions, contrary to the family $(\dL^p ; p\geq 2)$ used in
Gross theory, the extension of this theory to our framework is not
immediate. Our definitions are nevertheless similar to the ones used
in H{\o}egh-Krohn and Simon result relating hypercontractivity and
logarithmic Sobolev inequalities. A proof of H{\o}egh-Krohn and Simon
Theorem using semi-group techniques is contained in \cite[Theorem
3.6]{Bak94}. Another proof is given in \cite[Corollary 2.8]{Cat02}.

We  follow the route in \cite{Cat02} in order to get some
functional inequality for an Orlicz  hyperbounded
semi-group. The starting point is the following particular case of
Inequality (2.4) in \cite{Cat02}:
for all non-negative $f\in \dD$ (a nice core algebra see \cite{Cat02},
in the usual $\dR^n$ case we may choose the smooth compactly supported
functions plus constants) such that $\int f^2 d\mu = 1$, all positive
and bounded $h$ and all $t\ge 0$,
\begin{equation}\label{4.2}
\int f^2 \log h d\mu
\leq
\frac t2 \mathcal E(f,f) + \log \int f h \PT{t} f d\mu .
\end{equation}
Recall that for \eqref{4.2} to hold, $\PT{s}$ has to be
$\mu$-symmetric.
If $\PT{t}$ maps continuously $\dL^2$ in some $(\dL_{\tau},N_{\tau})$
with operator norm $C(t,\tau)$, applying H\"{o}lder inequality we
obtain
$$
\int f^2 \log h d\mu \leq \frac t2 \mathcal E(f,f) + \log
(C(t,\tau)) + \log (\NRM{|f| h}_{\tau^*}).
$$
Hence if we can choose some $h$ such that the last term in the
above sum is less than $0$, we will obtain some functional inequality
reminding the (defective) logarithmic Sobolev inequality.  A natural
choice is
$$
h(f) = \frac{(\tau^*)^{-1}(f^2 \tau^*(1))}{|f|},
$$
since in this case $I_{\tau^*}(|f|h) = \int \tau^*(1) f^2 d\mu =
\tau^*(1)$.  It follows that $N_{\tau^*}(|f|h)=1$ and we may apply
\eqref{2.8}.  This choice is allowed provided $f^2 \log h$ is $\mu$
integrable and interesting provided $\log h$ is a non-negative
function growing to infinity with $f$. Note that with our choices of
$(\tau,\tau^*)$, $(\tau^*)^{-1}(y) \gg \sqrt{y}$ and thus $h(y)
\rightarrow +\infty$ as $y \rightarrow +\infty$.  We have shown

\begin{theorem}\label{4.3}
If the $\mu$-symmetric semi-group $(\PT{t})_{t \geq 0}$ is
$\tau$-Orlicz hyperbounded (with operator norm $C(t_0,\tau)$ for some
$t_0>0$) then for all $f \in \mathbb D$ the following (defective) logarithmic
Orlicz Sobolev inequality holds
\begin{equation*}
(DLOSI) \qquad
\ent_{\tau}(f) \leq  a \mathcal E(f,f)  +  b \NRM{f}_2^2 ,
\end{equation*}
with $\ent_{\tau}(f) \stackrel{\rm def}{=} \int f^2 \log \Big(
\frac{(\tau^*)^{-1}(\tau^*(1) (f/\NRM{f}_2)^2)}{|f/\NRM{f}_2|} \Big) d\mu$,
$a=\frac {t_0}{2}$ and $b= \log (C(t_0,\tau))$ provided the function
$$
y \mapsto y^2 \log \Big( \frac{(\tau^*)^{-1}(y^2 \tau^*(1))}{|y|} \Big)
$$
can be continuously extended up to the origin
(here $(\tau^*)^{-1}$ is the inverse function of $\tau^*$ and not $1/\tau^*$).

In particular if $(\PT{t})_{t \geq 0}$ is $\tau$-Orlicz hypercontractive,
(DLOSI) is tight {\it i.e.} becomes
$$
(TLOSI) \qquad
\ent_{\tau}(f) \leq a \mathcal E(f,f) .
$$
\end{theorem}

\begin{remark}
If we formally replace $\tau(y)$ by $y^p$ for some $p>2$, then $(\tau^*)^{-1}(y)$ behaves like
$y^{1/q}$ for the conjugate $q$ of $p$. Hence we recover the usual logarithmic Sobolev inequality as
in H{\o}egh-Krohn and Simon theorem. This is not surprising since the previous proof is the exact
analogue of the one in \cite{Cat02}.
\end{remark}

\begin{remark}\label{4.4}
Remark that our choice of $h$ is such that $N_{\tau^*}(|f|h) \leq 1$.
Hence we may replace the operator norm of $\PT{t}$ as a linear
operator between $\dL^2$ and $(\dL_{\tau},N_{\tau})$ by the similar
operator norm with $(\dL_{\tau},\NRM{\cdot}_{\tau})$ instead.
This should be interesting if $\PT{t}$ becomes a contraction for this norm,
while it is not for the previous one (recall \eqref{2.8}).
\end{remark}

In view of Theorem \ref{4.3} it is now natural to ask for a converse,
hence a Gross-Orlicz Theorem.

Actually an inequality like (DLOSI) was already discussed in the literature,
where it appears as a particular $F$-Sobolev inequality (see below).
In addition the explicit form of $\ent_{\tau}$ is not
easily tractable as it stands.
For instance we cannot obtain an explicit form of $\ent_{\tau}$ for
the $|x|^{\alpha}$ Boltzmann measure, but only an asymptotic behavior,
{\it i.e.} (recall Example \ref{3.15} and see
Example \ref{examplealpha} in this section)
\begin{equation} \label{equivalpha}
\log \left(\frac{(\tau^*)^{-1}(y^2 \tau^*(1))}{|y|}\right)
\approx
\left(\log (|y|)\right)^{2(1-\frac{1}{\alpha})},
\end{equation}
as $y \to +\infty$ (where $a \approx b$ means that $c a \leq b \leq C a$ for
some universal constant $c$ and $C$).
It is thus natural to ask whether one can replace one by the other in Theorem
\ref{4.3}. All these reasons lead to the study of
Orlicz-hyperboundedness in connection with general $F$-Sobolev inequalities.

\subsection{A Gross-Orlicz Theorem.}\label{IV2}

Our main Theorem is Theorem \ref{th:go} below.
This Theorem gives the equivalence between the
homogeneous $F$-Sobolev inequality and the
Orlicz hypercontractivity and gives a generalization of
the standard Gross Theorem \cite{gross75}.

Recall that the probability measure $\mu$ satisfies a
log-Sobolev inequality if there exists a constant $C_{LS}$
such that for any smooth enough function $f$,
\begin{equation}\label{in:lsob}
\int f^2 \log \left( \frac{f^2}{\mu(f^2)} \right) d\mu
\leq
C_{LS} \int  |\nabla f|^2 d\mu,
\end{equation}
where $\mu(f^2)$ is a short hand notation for
$\int f^2 d\mu$ and $|\nabla f|^2$ stands for $\Ga(f,f)$.
The following theorem is the celebrated Gross Theorem (\cite{gross75},
see also \cite{Ane}) relating this property to the
hypercontractivity of the semi-group $(\PT{t})_{t \geq 0}$.

\begin{theorem}[\cite{gross75}] \label{th:gross}
Let $\mu$ be a probability measure. The following holds:\\
$(i)$ Assume that $\mu$ satisfies a log-Sobolev inequality \eqref{in:lsob}
with constant $C_{LS}$, then, for any function $f$, any $q(0)>1$,
$$
\NRM{\PT{t} f}_{q(t)} \leq \NRM{f}_{q(0)} ,
$$
where $q(t)=1+(q(0)-1)e^{4t/C_{LS}}$.\\
$(ii)$ Assume that for any function $f$,
$$
\NRM{\PT{t}f}_{q(t)} \leq \NRM{f}_{2}
$$
with $q(t)=1+e^{4t/c}$ for some $c>0$. Then the probability measure
$\mu$ satisfies a log-Sobolev inequality \eqref{in:lsob} with constant $c$.
\end{theorem}

A natural extension of the log-Sobolev inequality is the homogeneous
$F$-Sobolev inequality. Let $F: \dR^+ \rightarrow \dR$ be a non-decreasing
function satisfying $F(1)=0$. A probability measure
$\mu$ satisfies an homogeneous $F$-Sobolev inequality if there exist
two constants $C_F$ and $\widetilde C_F$ such that for any smooth
enough function $f$,
\begin{equation} \label{in:Fsob}
  \int f^2 F \left( \frac{f^2}{\mu(f^2)} \right) d \mu
  \leq
  C_F \int |\nabla f|^2 d\mu + \widetilde C_F \int f^2 d \mu .
\end{equation}
If $\widetilde C_F=0$ (resp. $\neq 0$) the inequality is \emph{tight} (resp. \emph{defective}). We
shall use this terminology only when it is necessary.

We have the following result


\begin{theorem}[Gross-Orlicz] \label{th:go}
Fix $p >1$. Let $F: \dR^+ \rightarrow \dR$ be a $\C{2}$ non-decreasing
function satisfying $F(1)=0$.
Define for all $q \geq 0$, $\tau_q(x):=x^p e^{qF(x^p)}$.

$(i)$ Assume that
\begin{itemize}
\item there exists a non negative function $k$ on $\dR^+$
such that for all $q \geq 0$:
$\tau_q'' \tau_q \geq \frac{k(q)}{4} {\tau_q'}^2$ (hence
$\tau_q$ is a Young function),
\item there exists a non
negative function $\ell$ on $\dR^+$ and a constant $m \geq 0$
such that $\tau_q(x) F(x^p) \leq \ell (q) \tau_q(x) F(\tau_q(x)) + m$,
for all $q \geq 0$ and all $x\geq 0$,
\item the measure $\mu$ satisfies the homogeneous
$F$-Sobolev inequality \eqref{in:Fsob} with constants $C_F$ and
$\widetilde C_F$.
\end{itemize}
Then, for all non-decreasing $\C{1}$ functions
$q:\dR^+ \rightarrow \dR^+$ with $q(0)=0$ and
satisfying $q' \leq \frac{k(q)}{\ell (q)C_F}$,
the following holds for all $f$,
$$
N_{\tau_{q(t)}}(\PT{t} f)
\leq
e^{\frac1p [ mq(t) + \widetilde C_F \int_0^{q(t)} \ell(u)du ]}  \NRM{f}_p.
$$

$(ii)$ Conversely assume that there exist two non-decreasing functions,
$q, r :\dR^+ \rightarrow \dR^+$, differentiable at $0$,
with $q(0)=0$, such that for any $f$,
\begin{equation}\label{eq:hyp}
N_{\tau_{q(t)}}(\PT{t} f) \leq e^{r(t)} \NRM{f}_p .
\end{equation}
Then $\mu$ satisfies the following homogeneous
$F$-Sobolev inequality: for all $f$ smooth enough
$$
\int f^2 F\left( \frac{f^2}{\mu(f^2)} \right) d \mu
\leq
\frac{4(p-1)}{p q'(0)} \int |\nabla f|^2 d\mu
+
\frac{pr'(0)e^{r(0)}}{q'(0)}  \int f^2 d \mu .
$$
\end{theorem}

\begin{remark}
Note that by our assumptions on $\tau_q$,
$N_{\tau_{q(t)}}(f)$ is well defined.

Furthermore when $m=0$ the previous result states that the Orlicz hypercontractivity is equivalent
to the tight homogeneous $F$-Sobolev inequality ($\widetilde C_F=0$).
\end{remark}

\begin{proof}
We follow the general line of the original proof by Gross \cite{gross75},
see also \cite{Ane}. It is based on differentiation.

Without loss of generality we can assume that $f$ is non negative.
Then, for a general $\C{1}$ non
decreasing function $q : \dR^+ \rightarrow \dR^+$
satisfying $q(0)=0$, let $N(t):= N_{\tau_{q(t)}}(\PT{t} f)$.
For simplicity, we set $T(x,p):=\tau_p(x)$. Then, by definition of the
gauge norm \eqref{2.6} we have
$$
\int T \left(\frac{\PT{t} f}{N(t)},q(t) \right) d\mu = 1
\qquad \forall t \geq 0.
$$
Thus, by differentiation, we get
\begin{eqnarray*}
\frac{N'(t)}{N^2(t)}
\! \int \! \PT{t} f \partial_1 T \left( \! \frac{\PT{t} f}{N(t)},q(t) \right) d\mu
& = &
\frac{1}{N(t)}
\!\int\! \GI \PT{t}f \partial_1 T\left(\!\frac{\PT{t} f}{N(t)},q(t) \right) d\mu \\
&&
\quad + q'(t) \int \partial_2 T \left( \frac{\PT{t} f}{N(t)},q(t) \right) d\mu,
\end{eqnarray*}
or equivalently, if $g:=\frac{\PT{t} f}{N(t)}$,
\begin{equation}\label{eq:start}
\frac{N'}{N} \int g \; \partial_1 T(g,q) d\mu =
\int \GI g \; \partial_1 T(g,q) d\mu + q' \int  \partial_2 T(g,q) d\mu  .
\end{equation}
Here $\partial_1$ and $\partial_2$ are short hand notations for the partial derivative
with respect to the first and second variable respectively.

\bigskip

Let us start with the proof of the second part $(ii)$ of the Theorem.
For simplicity, assume that $N(0)=\NRM{f}_p=1$.
Take $t=0$ in the latter equality gives
$$
p N'(0) \int f^p d\mu =
p \int \GI f \cdot f^{p-1} d\mu + q'(0) \int  f^p F(f^p) d\mu ,
$$
because $g(0)=f$, $\partial_1 T(f,q(0))=p f^{p-1}$ and
$\partial_2 T(f,q(0))=f^p  F(f^p)$ (recall that $q(0)=0$ and $N(0)=1$).
Using the integration by parts formula
$\int \GI f \cdot \phi(f) d\mu = - \int |\nabla f|^2 \phi'(f) d\mu$, we get
$$
p \! \int\! \GI f \cdot f^{p-1} d\mu
=
- p(p-1)\! \int \!|\nabla f|^2 f^{p-2} d\mu
=
\frac{4(p-1)}{p} \!\int\! |\nabla f^{p/2}|^2 d\mu .
$$
Now, it follows from the bound \eqref{eq:hyp} that
$N'(0) \leq r'(0)e^{r(0)} \NRM{f}_p$. This implies
$$
pr'(0)e^{r(0)} \NRM{f}_p \!\int\!\! f^p d\mu
\geq
- \frac{4(p-1)}{p} \!\int\! |\nabla f^{p/2}|^2 d\mu +
q'(0) \!\int\!\!  f^p F(f^p) d\mu .
$$
Since  $\NRM{f}_p=1$, this achieves the proof of $(ii)$.

\bigskip

The proof of part $(i)$ is more technical.
A simple computation gives
$x \partial_1 T(x,q) = pT(x,q) + pq x^{2p} F'(x^p) e^{qF(x^p)}$.
Since  $F$ is non decreasing and $g \geq 0$, we get when $N'(t)\ge 0$
$$
\frac{N'}{N} \int g \partial_1 T(g,q) d\mu
\geq
\frac{p N'}{N} \int T(g,q) d\mu
=
\frac{p N'}{N} .
$$
On the other hand, using once again the integration by parts formula
$\int \GI f \cdot \phi(f) d\mu = - \int |\nabla f|^2 \phi'(f) d\mu$,
and our assumption on $\tau_q$,
\begin{eqnarray*}
\int \GI g  \partial_1 T(g,q) d\mu & = &
- \int |\nabla g|^2 \partial_{11} T(g,q) d\mu \\
& \leq &
- k(q) \int  |\nabla g|^2 \frac{\partial_1 T(g,q)^2}{4 T(g,q)} d\mu \\
& =  &
- k(q)  \int |\nabla \sqrt{T(g,q)}|^2 d\mu .
\end{eqnarray*}
Next,
$\partial_2 T(x,q) = T(x,q) F(x^p) \leq \ell(q) T(x,q) F(T(x,q)) + m$
by hypothesis.
Thus, \eqref{eq:start} becomes
$$
\frac{p N'}{N}
\leq
- k(q)   \!\!\int  \!\! |\nabla \sqrt{T(g,q)}|^2 d\mu
+
\ell(q) q'   \!\!\int  \!\!  T(g,q) F(T(g,q))d\mu + m q' .
$$
Note that the right hand side of this inequality contains the
three terms appearing in the homogeneous
$F$-Sobolev inequality \eqref{in:Fsob} applied to $\sqrt{T(g,q)}$
(recall that $\int (\sqrt{T(g,q)})^2 d\mu =1$).
In consequence, applying the homogeneous $F$-Sobolev inequality
\eqref{in:Fsob} to $\sqrt{T(g,q)}$ gives
$$
\frac{p N'}{N}
\leq q'(m + \widetilde C_F \ell(q))
+  [-k(q)+ q' \ell(q) C_F] \int |\nabla \sqrt{T(g,q)}|^2 d\mu .
$$
If $q' \leq \frac{k(q)}{C_F \ell(q)}$, it follows that
$\frac{p N'}{N} \leq  q'(m + \widetilde C_F \ell(q))$. This we proved when
$N'(t)\ge 0$. It is obviously true when $N'(t)<0$. Thus by integration
$$
N(t) \leq N(0)
e^{\frac1p [ mq(t) + \widetilde C_F \int_0^{q(t)} \ell(u)du ]} .
$$
Noting that $N(0)= \NRM{f}_p$ achieves the proof.
\qed
\end{proof}

\begin{remark}
Since the  homogeneous $F$-Sobolev inequality \eqref{in:Fsob}
recover the log-Sobolev inequality
\eqref{in:lsob} (with $F=\log$ and $\widetilde C_F=0$),
it is natural to ask whether the previous
Theorem recover the classical Gross Theorem or not.

So, take $F=\log$. Then, $\tau_q(x)=x^{p(q+1)}$,
$$
\frac{\partial}{\partial_q} \tau_q(x)
=
p \tau_q(x) \log(x)
=
\frac{1}{q+1} \tau_q(x) F(\tau_q(x)) ,
$$
and thus, we can choose $\ell(q)=\frac{1}{q+1}$ and $m=0$.
Moreover, it is easy to see that
$\tau_q'' \tau_q  = \frac{p(q+1)-1}{p(q+1)} {\tau_q'}^2
\geq \frac{p-1}{p}{\tau_q'}^2$,
hence $k(q)=4\frac{p-1}{p}$.
Applying the Theorem, we get that if $\mu$ satisfies a log-Sobolev inequality
\eqref{in:lsob} with constant $C_{LS}$ ($\widetilde C_F =0$),
then, for any function $f$ and any $t \geq 0$,
$$
\NRM{\PT{t} f}_{p(\tilde q(t)+1)} \leq \NRM{f}_{p} ,
$$
where $\tilde q(t)=-1+e^{4\frac{(p-1)t}{p C_{LS}}}$.
The function $p( \tilde q(t) + 1) = p e^{4 \frac{(p-1)t}{p C_{LS}}}$
is less than $q(t)=1+(p-1)e^{4t / C_{LS}}$ of Theorem \ref{th:gross}.
\end{remark}

Let us make some additional remarks on the hypotheses of the Theorem.

\begin{remark}
Let $m_F := |\min_{x \in (0,1)} xF(x)|$ and assume that $m_F < \infty$.
With our choice of $\tau_q(x)$ in the Theorem,
one can choose $l \equiv 1$ and $m=m_F$ in
order to have
$\tau_q(x) F(x^p) \leq \ell(q) \tau_q(x) F(\tau_q(x)) + m$.

Moreover, if $F$ is non negative, then $m_F=0$. Thus, in that particular
case, the previous Theorem
states that the Orlicz hypercontractivity is equivalent to
the tight homogeneous $F$-Sobolev inequality.
\end{remark}

\begin{remark}
   The condition $\tau_q(x) F(x^p) \leq \ell(q) \tau_q(x) F(\tau_q(x))
   + m$ is technical. It comes  from our  choice of
   $\tau_q = x^p e^{qF(x^p)}$. In view of the proof of Theorem~\ref{th:go}
  the most natural choice for  $\tau_q$ would be the solution of
   $$
   \left\{
\begin{array}{rcl}
\frac{\partial}{\partial q}\tau_q (x) & = & \tau_q(x) F(\tau_q(x)) \\
\tau_0(x) & = & x^p \;.
\end{array}
\right.
$$
Unfortunately, it  is not explicit in general. This is
why we preferred the expression  $ x^p e^{qF(x^p)}$ which has the
same asymptotics.
\end{remark}

\begin{remark}\label{rem:Fcond}
The hypothesis $\tau_q'' \tau_q \geq \frac{k(q)}{4} {\tau_q'}^2$
can be read as: $\tau_q^{1-\frac{k(q)}{4}}$ is a convex function.
Note that if $x F'(x) \rightarrow 0$ and
$x F''(x) \rightarrow 0$ when $x \rightarrow 0$,
$\tau_q^{1-\frac{k(q)}{4}}$ is no more convex if $k(q) >4(p-1)/p$.
Thus, we cannot hope for a better exponent than $k(q)=4(p-1)/p$
({\it i.e.} $1-\frac{k(q)}{4}=\frac1p$).
\end{remark}

Now, we give a  condition on $F$ which 
ensures  that $\tau_q$ satisfies the above  hypothesis.

\begin{proposition}\label{prop:Fcond}
Let $F: \dR^+ \rightarrow \dR$ be a $\C{2}$ non decreasing function
satisfying $F(1)=0$. Fix $p>1$. Define for all $q \geq 0$,
$\tau_q(x)=x^p e^{qF(x^p)}$. Assume that there exists a constant
$k \leq 4(p-1)/p$ such that for any $x \geq 0$,
$$
x F''(x) + ( 2 + \frac1p - \frac{k}{2}) F'(x) \geq 0 ,
$$
Then, for any $q \geq 0$, $\tau_q$ satisfies $\tau_q'' \tau_q \geq \frac{k}{4} {\tau'_q}^2$.
\end{proposition}

\begin{proof}
Note that for $q=0$ the conclusion is clearly true. Suppose $q>0$.
It is not difficult to check that for all $x > 0$,
$$
\frac{\tau_q''(x) \tau_q (x)}{{\tau'_q}^2(x)}
=
1+ \frac1p \frac{-1+qx^p F'(x^p) + pq x^{2p} F''(x^p)}{(1+qx^pF'(x^p))^2} .
$$
Thus, it is enough to prove that for any $x >0$,
\begin{eqnarray*}
&& -\frac1p - \frac k4 + 1 +
\left(\frac{q}{p} -2q\left( \frac{k}{4} -1 \right) \right) xF'(x) \\
&& \qquad + q x^2 F''(x) - \left( \frac{k}{4} -1 \right) q^2x^2F'(x)^2) \geq 0.
\end{eqnarray*}
Note that $-\frac1p \geq  \frac{k}{4} -1$ because $k \leq 4(p-1)/p$, hence, it
is sufficient to have
$$
\left(\frac{q}{p}  -2q\left( \frac{k}{4} -1 \right) \right) xF'(x)
+ q x^2 F''(x)
- \left( \frac{k}{4} -1 \right) q^2x^2F'(x)^2 \geq 0 .
$$
Since $x>0$, $\frac k4 -1 \leq 0$ and $F'(x)^2 \geq 0$, it is
satisfied when
$$
\left( \frac1p -2 \left( \frac{k}{4} -1 \right) \right) F'(x)
+ x F''(x) \geq 0
$$
which is our condition. This achieves the proof.
\qed
\end{proof}

\subsection{$\tau$-Entropy and $F$-Sobolev inequalities.}\label{IV3}

We developed in the previous two subsections two separate
versions of $F$-Sobolev inequalities
related to some hyperboundedness property.
Recall that the $\tau$-entropy involves
$$
F(y^2)  = \log \left(\frac{(\tau^*)^{-1}(y^2 \tau^*(1))}{|y|}\right) .
$$
It is necessary to relate one
to the other, in particular, since our criteria in
Section \ref{III} are written in terms of $\tau$
(or $\psi$), we have to link them to our general Gross-Orlicz theory.


First, for $\tau(y) = y^2 \psi(y) = y^2 e^{F(y^2)}$,
as above, it is easily seen that
$$
(\tau^*)^{-1}(y) = \sqrt{2y} \eta(y) ,
$$
where $\eta$ goes to $+\infty$ at infinity.
Furthermore (see \cite[Proposition 1 $(ii)$ p.14]{RR}), for all $y > 0$
\begin{equation}\label{4.16}
y \leq \tau^{-1}(y) (\tau^*)^{-1}(y) \leq  2y .
\end{equation}
Apply the inequality \eqref{4.16} with $y = \tau(z)$.
The first inequality yields
$(\psi(\tau^{-1}(z)))^{\frac 12} \leq \sqrt{2} \eta(z)$.
Since $\tau^{-1}$ is a non decreasing function and
$\tau^{-1}(z) \geq c_{\varepsilon} z^{\frac 12 - \varepsilon}$
for all $\frac 12 > \varepsilon > 0$ for some $c_{\varepsilon}$,
we certainly have
$\sqrt 2 \eta(z) \geq
(\psi(c_{\varepsilon} z^{\frac 12 - \varepsilon}))^{\frac 12}$.
Hence provided
\begin{equation}\label{4.18}
\psi(y)
\geq d_{\varepsilon}
\big(\psi(y^{\frac{2}{1-2\varepsilon}})\big)^{k_{\varepsilon}},
\end{equation}
for some positive $k_{\varepsilon}$ and $d_{\varepsilon}$,
we get that at least for large $|y|$ (using condition $\nabla_2$),
$\log(\psi(y)) \leq K_\varepsilon \log(\eta(y))$.

Also note that \eqref{4.16} furnishes
$\eta(z) \leq \psi(\sqrt{2z}) \leq C \psi(z)$.
Hence provided \eqref{4.18} holds, we have that at least for large $|y|$,
there exists two constant $c$ and $C$ such that
$c \log(\psi) \leq \log(\eta) \leq C \log(\psi)$.
In addition
$$
\log \left(\frac{(\tau^*)^{-1}(y^2 \tau^*(1))}{|y|}\right)
\approx
\log\big(\eta(\tau^*(1) y^2)\big),
$$
so that for $|y|$ large enough
\begin{equation*}
\log \left(\frac{(\tau^*)^{-1}(y^2 \tau^*(1))}{|y|}\right)
\approx
\log\big(\psi(c y^2)\big)
\end{equation*}
where we recall that $a \approx b$ if there exist some universal constants
$c_1,c_2$ such that $c_1a \leq b \leq c_2a$.

Finally note that for a defective $F$-Sobolev inequality
we may replace $F$ by $\widetilde F$ that behaves
like $F$ at infinity, up to the modification of both constants
$C_F$ and $\widetilde C_F$ in \eqref{in:Fsob}.
Hence provided $\psi$ satisfies \eqref{4.18} we may choose
$$
F =  \log(\psi) \quad {\textrm or } \quad F = \log(\eta) .
$$

\begin{example}\label{examplealpha}
Consider the $|x|^{\alpha}$ Boltzmann measure.
According to Example \ref{3.15}, Theorem \ref{4.3}
and the discussion above, we have obtained:
there exist some $A$ and $B$ such that for
$\int f^2 e^{-2|y|^{\alpha}}  dy  =  1$,
\begin{equation}\label{Rosen}
\int f^2(y)
\Big(\log_+(|f(y)|)\Big)^{2(1-\frac{1}{\alpha})} e^{-2|y|^{\alpha}} dy
\leq
A \int |\nabla f|^2 e^{-2|y|^{\alpha}} dy  + B .
\end{equation}

The latter \eqref{Rosen} is exactly the
inequality shown by Rosen (see \cite[Theorem 1]{Ro76}).
Rosen's proof lies on Sobolev inequalities in $\dR^n$
and results on monotone operators. Of course
$F=(\log^+)^{\beta}$ does not satisfy the regularity assumptions
in Theorem \ref{th:gross}, so that we cannot apply it.
But smoothing this function, we may obtain
similar (defective) inequalities.
This will be discussed in details in the section \ref{examples}.
\end{example}

\begin{example}\label{examplecritere}
In the more general situation studied in
Theorem \ref{OHGB} we may take $F(y) = G(\log (y))$
provided, in addition to assumption OB(V), $G$ is such that
\begin{equation}\label{excritint}
\int e^{- q G(|V|)} d\mu  <  +\infty ,
\end{equation}
for some $q>0$ and
\begin{equation}\label{excritasymp}
G(u) \geq k G(l u) ,
\end{equation}
for some $k>0$, $l>2$ and all $u$ large enough
(this is sometimes called condition $\nabla_2$).
\end{example}

\subsection{From defective to tight inequalities.}\label{IV.4}
%
It is well known that a
defective log-Sobolev inequality and a Poincar\'e inequality
together are equivalent to a tight log-Sobolev inequality. We shall finish
this section with the proof of a similar statement for
$F$-Sobolev inequalities
(we refer to section \ref{Condfsob} for additional results).

The first statement is straightforward
\begin{lemma}\label{spectralgap}
Let $\mu$ be a probability measure on $\dR^n$.
Let $F: (0,+ \infty) \to \dR$ be $\C{2}$
on a neighborhood of $1$.
Assume that $F(1)=0$ and that every smooth function $f$
satisfies
$$
\int f^2 F\left(\frac{f^2}{\int f^2d\mu}\right)
\le
\int |\nabla f|^2 d\mu.
$$
Then  for every smooth function $g$
$$
(4F'(1)+2F''(1)) \int \left(g-\int g d\mu\right)^2d\mu
\le
\int |\nabla g|^2d\mu.
$$
In particular, setting $\Phi(x)=xF(x)$, if $\Phi''(1)>0$ one has
$C_P(\mu)\le 1/(2\Phi''(1))$ where
$C_P(\mu)$ denotes the Poincar\'e constant.
\end{lemma}

\begin{proof}
We apply the $F$-Sobolev inequality to $f=1+\varepsilon g$ where
$g$ is bounded and $\int g d\mu = 0$ and we let $\varepsilon$ to zero.
\qed
\end{proof}

Conversely, we   prove two inequalities which allow to turn defective $F$-Sobolev
inequalities into tight ones, under a spectral gap hypothesis. The first one
deals with functions $F$ which vanish for small values. The second one, an analogue of
Rothaus inequality \cite{Rot85}, is suited to concave functions $F$. The two 
cases will be used in forthcoming arguments.

\begin{lemma}\label{lem:tighten}
Let $\rho>0$ and let $F:\dR^{+}\to \dR^{+}$ be non-decreasing and such that
$F_{|[0,\rho^{2}]}=0$. Set $G(t)=F(t\rho^{2}/(\rho+1)^{2})$ for $t\ge 0$.
Let $\mu$ be a probability measure on a space $X$ and $f:X\to \dR^{+}$ be square 
integrable.
Then 
$$ \int f^{2}G\left(\frac{f^{2}}{\mu(f^{2})} \right) \, d\mu
 \le \left(\frac{\rho+1}{\rho}\right)^{2} 
 \int \tilde{f}^{2}F\left(\frac{\tilde{f}^{2}}{\mu(\tilde{f}^{2})} \right) \, d\mu,$$
 where $\tilde{f}=f-\mu(f)$.
\end{lemma}

\begin{proof}
   Note that $G_{|[0,(\rho+1)^{2}]}=0$. Hence, 
the left integrand is non-zero only when $f^{2}\ge (\rho+1)^{2}\mu(f^{2})$.
This condition implies that $f\ge (\rho+1)\mu(f)$ and consequently $\tilde{f}\ge
 \frac{\rho}{\rho+1}f$. Combining this inequality with the classical $\mu(\tilde{f}^{2})
\le \mu(f^{2})$, we get
$$\int f^{2}
G\left(\frac{f^{2}}{\mu(f^{2})} \right) \, d\mu
 \le \left(\frac{\rho+1}{\rho}\right)^{2} 
 \int \tilde{f}^{2}G\left(\left(\frac{\rho+1}{\rho}\right)^{2}
\frac{\tilde{f}^{2}}{\mu(\tilde{f}^{2})} \right) \, d\mu.$$
The proof is complete.
\qed
\end{proof}

\begin{theorem}\label{thm:tightfsob0}
Let $\mu$ be a probability measure on a set $X$.
Let $F$ be as in Lemma~\ref{lem:tighten} with $\rho>0$ and let 
$G(t)=F(t\rho^{2}/(\rho+1)^{2})$. If $\mu$
satisfies a defective $F$-Sobolev inequality
\emph{and} a Poincar\'e inequality, {\it i.e.}
$$
\int f^2 F \left( \frac{f^2}{\mu(f^2)} \right) d \mu
\leq
C_F \int |\nabla f|^2 d\mu + \widetilde C_F \int f^2 d \mu ,
$$
and
$$
\int \left(f-\int f d\mu \right)^2 d\mu \le C_P \int |\nabla f|^2 d\mu ,
$$
then $\mu$ satisfies a tight $G$-Sobolev inequality,
more precisely
$$
\int f^2 G \left( \frac{f^2}{\mu(f^2)} \right) d \mu \leq
\left(\frac{\rho+1}{\rho}\right)^{2}(C_F+C_P  \widetilde C_F) \int |\nabla f|^2 d\mu .
$$ 
\end{theorem}

\begin{proof}
It is enough to consider non-negative functions $f$. Combining  Lemma~\ref{lem:tighten}
and the above hypotheses yields
\begin{eqnarray*}
    \int f^{2}G\left(\frac{f^{2}}{\mu(f^{2})} \right) \, d\mu
 &\le& \left(\frac{\rho+1}{\rho}\right)^{2} 
 \int \tilde{f}^{2}F\left(\frac{\tilde{f}^{2}}{\mu(\tilde{f}^{2})} \right) \, d\mu \\
&\le & \left(\frac{\rho+1}{\rho}\right)^{2}  \left( C_F \int |\nabla \tilde f|^{2}d\mu
  +\widetilde C_F \int \tilde{f}^{2}d\mu\right)  \\
&\le &
\left(\frac{\rho+1}{\rho}\right)^{2}(C_F+C_P  \widetilde C_F) \int |\nabla f|^2 d\mu.
\end{eqnarray*}
\qed
\end{proof}

\begin{lemma}[Rothaus-Orlicz inequality]\label{lem:Rothaus}
For any bounded function $f$, denote by $\tilde f$ the centered
$f  - \int f  d\mu$. If $F$ is $\C{2}$ on
$(0,+\infty)$ with $F(1)=0$ and satisfies\\
$(i)$ $F$ is concave non decreasing, goes to infinity at $+\infty$,\\
$(ii)$ $u F'(u)$ is bounded by $K(F)$.

Then it holds
$$
\int f^2 F \left( \frac{f^2}{\mu(f^2)} \right) d \mu
\leq
\int \tilde f^2 F \left( \frac{\tilde f^2}{\mu(\tilde f^2)} \right) d \mu
+ C_{Rot}(F) \lVert \tilde f \rVert^2_2.
$$
\end{lemma}

\begin{proof}
We follow the proof in \cite{Bak94}.
Again it is enough to prove the result for functions $f$
written as $f  =  1 +  t g$ for some bounded function
$g$ such that $\int  g d\mu = 0$ and $\int  g^2 d\mu = 1$.
We introduce
$$
u(t) \!=\! \frac{f^2}{\mu(f^2)}\! =\! \frac{(1+tg)^2} {1+t^2} ,
\; \log A(t)\! =\! F(u(t)+\varepsilon^2),
\; \log A \!=\! F(g^2+\varepsilon^2),
$$
for some $\varepsilon>0$ and define
\begin{eqnarray*}
\varphi(t) & = &
\int f^2 F\left(\frac{f^2}{\mu(f^2)} + \varepsilon^2 \right) d\mu
-
\int {\tilde f}^2
F\left(\frac{{\tilde f}^2}{\mu({\tilde f}^2)} + \varepsilon^2 \right) d\mu \\
& = &
\int (1+tg)^2 \log A(t) d\mu - t^2 \int g^2 \log A d\mu .
\end{eqnarray*}
The introduction of $\varepsilon$ is necessary for avoiding problems near 0. Simple calculations
yield
$$
\varphi'(t) =
\int \Big(
2g(1+tg) \log A(t) - 2t g^2 \log A +(1+tg)^2 \frac{A'(t)}{A(t)}
\Big) d\mu,
$$
and
\begin{eqnarray*}
\varphi''(t) & = &
\int \Big( 2g^2 \log \frac{A(t)}{A} + 4g(1+tg) \frac{A'(t)}{A(t)} \\
&& \quad + (1+tg)^2 \frac{A''(t) A(t) - {A'}^2(t)}{A^2(t)} \Big) d\mu .
\end{eqnarray*}
It is then easy to see that $\varphi(0) = F(1+\varepsilon^2)$ and
$\varphi'(0) = 0$.
Thanks to Taylor Lagrange formula,
$$
\varphi(t) = F(1+\varepsilon^2) + \frac{t^2}{2} \varphi''(s) ,
$$
for some $s$ so that what we need is an upper
bound for the second derivative, since $t^2 = \lVert \tilde f \rVert^2_2$.

On one hand one has for all $t\geq 0$
\begin{eqnarray*}
\log \frac{A(t)}{A}
=
F(u(t)+\varepsilon^2) - F(g^2+\varepsilon^2)
\leq 0
\end{eqnarray*}
if $u(t)\leq g^2$, and
\begin{eqnarray*}
\log \frac{A(t)}{A}
=
F(u(t)+\varepsilon^2) - F(g^2+\varepsilon^2)
\leq
F'(g^2+\varepsilon^2) (u(t)-g^2) ,
\end{eqnarray*}
if $u(t) > g^2$ since $F'$ is non-increasing.
For $|g|\geq 1$ it is easy to check that
$u(t)-g^2\leq 2 g^2$.
In this case we thus have $\log(A(t)/A) \leq 2 K(F)$.
If $|g|<1$, $u(t)-g^2 \leq 2$, hence
\begin{eqnarray*}
\int 2 g^2 \log \frac{A(t)}{A} d\mu
& \leq &
\int_{\genfrac{}{}{0pt}{}{|g|<1,}{u(t)>g^2}}
\!\! 2g^2 \log \frac{A(t)}{A} d\mu
+
\int_{\genfrac{}{}{0pt}{}{|g|\geq 1,}{u(t)>g^2}}
\!\! 2g^2 \log \frac{A(t)}{A} d\mu \\
& \leq &
\int_{\genfrac{}{}{0pt}{}{|g|<1,}{u(t)>g^2}}
2g^2 F'(g^2+\varepsilon^2) (u(t)-g^2) d \mu
+4 K(F) \\
& \leq &
8 K(F) .
\end{eqnarray*}

On the other hand,
$$
A_1 = \int \Big(
4g(1+tg) \frac{A'(t)}{A(t)} - (1+tg)^2 \left(\frac{A'(t)}{A(t)}\right)^2
\Big) d\mu
\leq 4 .
$$
Indeed define $Z = \left(\int (1+tg)^2
\left(\frac{A'(t)}{A(t)}\right)^2 d\mu \right)^{\frac 12}$
and remark that, just using Cauchy-Schwarz and $\int g^2 d\mu = 1$,
$A_1 \leq 4Z - Z^2$ which is less than 4.

It remains to control the final term
$A_2 = \int (1+tg)^2 \frac{A''(t)}{A(t)}  d\mu$.
This term may be written with the help of $F$, namely
\begin{equation*}
A_2 \!= \!\!\int \!
(1+tg)^2 \! \big[ {u'}^2(t)\big( F''+{F'}^2 \big)(u(t)+\varepsilon^2)
+ u''(t) F'(u(t)+\varepsilon^2) \big] d\mu .
\end{equation*}
Since $F''\leq 0$ we only look at terms involving $F'$.
But
$$
u'(t) = \frac{2(1+tg)(g-t)}{(1+t^2)^2}
$$
and
$\displaystyle (1+tg)^2 u''(t) =
2u(t) \big[ \frac{(3t^2-1) - 2gt(t^2-3) - g^2(3t^2-1)}{(1+t^2)^2} \big]$.
According to assumption $(ii)$,
\begin{eqnarray*}
(1+tg)^2 {u'}^2(t) {F'}^2(u(t)+\varepsilon^2)
& = &
u(t)^2  {F'}^2(u(t)+\varepsilon^2) \frac{4(g-t)^2}{(1+t^2)^2} \\
& \leq &
K^2(F) \frac{4(g-t)^2}{(1+t^2)^2} \leq 4K^2(F) (1+g^2)
\end{eqnarray*}
while
\begin{eqnarray*}
&& (1+tg)^2 u''(t) F'(u(t)+\varepsilon^2) \\
&& \qquad\qquad
\leq
2 K(F) \Big[ \frac{3t^2+1}{(1+t^2)^2} +  \frac{2|g|t(t^2+3)}{(1+t^2)^2}
+ g^2  \frac{(3t^2+1)}{(1+t^2)^2} \Big] \\
&& \qquad \qquad\leq  6 K(F)(1+|g|+ g^2) \leq 12 K(F) (1+g^2) .
\end{eqnarray*}
Integrating with respect to $\mu$ yields that
$A_2$ is uniformly bounded from above, with
a bound that does not depend on $\varepsilon$.
It remains to let $\varepsilon$ go to $0$.
\qed
\end{proof}

\begin{remark}\label{rem:rotalpha}
Remark that a smoothed version of
$F=(\log_+)^{2(1-\frac{1}{\alpha})}$ will satisfy the hypotheses
of the Lemma, for $1\leq \alpha \leq 2$ (see section~\ref{examples}).
\end{remark}

\begin{remark}\label{rem:rottau}
Using the notations in the previous subsection,
we have seen conditions for $F=\log(\eta)$ to be an
appropriate choice.
In this case using the fact that
$y \mapsto \sqrt{y} \eta(y)$ is concave and
non decreasing, it is easy to check that
$y^2  (\eta''/\eta)(y)  \leq  3/4$ and
$y (\eta'/\eta)(y) \leq  (1/2) \sqrt{y}$.
Though we are not exactly in the situation  of the
Lemma one can however check  with  more efforts
that a similar statement for $\ent_{\tau}$ is
available.
\end{remark}

To conclude this section we may state

\begin{theorem}\label{thm:tightfsob}
Let $\mu$ be a probability measure on a set $X$.
Let $F$ be as in Lemma \ref{lem:Rothaus}. If $\mu$
satisfies a defective $F$-Sobolev inequality
\emph{and} a Poincar\'e inequality, {\it i.e.}
$$
\int f^2 F \left( \frac{f^2}{\mu(f^2)} \right) d \mu
\leq
C_F \int |\nabla f|^2 d\mu + \widetilde C_F \int f^2 d \mu ,
$$
and
$$
\int \left(f-\int f d\mu \right)^2 d\mu \le C_P \int |\nabla f|^2 d\mu ,
$$
then $\mu$ satisfies a tight $F$-Sobolev inequality,
more precisely
$$
\int f^2 F \left( \frac{f^2}{\mu(f^2)} \right) d \mu \leq
C'_F \int |\nabla f|^2 d\mu
$$
with $C'_F=C_F + C_P(\widetilde C_F+C_{Rot}(F))$. 
\end{theorem}

\begin{proof}
Using the notation $\tilde f$ as in the previous Lemma, we have
\begin{eqnarray*}
\int f^2 F \left( \frac{f^2}{\mu(f^2)} \right) d \mu
& \leq &
\int \tilde f^2 F \left( \frac{\tilde f^2}{\mu(\tilde f^2)} \right) d \mu
+ C_{Rot}(F) \lVert \tilde f \rVert_2^2 \\
& \leq &
C_F \mathcal E(\tilde f, \tilde f)
+ (\widetilde C_F+C_{Rot}(F))\lVert \tilde f \rVert_2^2 \\
& \leq &
\big(C_F + C_P(\widetilde C_F+C_{Rot}(F))\big) \mathcal E (f,f).
\end{eqnarray*}
\qed
\end{proof}


\section{Sobolev inequalities}\label{sec:fsobolev}

A  measure $\mu$ on $\dR^n$ satisfies a logarithmic Sobolev inequality
for the usual Dirichlet
form if  there exists a constant $C>0$ such that for every smooth function
$$
\int f^{2} \log\left( \frac{f^{2}}{\int f^2 d\mu}\right)  d\mu
 \le
C \int |\nabla f|^2 d\mu.
$$
The latter can be rewritten as
$$
\int f^{2} \log f^2  d\mu -\left(\int f^2\, d\mu\right)
\log\left(\int f^2\, d\mu\right)
\le
C
\int |\nabla f|^2 d\mu,
$$
and also as
$$
\lim_{p\to 2^-}
\frac{\int f^2d\mu -\left( \int |f|^p d\mu\right)^{\frac2p}}{2-p}
\le
2C \int |\nabla f|^2 d\mu.
$$
Each of these forms naturally leads to considering more general inequalities. We present them before
studying their properties in details. We shall say that $\mu$ satisfies a \emph{homogeneous}
$F$-Sobolev inequality when every smooth function satisfies
\begin{equation}\label{Fsob}
\int f^2 F\left(\frac{f^2}{\int f^2 d\mu}\right)d\mu
\le
\int |\nabla f|^2 d\mu .
\end{equation}
 i.e. in this section we only consider the tight $F$-Sobolev inequality introduced in
 \eqref{in:Fsob}.

The measure $\mu$ is said to verify an \emph{additive}  $\varphi$-Sobolev inequality when for all
$f$'s
\begin{equation}\label{PHIsob}
\int f^2\varphi\Big(f^2\Big) d\mu -
\left(\int f^2d\mu\right) \varphi\left( \int f^2d\mu\right)
\le
\int |\nabla f|^2 d\mu.
\end{equation}
Finally we consider the following generalization of
Beckner's inequality: for every smooth $f$
\begin{equation}\label{Bec}
\sup_{p\in(1,2)}
\frac{\int f^2d\mu -\left( \int |f|^p d\mu\right)^{\frac2p}}{T(2-p)}
\le
\int |\nabla f|^2 d\mu.
\end{equation}
This property was introduced by Beckner \cite{Bec89} for the
Gaussian measure and $T(r)=r$. It was
considered by Lata{\l}a and Oleszkiewicz \cite{LO00} for
$T(r)=C\, r^a$. A recent independent paper
by Wang \cite{Wa03b}   studies the general case and
gives correspondences between certain
homogeneous $F$-Sobolev inequalities and generalized
Beckner-type inequalities (and actual
equivalences for $T(r)=C\, r^a$).

\subsection{First remarks, tightness and tensorisation}

Using the homogeneity property, Inequality \eqref{Fsob} above equivalently
asserts that for every smooth
function $f$ satisfying $\int f^2d\mu=1$, one has
$\int f^2F(f^2) d\mu\le  \int |\nabla f|^2 d\mu$. It is then obvious that
when $\mu$ verifies an additive   $\varphi$-Sobolev
inequality as (\ref{PHIsob}) then it satisfies a homogeneous
 $F$-Sobolev inequality with $F=\varphi-\varphi(1)$.

Inequality  (\ref{Fsob}) is tight
(it is an equality for constant functions) whenever $F(1)=0$.
Inequalities (\ref{PHIsob}) and (\ref{Bec}) are
tight by construction. Big differences appear about
tensorisation. The homogeneous $F$-Sobolev inequality
need not tensorise in general. The generalized
Beckner inequality (\ref{Bec}) has the tensorisation property.
This is established in \cite{LO00} as
a consequence of the following

\begin{lemma}
Let $\Phi: [0,\infty)\to \dR$ having a strictly positive second
derivative and such that $1/\Phi''$ is concave. Let $(\Omega_1,\mu_1),
(\Omega_2,\mu_2)$ be probability spaces. Then for any non-negative
random variable $Z$
defined on the product space
$(\Omega,\mu)=(\Omega_1\times \Omega_2,\mu_1\otimes \mu_2)$
with finite expectation one has
$$
\mathbb E_\mu\Phi(Z)-\Phi(\mathbb E_\mu Z)
\le
\mathbb E_\mu\left( \mathbb E_{\mu_1} \Phi(Z)\!-\!\Phi(\mathbb E_{\mu_1} Z)
\! +\!
\mathbb E_{\mu_2}\Phi(Z)\!-\!\Phi(\mathbb E_{\mu_2} Z)\right).
$$
\end{lemma}

When $\Phi(x)=x\varphi(x)$ satisfies the hypothesis of the lemma,
one can prove that
the corresponding additive $\varphi$-Sobolev inequality tensorises, even for
very general Dirichlet forms.
In our case, we can use the properties of the square
of the gradient to prove the tensorisation property for arbitrary $\Phi$.

\begin{lemma}\label{lem:tensorisation}
Consider for $i=1,2$ probability spaces $(\dR^{n_i},\mu_i)$.
Assume that for $i=1,2$
and every smooth function $f:\dR^{n_i}\to \dR$ one has
\begin{equation}\label{eq:LemTens}
\int \Phi\Big(f^2\Big) d\mu_i-\Phi\left( \int f^2d\mu_i\right)
\le
\int |\nabla f|^2 d\mu_i,
\end{equation}
then the measure $\mu_1\otimes \mu_2$ enjoys exactly the same property.
\end{lemma}

\begin{proof}
Let $f:\dR^{n_1+n_2}\to\dR$. We start with applying
Inequality (\ref{eq:LemTens}) in the second variable.
This gives
\begin{eqnarray*}
&&
\int \Phi(f^2) d\mu_1d\mu_2
=
\int \left(\int \Phi(f^2(x,y)) d\mu_2(y) \right) d\mu_1(x) \\
&&
\quad \le
\int \left(\Phi\left(\int f^2(x,y) d\mu_2(y)\right)
+
\int |\nabla_y f|^2(x,y) d\mu_2(y)  \right)  d\mu_1(x) \\
&&
\quad = \int \Phi(g^2) d\mu_1+\int |\nabla_y f|^2d\mu_1d\mu_2,
\end{eqnarray*}
where we have set  $g(x)=\sqrt{\int f^2(x,y) d\mu_2(y)}$.
Next we apply \eqref{eq:LemTens} on the first space to $g$.
Note that $\int g^2 d\mu_1=\int f^2 d\mu_1d\mu_2$ and that by the
Cauchy-Schwartz inequality
$$
|\nabla g|^2(x) =
\frac{\left|\int f(x,y)\nabla_xf(x,y) d\mu_2(y)\right|^2}
{\int f^2(x,y)d\mu_2(y)}
\le
\int  |\nabla_xf|^2(x,y) d\mu_2(y).
$$
Thus we get that
$\int \Phi(g^2) d\mu_1\le \Phi\left(\int f^2 d\mu_1d\mu_2\right)+
\int  |\nabla_xf|^2 d\mu_1d\mu_2$.
Combining this with the former inequality
yields the claimed $\Phi$-Sobolev inequality on the product space.
\qed
\end{proof}

\subsection{The notion of capacity of a set with
respect to a probability measure}

There exists a wide variety of Sobolev-type inequalities in the literature. It
is  natural to  analyze connections between them. To do so, one tries to define
for each inequality an equivalent ``reduced inequality'',
in such a way that it is easy to decide equivalences on the
reduced forms. For example it is known that Sobolev inequalities involving
the $\dL^1$-norm of the gradient are equivalent to isoperimetric inequalities.
 There exists a corresponding tool for
Sobolev inequalities involving $\dL^2$-norms (and even $\dL^p$-norms)
of gradients: capacities.
We refer to the book of Maz'ya \cite{mazy85ss} for more details.
The classical \emph{electrostatic capacity} of a set $A\subset \dR^n$ is
$$
\mathrm{Cap}(A) \stackrel{\rm def}{=} \inf
\left\{ \int |\nabla f(x)|^2 dx ;  f_{|A}=1
\mbox{ and } f \mbox{ has compact support}
\right\}
$$
where from now on the functions appearing in the infimum are locally Lipschitz.
The usual $\dL^2$-Sobolev inequalities on $\dR^n$
can be reduced to an inequality
relating the capacity of sets to their volume.
This was extended to more general measures and inequalities
(see \cite{mazy85ss}).
However,  if one replaces
the $dx$ in the latter formula by
$d\nu(x)$ where $\nu$ is a finite measure, then the above capacity
is zero. The appropriate notion was introduced in \cite{BR03}.
We recall it after a few definitions.
Let $\nu$ be an absolutely continuous measure on $\dR^n$.
Let $A \subset \Omega$ be Borel sets, we
write
\begin{eqnarray*}
\mathrm{Cap}_\nu(A,\Omega)
& \stackrel{\rm def}{=} &
\inf\left\{\int |\nabla f|^2d\nu ; f_{|A}\ge 1
\mbox{ and } f|_{\Omega^c}=0 \right\} \\
&=&
\inf\left\{\int |\nabla f|^2d\nu; \ind_A\le f\le \ind_\Omega \right\},
\end{eqnarray*}
where the equality follows from an easy truncation.
If $\mu$ is a probability measure on $\dR^n$,
then we set for $A$ with $\mu(A) \le 1/2$
\begin{eqnarray*}
\mathrm{Cap}_\nu(A,\mu)
&\stackrel{\rm def}{=} &
\inf\left\{\int |\nabla f|^2d\nu; f_{|A}\ge 1
\mbox{ and } \mu(f=0) \ge \frac12 \right\}\\
&=&
\inf\left\{\mathrm{Cap}_\nu(A,\Omega); A\subset \Omega
\mbox{ and } \mu(\Omega)\le \frac12  \right\}.
\end{eqnarray*}
If $\mu$ is absolutely continuous,
then since $\mathrm{Cap}_\nu(A,\Omega)$ is
non-increasing in $\Omega$,
 $$
\mathrm{Cap}_\nu(A,\mu)=
\inf\left\{\mathrm{Cap}_\nu(A,\Omega) ;  A\subset \Omega
\mbox{ and } \mu(\Omega)= \frac12  \right\}.
$$
We write
$\mathrm{Cap}_\mu(A)$ for $\mathrm{Cap}_\mu(A,\mu)$.

The reduction of an $\dL^2$-Sobolev
inequality to an inequality between capacity and
measure of sets is done via level-sets decomposition.
For completeness we illustrate
this on the simplest possible inequality (see \cite{mazy85ss}).

\begin{proposition} \label{pr:hardy-n}
Let $\mu,\nu$ be absolutely continuous measures on $\dR^n$ and let
$\Omega\subset \dR^n$.
Let $C$ denote the smallest constant so that every locally Lipschitz function
vanishing on $\Omega^c$ verifies
$$
\int f ^2 d\mu\le C\int |\nabla f|^2d\nu.
$$
Then $B\le C\le 4 B$, where $B$ is the smallest constant
so that for all $A\subset \Omega$ one has
$\mu(A)\le B\, \mathrm{Cap}_\nu(A,\Omega)$.
\end{proposition}

\begin{remark}
The constant $4$ in the above result is best possible, and is obtained
by using a result of page 109 in \cite{mazy85ss}. We shall prove  the result
with a worse constant. We follow a simplified proof, written in page 110 of
this book (this paragraph contained a small mistake which we correct
below).
\end{remark}

\begin{proof}
The fact that $B\le C$ is obvious from the definition of capacity. The other
bounds requires level-sets decomposition.
First note that replacing $f$ by $|f|$
makes the inequality tighter. So we may restrict to $f\ge 0$ vanishing outside
$\Omega$. Let $\rho>1$ and consider for $k\in \mathbb Z$,
$\Omega_k=\{f^{2}\ge \rho^k\}$. Then
\begin{eqnarray*}
\int f^2d\mu
&\le&
\sum_{k\in\mathbb Z} \rho^{k+1} \mu(\{\rho^{k}\le f^{2}<\rho^{k+1}\}) \\
& = &
\sum_{k\in\mathbb Z} \rho^{k+1} \Big(\mu(\Omega_k)-\mu(\Omega_{k+1})\Big)
= \frac{\rho-1}{\rho}  \sum_{k\in\mathbb Z} \rho^{k+1} \mu(\Omega_k).
\end{eqnarray*}
We estimate the latter measures as follows:
$$
\mu(\Omega_k)\le B \mathrm{Cap}_\nu(\Omega_k,\Omega)
\le
B \int |\nabla g_k|^2 d\nu,
$$
where we have set
$g_k=\min
\left(1,\left(\frac{f-\rho^{(k-1)/2}}{\rho^{k/2}-\rho^{(k-1)/2}}\right)_+
\right)$.
Indeed this function is $1$ on $\Omega_k$ and vanishes outside $\Omega_{k-1}$
so outside $\Omega$. Note that
\begin{eqnarray*}
\int |\nabla g_k|^2d\nu
&=&
\int_{\Omega_{k-1}\setminus \Omega_k}
\frac{|\nabla f|^2}{(\rho^{k/2}-\rho^{(k-1)/2})^2} d\nu \\
&& \quad +
\int_{f=\rho^{k/2}}  \frac{|\nabla f|^2}{(\rho^{k/2}-\rho^{(k-1)/2})^2} d\nu.
\end{eqnarray*}
Since $f$ is locally Lipschitz, the sets
$\{f=\rho^{k/2}\}\cap\{\nabla f\neq 0\}$
are Lebesgue negligible. So the latter integral vanishes (in the rest of the
paper, similar arguments are sometimes needed but we omit them).
Thus
\begin{eqnarray*}
\int f^2 d\mu
& \le &
\frac{\rho-1}{\rho}B \sum_k \frac{\rho^2}{(\sqrt\rho-1)^2}
\int_{\Omega_{k-1}\setminus \Omega_k} |\nabla f|^2d\nu \\
& \le &
B \rho \frac{\sqrt\rho+1}{\sqrt\rho-1}
\int |\nabla f|^2d\nu.
\end{eqnarray*}
The best choice of $\rho$ leads to  a constant
$(11+5\sqrt 5)/2<11.1$.
\qed
\end{proof}

\begin{remark}\label{rem:rayleigh}
Let us mention another possible reduction of Sobolev type inequalities
to inequalities of the form $R(\Omega)\ge \psi(\mu(\Omega))$ where
$R(\Omega)$ is the infimum over functions $f$ with compact support in
$\Omega$ of $\int |\nabla f|^2d\nu/\int f^2d\mu$ (Rayleigh quotient).
See e.g. \cite{Bakrcls95sid,coulgl03ipps}
where the focus is on infinite measures.
Note that by Proposition~\ref{pr:hardy-n} this criterion amounts to
inequalities of the form
$$
\mu(A) \psi(\mu(\Omega))\le \mathrm{Cap}_\nu(A,\Omega)
$$
for $A\subset \Omega$. Here the interest is in the behavior of the capacity
in terms of the outer set. We shall be rather interested in estimates of the
form $G(\nu(A))\le   \mathrm{Cap}_\nu(A,\Omega)$, that is in the dependence on
the measure of the inner sets. These two approaches are rather different, and
seem to be efficient in different settings.
\end{remark}

\begin{remark}\label{rem:hardy}
Proposition~\ref{pr:hardy-n} appears as a $n$-dimensional
version of the generalized Hardy inequality
(see Muckenhoupt \cite{muck72hiw}), which asserts that the
best $A$ so that every smooth $f$ on $\dR$ with $f(0)=0$ one has
$$
\int_0^{+\infty}f^2 d\mu\le A\int_0^{+\infty} {f'}^2 d\nu,
$$
verifies $B\le A\le 4 B$ where
$B=\sup_{x>0} \mu([x,+\infty)) \int_0^x \rho_\nu^{-1}$,
and $\rho_\nu$ is the density of the absolute continuous part of $\nu$.
Note that
$\mathrm{Cap}_\nu([x,+\infty),[0,\infty))= (\int_0^x \rho_\nu^{-1})^{-1}$,
so $B$ is the smallest constant so that
$$
\mu([x,+\infty))\le B \mathrm{Cap}_\nu([x,+\infty),[0,\infty))
$$
for all $x>0$.
This criterion is simpler than the one in $n$ dimensions, because one can
reduce to non-decreasing functions, for which level sets are half-lines.
\end{remark}

\begin{remark}\label{rem:sgcap}
It is shown in \cite{BR03} that the Poincar\'e constant of a measure $\mu$
verifies $C/2\le C_P\le K C$ where $C$ is the best constant in:
$\mu(A)\le C \mathrm{Cap_\mu(A,\mu)}$
for all $A$ with $\mu(A)<1/2$, and $K$ is a universal
constant. Proposition~\ref{pr:hardy-n} shows that one can take $K=4$.
\end{remark}

\subsection{A criterion for general Beckner-type inequalities}

The aim of this section is to give a sharp criterion for
inequalities of the form (\ref{Bec}). Since
they appear as a collection of Sobolev inequalities,
the first step consists in finding a criterion
for each Sobolev inequality. This was done by
the first and last-named authors in the case of
measures on the line.
We present here a slightly weaker but more convenient formulation of
Theorem~11  in  \cite{BR03} and its extension to arbitrary dimension.
\begin{theorem}\label{th:crit-sob}
Let $p \in (1,2)$, $\mu, \nu$ be  Borel measures on $\dR^n$,
with $\mu(\dR^n)=1$  and $d\nu(x)=\rho_\nu(x) dx$.
Let $C$ be the optimal constant
such that for every smooth $f: \dR^n\to \dR$ one has
\begin{equation} \label{in:sob-p}
\int f^{2}d\mu -\left(\int |f|^{p} d\mu\right)^{\frac2p}
\le
C \int |\nabla f|^{2}d\nu.
\end{equation}
Then   $\frac12 B(p)\le C\le 20  B(p),$ where
$B(p)$ is the optimal constant so that every Borel set $A\subset \dR^n$ with
$\mu(A)\le 1/2$ satisfies
$$
\mu(A) \left(1- \left(1+ \frac{1}{\mu(A)}  \right)^\frac{p-2}{p} \right)
\le B(p) \mathrm{Cap}_\nu(A,\mu).
$$
If $n=1$, one has
$\frac12 \max(B_-(p),B_+(p))\le C\le 20 \max(B_-(p), B_+(p))$ where
$$
B_+(p) = \sup_{x>m} \mu([x,+\infty))
       \left(1- \left(1+ \frac{1}{\mu([x,\infty))}
       \right)^\frac{p-2}{p} \right) \int_m^{x}\frac{1}{\rho_\nu},
$$
$$
B_-(p) = \sup_{x<m} \mu((-\infty,x])
       \left(1- \left(1+ \frac{1}{\mu((-\infty,x])}
       \right)^\frac{p-2}{p} \right) \int_x^{m}\frac{1}{\rho_\nu},
$$
and $m$ is a median of $\mu$.
\end{theorem}

\begin{proof}
The one dimensional result follows from
\cite[Theorem~11 and Remark~12]{BR03} which involve
$1+1/(2\mu([x,\infty))$.
In order to derive the result presented here we  have used the following
easy inequality, valid for  $ y \ge 2$, and $p \in (1,2)$,
\begin{equation}\label{in:quo}
\frac{1-(1+y)^\frac{p-2}{p}} {1- (1+ y/2)^\frac{p-2}{p}}
\le
\frac{\log 3}{\log 2}.
\end{equation}
Note that the left hand side is monotonous in $y$ and  $p$.

We turn to the $n$-dimensional part of the theorem.
We use three lemmas from \cite{BR03} which we
recall just after this proof.
We start with the lower bound on the constant $C$. Assume that the
Sobolev inequality (\ref{in:sob-p}) is satisfied for all functions.
Let $A\subset \dR^n$ with $\mu(A)\le 1/2$, and let
$f:\dR^n\to \dR$ be locally Lipschitz, with $f\ge \ind_A$ and
$\mu(f=0)\ge 1/2$.
Denote $S=\{x; f(x) \neq 0\}$. By Inequality  (\ref{in:sob-p}) and
Lemma~\ref{sup}, one has
$$
C \!\! \int \! |\nabla f|^2d\nu \ge \sup \! \left\{ \! \int \! f^2 g d\mu;
g:\dR^n \! \to (-\infty,1), \!\int (1-g)^{\frac{p}{p-2}} d\mu\le 1 \! \right\}
\!.
$$
In the latter supremum, the values of $g$ on $\{f=0\}$ have no incidence on the
integral, but they have an incidence on the constraint.
 So the supremum is achieved for $g$'s being $-\infty$ on   $\{f=0\}$. Thus
\begin{eqnarray*}
C \int |\nabla f|^2d\nu
& \ge &
\sup\left\{\int_S f^2 gd\mu;g\in(-\infty,1),
\int_S (1-g)^{\frac{p}{p-2}} d\mu\le 1  \right\}\\
& \ge &
\sup\left\{\int_S f^2 gd\mu; g\in[0,1),
\int_S (1-g)^{\frac{p}{p-2}} d\mu\le 1  \right\}\\
& \ge &
\sup\left\{\int \ind_A g \ind_S d\mu; g\in[0,1),
\int (1-g)^{\frac{p}{p-2}} \ind_S d\mu \le 1  \right\}\\
&=&
\mu(A)\left(1-\left(1+\frac{1-\mu(S)}{\mu(A)} \right)^{\frac{p-2}{p}} \right),
\end{eqnarray*}
where we have used $f\ge \ind_A$ and Lemma~\ref{supA}
for the measure $dQ=\ind_S d\mu$.
Since $\mu(S)\le 1/2$ and this is valid for any $f$ larger than 1 on $A$
and vanishing for probability $1/2$ one gets
$$
\mu(A) \left(1- \left(1+ \frac{1}{2\mu(A)}  \right)^\frac{p-2}{p} \right)
\le
C \mathrm{Cap}_\nu(A,\mu).
$$
One concludes with Inequality (\ref{in:quo}).

Next we prove the upper bound on $C$.
Let $f$ be a locally Lipschitz function. Let
$m$  be a median of the law of $f$ under $\mu$.
Set $F=f-m$, $\Omega_+=\{f>m\}$,
$\Omega_-=\{f<m\}$, $F_+=F \ind_{\Omega_+}$ and $F_-=F \ind_{\Omega_-}$.
Note that $\mu(\Omega_+),\mu(\Omega_-)\le 1/2$.
 We define the class
of functions $\mathcal I$ by
$$
\mathcal I= \left\{
g:\dR^n\to [0,1); \int (1-g)^{\frac{p}{p-2}}d\mu\le 1+(p-1)^{\frac{p}{p-2}}
\right\}.
$$
Combining Lemmas~\ref{trans} and \ref{sup}
and observing that $F^2=F_+^2+F_-^2$ gives
\begin{eqnarray*}
\int f^2d\mu-\left(\int|f|^p d\mu\right)^{\frac{2}{p}}
&\le&
\int F^2d\mu-(p-1)\left(\int|F|^p d\mu\right)^{\frac{2}{p}}\\
&\le&
\sup_{g\in \mathcal I} \int( F_+^2+F_-^2)g d\mu \\
&\le&
\sup_{g\in \mathcal I} \int F_+^2g d\mu +
\sup_{g\in \mathcal I} \int F_-^2g d\mu.
\end{eqnarray*}
Applying Proposition~\ref{pr:hardy-n} with the measures $g\,d\mu$ and $d\nu$
(it is crucial here that $g\ge 0$)  gives
$\int F_+^2 gd\mu \le 4 B_g \int |\nabla F_+|^2 d\nu$, where
\begin{eqnarray*}
B_g &:=&
\sup_{A\subset \Omega_+}
\frac{\int \ind_A g d\mu}{\mathrm{Cap}_\nu(A,\Omega_+)}
\le
\sup_{\mu(A)\le \frac12} \frac{\int \ind_A g d\mu}{\mathrm{Cap}_\nu(A,\mu)}\\
&\le&
\sup_{\mu(A)\le \frac12}
\frac{\sup\left\{\int \ind_A g d\mu; g \in \mathcal I\right\}}
{\mathrm{Cap}_\nu(A,\mu)}\\
&=&
\sup_{\mu(A)\le \frac12} \frac{\mu(A)
\left(1-\left(1+\frac{(p-1)^{\frac{p}{p-2}}}{\mu(A)} \right)^{\frac{p-2}{p}}
\right)}
{\mathrm{Cap}_\nu(A,\mu)}\le 5 B(p).
\end{eqnarray*}
In the preceding lines we have used
Lemma~\ref{supA} and the inequality
$$
1-\left(1+x(p-1)^{\frac{p}{p-2}}\right)^{\frac{p-2}{p}}
\le
5\left(1-(1+x)^{\frac{p-2}{p}}\right),
\; x\ge 2, p\in(1,2),
$$
which follows from Remark~12 of \cite{BR03}.  We have shown that
$$
\sup_{g\in \mathcal I} \int F_+^2 g d\mu
\le
20 B(p) \int |\nabla F_+|^2 d\nu.
$$
 Adding up with a similar
relation for $F_-$ leads to
\begin{eqnarray*}
\int f^2d\mu-\left(\int|f|^p d\mu\right)^{\frac{2}{p}}
& \le &
20 B(p) \left(\int |\nabla F_+|^2 d\nu+\int |\nabla F_-|^2 d\nu \right)\\
& = &
20 B(p) \int |\nabla f|^2 d\nu.
\end{eqnarray*}
\qed
\end{proof}

We list the three lemmas from \cite{BR03} that we used in the previous proof.

\begin{lemma}\label{trans}
Let $p\in(1,2)$. Let $f: X \to \dR$ be square integrable
function on a probability space $(X,Q)$. Then
for all $a\in \dR$ one has
$$
\int f^{2}dQ -\left(\int|f|^{p}dQ\right)^{\frac2p}
\le
\int (f-a)^{2}dQ-(p-1) \left(\int |f-a|^{p}dQ\right)^{\frac2p}.
$$
\end{lemma}

\begin{lemma}\label{sup}
Let  $\varphi$ be a non-negative integrable function on a probability space
$(X,P)$. Let $A>0$ and $a\in (0,1)$, then
\begin{eqnarray*}
&&\int \varphi dP- A \left(\int \varphi^{a}dP \right)^{\frac1a} \\
& =&
\sup \left\{ \int \varphi g  dP;  g: X\to(-\infty,1)
\mbox{ and } \int (1-g)^{\frac{a}{a-1}}dP\le A ^{\frac{a}{a-1}}
\right\} \\
& \le &
\sup\left\{\int \varphi g  dP;  g: X\to[0, 1) \mbox{ and }
\int (1-g)^{\frac{a}{a-1}}dP\le 1+A^\frac{a}{a-1} \right\}.
\end{eqnarray*}
\end{lemma}

\begin{lemma}\label{supA}
Let $a\in(0,1)$.
Let $Q$ be a finite measure on a space $X$ and let $K>Q(X)$.
Let $A\subset X$ be  measurable with $Q(A)>0$. Then
\begin{eqnarray*}
&& \sup\left\{\int_X \ind_A g dQ; g: X\to[0,1) \mbox{ and }
\int_X (1-g)^{\frac{a}{a-1}}dQ \le K \right\} \\
& = &
Q(A)\left(1-\left( 1+ \frac{K-Q(X)}{Q(A)}\right)^{\frac{a-1}{a}}\right).
\end{eqnarray*}
\end{lemma}

Theorem~\ref{th:crit-sob} readily implies a sharp
criterion for inequalities generalizing the ones
of Beckner and Lata{\l}a-Oleszkiewicz.

\begin{theorem}\label{th:crit-bec}
Let $T: [0,1] \to \dR^+$.
Let  $\mu, \nu$ be a Borel measures on $\dR^n$,
with $\mu(\dR^n)=1$  and $d\nu(x)=\rho_\nu(x) dx$.
Let $C$ be the optimal constant
such that for every smooth $f: \dR^n \to \dR$ one has
\begin{equation} \label{in:loth}
\sup_{p\in(1,2)}
\frac{\int f^{2}d\mu -\left(\int |f|^{p} d\mu\right)^{\frac2p}}{T(2-p)}
\le C
\int |\nabla f|^{2}d\nu.
\end{equation}
Define the function
$$
\widetilde{T}(x)=\sup_{p\in(1,2)}\frac{1-x^{\frac{p-2}{p}}}{T(2-p)}.
$$
Then   $\frac12 B(T)\le C\le 20  B(T)$,
where  $B(T)$ is the smallest  constant
so that every Borel set $A\subset \dR^n$ with
$\mu(A)<1/2$ satisfies
$$
\mu(A) \widetilde{T}\left(1+ \frac{1}{\mu(A)} \right)
\le
B(T) \mathrm{Cap}_\nu(A,\mu).
$$
If the dimension $n=1$, then
$$
\frac12 \max(B_+(T),B_-(T))\le C\le 20 \max(B_+(T),B_-(T)),
$$
where
$$
B_+(T)=\sup_{x>m}\mu([x,+\infty))
\widetilde{T}\left(1+\frac{1}{\mu([x,+\infty))} \right)
\int_m^x \frac{1}{\rho_\nu},
$$
$$
B_-(T)=\sup_{x>m}\mu((-\infty,x])
\widetilde{T}\left(1+\frac{1}{\mu((-\infty,x])} \right)
\int_x^m \frac{1}{\rho_\nu},
$$
and $m$ is a median of $\mu$.
\end{theorem}

Under fairly reasonable assumptions,
the following lemma gives a simple expression
of $\widetilde{T}$ in terms of $T$. In particular the lemma
and the theorem recover
the criterion for the Lata{\l}a-Oleszkiewicz on the real
line and extends it to any dimension.

\begin{lemma} \label{lem:sup}
Let $T: [0,1]\to\dR_+$  be non-decreasing. Then, for any $X \ge  e$,
$$
\sup_{p \in (1,2)} \frac{1-X^{\frac{p-2}{p}}}{T(2-p)}
\ge
\frac{1}{3\,T\left(\frac{1}{\log X}\right)} .
$$
If one also assumes that  $x\mapsto T(x)/x$
is non-increasing, then for  $X \ge e$
$$
\sup_{p \in (1,2)} \frac{1-X^{\frac{p-2}{p}}}{T(2-p)}
\le
\frac{1}{T\left(\frac{1}{\log X}\right)} .
$$
\end{lemma}

\begin{proof}
Let $b=\frac{2-p}{p}$, $c=b \log X$ and note that
$2-p = \frac{2b}{b+1} \le 2b$. Since $T$ is
non-decreasing, one has
\begin{eqnarray*}
\sup_{p \in (1,2)} \frac{1-X^{\frac{p-2}{p}}}{T(2-p)}
& = &
\sup_{b \in (0,1)} \frac{1-e^{-b \log X}}{T \left(\frac{2b}{b+1}\right)}
\ge
\sup_{b \in (0,1/2)} \frac{1-e^{-b \log X}}{T(2b)} \\
& \ge &
\frac{1-\sqrt e}{T \left(\frac{1}{\log X}\right)} ,
\end{eqnarray*}
by choosing $b=1/(2 \log X) \le 1/2$.
Finally $1-\sqrt e \simeq 0.393 \ge \frac13$.

For  the second assertion, let $b=\frac{2-p}{p}\in(0,1)$,
$c=b \log X$ and note that $2-p =
\frac{2b}{b+1} \ge b$. Since $T$ is non-decreasing,
\begin{eqnarray*}
\sup_{p \in (1,2)} \frac{1-X^{\frac{p-2}{p}}}{T(2-p)}
& = &
\sup_{b \in (0,1)} \frac{1-e^{-b \log X}}{T \left(\frac{2b}{b+1}\right)}
\le
\sup_{b \in (0,1)} \frac{1-e^{-b \log X}}{T(b)} \\
& \le &
\max \left[
\sup_{c \in (0,1]} \frac{1-e^{-c}}{T \left(\frac{c}{\log X}\right)},
\sup_{c \in (1,\log X)} \frac{1-e^{-c}}{T \left(\frac{c}{\log X}\right)}
\right] .
\end{eqnarray*}
Recall that  $T(x)/x$ is non-increasing. So for $c\in(0,1]$,
$T \left(\frac{c}{\log X}\right) \ge c T \left(\frac{1}{\log X}\right)$.
Hence,
$$
\sup_{c \in (0,1]} \frac{1-e^{-c}}{T \left(\frac{c}{\log X}\right)}
\le
\frac{1}{T\left(\frac{1}{\log X}\right)} \sup_{c \in (0,1]} \frac{1-e^{-c}}{c}
=
\frac{1}{T\left(\frac{1}{\log X}\right)} .
$$
When $c\ge 1$, one has
$T\left(\frac{c}{\log X}\right) \ge  T \left(\frac{1}{\log X}\right)$
since  $T$ is non-decreasing.
Thus
$$
\sup_{c \in (1,\log X)} \frac{1-e^{-c}}{T \left(\frac{c}{\log X}\right)}
\le
\frac{1}{T\left(\frac{1}{\log X}\right)} \sup_{c > 1 } (1-e^{-c})
\le
\frac{1}{T\left(\frac{1}{\log X}\right)} .
$$
This achieves the proof.
\qed
\end{proof}

\subsection{Homogeneous $F$-Sobolev inequalities}\label{Condfsob}

In the next statement, we show how to derive special homogeneous
$F$-Sobolev inequalities, which ignore the behavior of functions
close to their average.  Such inequalities appear in the work of Wang.
Let us also note that any behavior of $F$ at infinity may occur.

\begin{theorem}\label{th:crit-fsob}
   Let $D>0$ and $\rho>1$. Let $F:[0,+\infty)\to [0,+\infty)$ be a
   non-decreasing function.  Assume that $F(x)= 0$ if $x\le 2\rho$.
   Let $\mu$ be a probability measure on $\dR^n$ such that every
   $A\subset \dR^n$ with $\mu(A)\le 1/(2\rho)<1/2$
   $$
   \mu(A)F\left(\frac{\rho}{\mu(A)}\right)
   \le D
   \mathrm{Cap}_\mu(A).
   $$
   Then for every smooth $f: \dR^n \to \dR$ one
   has
   $$
   \int f^2F\left(\frac{f^2}{\int f^2d\mu} \right)d\mu
   \le D
   \left(\frac{\rho}{\sqrt\rho-1 }\right)^{2} \int |\nabla f|^2d\mu.
   $$
\end{theorem}

\begin{proof}
   For $k\ge 1$, set $\Omega_k=\{x; f^2(x)\ge 2\rho^k \mu(f^2)\}$.
   Chebichev inequality gives $\mu(\Omega_k)\le 1/(2\rho^k)$.  Next,
   since $F$ vanishes on $[0,2\rho]$
\begin{eqnarray*}
   \int f^2F\left(\frac{f^2}{\int f^2d\mu} \right)d\mu
   &\le&
   \sum_{k\ge 1} \int_{\Omega_k\setminus \Omega_{k+1}}
   f^2F\left(\frac{f^2}{\int f^2d\mu} \right)d\mu\\
   &\le&
   \sum_{k\ge 1} \mu(\Omega_k) 2\rho^{k+1}\mu(f^2) F(2\rho^{k+1}).
\end{eqnarray*}
Since $k\ge 1$ and $F$ is non-decreasing, we have
$$
\mu(\Omega_k)F(2\rho^{k+1})
\le
\mu(\Omega_k)F\left(\frac{\rho}{\mu(\Omega_k)}\right)
\le D
\mathrm{Cap}_\mu(\Omega_k).
$$
Let us consider the function
$$
h_k=\min\left(1,\left(\frac{|f|-\sqrt{2\rho^{k-1}\mu(f^2)}}
      {\sqrt{2\rho^{k}\mu(f^2)}-\sqrt{2\rho^{k-1}\mu(f^2)}}\right)_+\right),
$$
it is equal to $1$ on $\Omega_k$ and zero outside $\Omega_{k-1}$.
Since for $k\ge 1$, $\mu(\Omega_{k-1})\le 1/2$, $h_k$ vanishes with
probability at least $1/2$.  Thus
$$
\mathrm{Cap}_\mu(\Omega_k)\le \int |\nabla h_k|^2d\mu
=
\frac{\int_{\Omega_{k-1}\setminus\Omega_k} |\nabla
f|^2d\mu}{2\rho^{k-1} \left(\sqrt\rho-1\right)^2\mu(f^2)}.
$$
Combining these estimates gives
\begin{eqnarray*}
\int f^2F\left(\frac{f^2}{\int f^2d\mu} \right)d\mu
& \le &
D \sum_{k\ge 1}  2\rho^{k+1}\mu(f^2) \mathrm{Cap}_\mu(\Omega_k) \\
& \le &
D \left(\frac{\rho}{\sqrt\rho-1 }\right)^{2} \int |\nabla f|^2d\mu.
\end{eqnarray*}
\qed
\end{proof}

In the following we briefly study homogeneous $F$-Sobolev inequalities
which are tight but do not ignore the values of functions close to
their $\dL^2$-norm. In this case the behavior of $F$ at $1$ is crucial.
We have already seen the next Lemma in section \ref{IV.4}

\begin{lemma}
   Let $\mu$ be a probability measure on $\dR^n$. Let
   $F: [0, + \infty) \to \dR$ be $C^2$ on a neighborhood of $1$. Assume that
   $F(1)=0$ and that every smooth function $f$ satisfies
   $$
   \int f^2 F\left(\frac{f^2}{\int f^2d\mu}\right) \le \int |\nabla
   f|^2 d\mu.
   $$
   Then for every smooth function $g$
   $$
   (4F'(1)+2F''(1)) \int \left(g-\int g\,d\mu\right)^2d\mu \le \int
   |\nabla g|^2d\mu.
   $$
   In particular, setting $\Phi(x)=xF(x)$, if
   $\Phi''(1)>0$ one has $C_P(\mu)\le 1/(2\Phi''(1))$.
\end{lemma}

If a measure satisfies a Poincar\'e inequality,
and a tight homogeneous $F$-Sobolev
inequality which ignores small values of functions, then one
can modify $F$ on small values in an almost arbitrary way:

\begin{lemma} \label{th:crit-fsob2}
   Let $D>0$ and $\rho>1$. Let $F: [0,+\infty) \to \dR$ be a
   non-decreasing function, such that $F=0$ on $[0,2\rho)$.  Let $\mu$
   be a probability measure on $\dR^n$ with Poincar\'e constant
   $C_P<\infty$ and such that every smooth function $f$ on $R^n$
   satisfies
   $$
   \int f^2F\left(\frac{f^2}{\int f^2d\mu} \right)d\mu \le D \int
   |\nabla f|^2d\mu.
   $$
   Let $\widetilde{F}: [0,+\infty)\to \dR$ be
   non-decreasing such that $\widetilde{F}(1)=0$, $\widetilde{F}$ is $C^{2}$
   on $[0,2\rho]$ and $\widetilde{F}(x)=F(x)$ for $x\ge 2\rho$.  Set
   $\Phi(x)=x\widetilde{F}(x)$.  Then for every smooth $f:\dR^n\to \dR$
   one has
   $$
   \int f^2 \widetilde{F}\Big(\frac{f^2}{\int f^2d\mu} \Big)d\mu
   \le
   \Big( (1+\sqrt{2\rho})^2 C_P \Big(\max_{[0,2\rho]} \Phi''\Big)_+
   + D \Big) \int |\nabla f|^2d\mu.
   $$
\end{lemma}

\begin{proof}
   Note that $\Phi(1)=0$ and $\Phi'(1)=\widetilde{F}'(1)\ge 0$.  We
   introduce the function $\Phi_1(x)=\Phi(x)-\Phi(1)-\Phi'(1)(x-1)$.
   Without loss of generality, we consider a function $f\ge 0$ with
   $\int f^2d\mu=1$. One has
\begin{equation}\label{eq:Phi1}
\int \Phi(f^2)d\mu = \int  \Phi_1(f^2)d\mu
=
\int_{f^2\le 2\rho}  \Phi_1(f^2)d\mu
+ \int_{f^2> 2\rho}  \Phi_1(f^2)d\mu.
\end{equation}
For the first term, using Taylor's formula and $0\le f\le
\sqrt{2\rho}$, we obtain
\begin{equation*}
\Phi_1(f^2)
\le
\left( \max_{[0,2\rho]} \Phi''\right)   \frac{(f^2-1)^2}{2}
\le
\frac{(1+\sqrt{2\rho})^2}{2}
\left(\max_{[0,2\rho]} \Phi''  \right)_+ (f-1)^2.
\end{equation*}
 Therefore
 $$
 \int_{f^2\le 2\rho} \Phi_1(f^2)d\mu
 \le
 \frac{(1+\sqrt{2\rho})^2}{2} \left(\max_{[0,2\rho]} \Phi'' \right)_+
 \int (f-1)^2 d\mu
 $$
 can be upper-bounded thanks to the Poincar\'e
 inequality. Indeed
\begin{eqnarray*}
 \int (f-1)^2d\mu&=& \int \Big(f-\mu(f^2)^{\frac12}\Big)^2\\
 & = &
 2 \left(
 \int f^2 d\mu- \int f\,d \mu \left(\int f^2 d\mu\right)^{\frac12}\right)\\
 & \le &
 2  \left( \int f^2 d\mu- \left(\int f d \mu\right)^2 \right)
 \le 2 C_P \int
 |\nabla f|^2d\mu.
\end{eqnarray*}
The second term in (\ref{eq:Phi1}) is easily handled by our hypothesis.
Indeed, since $\Phi'(1)\ge 0$
$$
\int_{f^2> 2\rho}  \Phi_1(f^2)d\mu
\le
\int_{f^2> 2\rho}  \Phi(f^2)d\mu
\le
\int f^{2} F(f^{2})  d\mu \le D \int |\nabla f|^2d\mu.
$$
\qed
\end{proof}

Finally, we show that an homogeneous $F$-Sobolev inequality implies
an inequality between capacity and measure. We believe that the
result should be true in more generality.

\begin{theorem}\label{th:fsobcap}
   Let $\mu$ be a probability measure on $\dR^n$. Let $F: \dR^+\to \dR^+$
   be a non-negative non-decreasing function such that there exists
   $\lambda\ge 4$ such that for $x\ge 2$, $F(x)/x$ is non-increasing
   and $F(\lambda x)\le \lambda F(x)/4$.  Assume that for every smooth
   function, one has
   $$
   \int f^2F\left(\frac{f^2}{\mu(f^2)} \right)d\mu \le D \int
   |\nabla f|^2 d\mu,
   $$
   then for all $A\subset \dR^n$ with $\mu(A)\le
   \frac12$ it holds
   $$
   \mu(A) F\left(\frac{1}{\mu(A)}\right)\le 4\lambda D
   \mathrm{Cap}_\mu(A).
   $$
\end{theorem}

\begin{proof}
   Let $A$ be a set of measure less than $1/2$. In order to estimate
   its capacity, we may consider non-negative functions
   $g\ge\mathbf1_A$ and $\mu(g=0)\ge 1/2$. For $k\in\mathbb Z$ we
   consider the function
   $$
   g_k= \min\left(\Big(g-2^k\sqrt{\mu(g^2)} \Big)_+
   ,2^k\sqrt{\mu(g^2)}\right).
   $$
   We also set $\Omega_k=\{x; g(x)\ge 2^k\sqrt{\mu(g^2)} \}$.
   Note that on $\Omega_{k+1}$,
   $g_k^{2}$ is constantly $2^{2k} \mu(g^2)$ and that
   $\int g_k^{2} d\mu \le \mu(\Omega_k) 2^{2k} \mu(g^2)$.
   Therefore, applying the
   $F$-Sobolev inequality (with $F\ge 0$) to $g_k$ yields
\begin{eqnarray*}
   D \int |\nabla g|^{2} d\mu
   &\ge &
   D \int |\nabla g_k|^{2} d\mu
   \ge
   \int_{\Omega_{k+1}} g_k^{2} F\left(\frac{g_k^{2}}{\mu(g_k^{2})}\right)d\mu\\
   &\ge &
   \mu(\Omega_{k+1}) 2^{2k} \mu(g^{2}) F\left(\frac{1}{\mu(\Omega_k)}\right).
\end{eqnarray*}
Setting $a_k=\mu(\Omega_k)$ and
$C=D \int |\nabla g|^{2}d\mu/\mu(g^{2})$, we have for $k\in \mathbb Z$
$$
2^{2k} a_{k+1}F(1/a_k)\le C.
$$
Lemma~\ref{lem:ak} guarantees that
$2^{2k} a_{k}F(1/a_k)\le \lambda C$ for every $k$ with $a_k>0$, that is
$$
2^{2k} \mu(g^{2}) \mu(\Omega_k) F\left(\frac{1}{\mu(\Omega_k)}\right)
\le
\lambda D \int |\nabla g|^{2} d\mu.
$$
We choose the largest $k$ with
$2^{k}\sqrt{\mu(g^{2})}\le 1$. Thus $2^{k+1}\sqrt{\mu(g^{2})}> 1$ and
$A\subset \Omega_k$. In particular $2\le 1/\mu(\Omega_k)\le 1/\mu(A)$, so these
ratios are in the range where $x\mapsto F(x)/x$ is non-increasing.
Combining these remarks with the above inequality yields
$$
\frac14 \mu(A) F\left(\frac{1}{\mu(A)} \right)
\le
\lambda D \int |\nabla g|^{2}d\mu.
$$
Since this is valid for every $g\ge \ind_A$
and vanishing on a set of measure at least $1/2$, we have shown that
$\mu(A) F\left(1/\mu(A) \right) \le 4\lambda D \mathrm{Cap}_\mu(A)$.
\qed
\end{proof}

The next lemma was inspired by the argument of Theorem~10.5 in
\cite{Bakrcls95sid}.

\begin{lemma}\label{lem:ak}
   Let $F: [2,+\infty)\to [0,+\infty)$ be a non-decreasing function
   such that $x\to F(x)/x$ is non increasing and there exists $\lambda
   \ge 4 $ such that for all $x\ge 2$ one has $F(\lambda x)\le \lambda
   F(x)/4$.  Let $(a_k)_{k\in \dZ}$ be a non-increasing (double-sided)
   sequence of numbers in $[0,1/2]$. Assume that for all $k\in\dZ$ with
   $a_k >0$ one has
   $$
   2^{2k}a_{k+1}F\left(\frac{1}{a_k}\right) \le C,
   $$
   then for all
   $k\in \dZ$ with $a_k>0$ one has
   $$
   2^{2k}a_k F\left(\frac{1}{a_k}\right) \le \lambda C.
   $$
\end{lemma}

\begin{proof} Discarding trivial cases where $F(1/a_k)$ is always zero,
   we observe that the sequence $2^{2k} F(1/a_k)$ tends to $+\infty$
   when $k$ tends to $+\infty$, and tends to zero when $k$ tends to
   $-\infty$. So we define $k_0$ as the largest integer such that
   $2^{2k}F(1/a_k)\le 2C$.  Let $k\le k_0$, then $2C\ge
   2^{2k}F(1/a_k)\ge 2^{2k}F(2)$ since $a_k\le 1/2$ and $F$ is
   non-decreasing. Moreover since $F(t)/t$ is non-increasing, we also
   have
   $$
   2^{2(k+1)}a_{k+1}F\left(\frac{1}{a_{k+1}}\right)\le 2^{2(k+1)}
   F(2)/2.
   $$
   Combining these two inequalities yields
   $$
   2^{2(k+1)}a_{k+1}F\left(\frac{1}{a_{k+1}}\right)
   \le 4 C\le \lambda C,
   $$
   so the claimed result is established for $k\le k_0+1$.
   For larger values we proceed by induction. Let $k\ge k_0+1$, for
   which the conclusion holds. If $a_{k+1}=0$ we have nothing to
   prove. Otherwise the hypothesis of the lemma gives
   $$
   \frac{1}{a_{k+1}}\ge
   \frac{2^{2k}F\left(\frac{1}{a_k}\right)}{C}.
   $$
   Since $k>k_0$ we
   know that the term on the right is larger than $2$. Using the fact
   that $t\ge 2\mapsto F(t)/t$ is non-increasing, we obtain
   $$
   a_{k+1}F\left(\frac{1}{a_{k+1}}\right)\le
   \frac{C}{2^{2k}F\left(\frac{1}{a_k}\right)}
   F\left(\frac{2^{2k}F\left(\frac{1}{a_k}\right)}{C} \right).
   $$
   Next,
   by the induction hypothesis for $k$ this is bounded from above by
   $$
   \frac{C}{2^{2k}F\left(\frac{1}{a_k}\right)}
   F\left(\frac{\lambda}{a_k} \right) \le \frac{C}{2^{2k}}\cdot
   \frac{\lambda}{4}
   $$
   where we have used $F(\lambda t)\le \lambda
   F(t)/4$.  So we have shown
   $$
   a_{k+1}F\left(\frac{1}{a_{k+1}}\right)\le2^{-2k-2} \lambda C,
   $$
   and the conclusion is valid for $k+1$.
\qed
\end{proof}

\begin{remark}
 The alternative reduction of Sobolev type inequalities to estimates
 on the Rayleigh quotient (see Remark~\ref{rem:rayleigh}) turns out
 to work better for homogeneous $F$-Sobolev inequalities. See Proposition~2.2
 in \cite{coulgl03ipps}, dealing with measures of infinite mass, but
 the proof of which extends to our setting.
\end{remark}

\begin{remark}
Applying Theorem~\ref{th:fsobcap} to the function $F=\ind_{[2,+\infty)}$
and $\lambda=4$ shows the following. If for every function one has
$$
\int_{f^{2}\ge 2\mu(f^{2})} f^{2}d\mu \le C \int |\nabla f|^{2}d\mu
$$
then for all $A\subset \dR^{n}$ with $\mu(A)\le 1/2$, one has
$\mu(A)\le 16C \mathrm{Cap}_\mu(A)$.
By Remark~\ref{rem:sgcap}, the measure $\mu$ satisfies
a Poincar\'e inequality with constant $C_P(\mu)\le 64 C$.

The converse implication also holds.
Assume that  $\mu$ satisfies for all $f$,
$\mathrm{Var}_\mu(f)  \le C_P(\mu)\int   |\nabla f|^{2}d\mu$.
Without loss of generality, we consider $f\ge 0$.
If $f^{2}\ge 2\mu(f^{2})$ then by Cauchy-Schwarz
one has $f\ge \sqrt2\, \mu(f)$ and consequently
$(f-\mu(f))^{2}\ge (1-1/\sqrt2)^{2}f^{2}$.
Hence
$\mathrm{Var}_\mu(f)\ge (1-1/\sqrt2)^{2}\int_{f^{2}\ge 2\mu(f^{2})} f^{2} d\mu$
and the Poincar\'e inequality implies
\begin{eqnarray*}
\int_{f^{2}\ge 2 \mu(f^{2})} f^{2 }d\mu
& \le &
\left(1-\frac{1}{\sqrt2}\right)^{-2}
C_P(\mu) \int |\nabla f|^{2}d\mu \\
& \le &
12 C_P(\mu) \int |\nabla f|^{2}d\mu .
\end{eqnarray*}

As a conclusion, Poincar\'e inequality enters the
framework of homogeneous $F$-Sobolev
inequalities and is (up to the constants) equivalent to
$\int f^{2}\mathrm{1}_{f^{2}\ge 2\mu(f^{2})} d\mu
\le C \int |\nabla f|^{2}d\mu$.
Note that the number $2$ is crucial in our argument.
\end{remark}

\begin{remark}
Let us present a convenient variant of Theorem~\ref{th:fsobcap}. Assume
that $\mu$  satisfies a Poincar\'e inequality  and a $F$-Sobolev inequality
as in  Theorem~\ref{th:fsobcap}. If $F$ verifies  the assumptions $F(x)/x$
non-increasing and $F(\lambda x)\le \lambda F(x)/4$ only for $x\ge x_0>2$ then
one can however conclude with a similar inequality
between capacity and measure.
To see this, introduce a function $\widetilde{F}$ on $\dR^{+}$ with
$\widetilde{F}(x):=F(x)$
for $x\ge x_0$, $\widetilde{F}(x):=F(x_0)$ for $x\in[2,x_0]$,
$\widetilde{F}(1)=0$ and
$F$ is $C^{2}$ and non-decreasing  on $[0,x_0]$.
Then by Lemma~\ref{th:crit-fsob2}, $\mu$ satisfies
a homogeneous $\widetilde{F}$-Sobolev inequality,
and $\widetilde{F}$ satisfies the assumptions
of Theorem~\ref{th:fsobcap}.
Therefore one obtains an inequality of the
form $\mu(A)\widetilde{F}(1/\mu(A))\le K \mathrm{Cap}_\mu(A)$.
In particular if $\mu(A)\le  1/x_0$ one has
$\mu(A)F(1/\mu(A))\le K \mathrm{Cap}_\mu(A)$.
 \end{remark}

\subsection{Additive $\phi$-Sobolev inequalities}

We present an extension of a method developed by Miclo and Roberto
\cite{miclr} for logarithmic Sobolev inequalities. Throughout this
section, we work with a function $\Phi(x)=x\varphi(x)$, where
$\varphi:(0,+\infty)\to \dR$ is non-decreasing, continuously
differentiable. We assume that $\Phi$ can be extended to $0$. For
$x,t>0$ we define the function
\begin{equation*}
   \Phi_t(x)= \Phi(x)-\Phi(t)-\Phi'(t) (x-t)
   = x( \varphi(x)-\varphi(t))-t \varphi'(t)(x-t).
\end{equation*}
We start with two preliminary statements about $\Phi$-entropies. The
first one is classical and easy, and we skip its proof (see also Lemma
3.4.2 in \cite{Ane}). For short, we write $\mu(g)$ for $\int g d\mu$.

\begin{lemma}\label{lem:phi}
   For every function $f$,
   $$
   \int \Phi(f^2)\,d\mu-\Phi\left(\int f^2d\mu\right)
   =\int \Phi_{\mu(f^2)}(f^2) d\mu.
   $$
   When $\Phi$ is convex, one has
   $$
   \int \Phi(f^2)\,d\mu-\Phi\left(\int f^2d\mu\right)
   =\inf_{t>0}\int \Phi_{t}(f^2) d\mu.
   $$
\end{lemma}

\begin{lemma}\label{borne-phi_t}
   Let the function $\varphi$ be non-decreasing and concave. Assume
   that there exists $\gamma\ge0$ such that $y\varphi'(y)\le \gamma$
   for all $y>0$.  Then for every $t>0$ and every $x\in[0,2t]$ one has
   $$
   \Phi_{t^2}(x^2)\le 9\gamma(x-t)^2.
   $$
\end{lemma}

\begin{proof}
   The concavity of $\varphi$ ensures that
   $\varphi(x^2)\le \varphi(t^2)+\varphi'(t^2)(x^2-t^2).$
   This yields
\begin{eqnarray*}
   \Phi_{t^2}(x^2)
   &\le& \varphi'(t^2) (x^2-t^2)^2=(x-t)^2 \varphi'(t^2) (x+t)^2\\
   &\le& (x-t)^2 \varphi'(t^2) (3t)^2 \le 9\gamma (x-t)^2,
\end{eqnarray*}
where we have used $x\le 2t$.
\qed
\end{proof}

\begin{theorem}\label{th:crit-phi}
   Let $\varphi$ be a non-decreasing, concave, $C^1$ function on
   $(0,+\infty)$ with $\varphi(8)>0$.  Assume that there exist
   constants $\gamma,M$ such that for all $x,y>0$ one has
   $$
   x\varphi'(x)\le \gamma\qquad \mathrm{ and } \qquad \varphi(xy)\le
   M+\varphi(x)+\varphi(y).
   $$
   Let $\mu$ be a probability measure on
   $\dR^n$ satisfying a Poincar\'e inequality with constant $C_P$ and
   the following relation between capacity and measure: there exists
   $D>0$ such that for all $A\subset\dR^n$ with $\mu(A)<1/4$
   $$
   \mu(A)\varphi\left(\frac{2}{\mu(A)}\right) \le D
   \mathrm{Cap}_\mu(A) ,
   $$
   then for every smooth function one has
   $$
   \int \! \Phi(f^2)\, d\mu- \Phi\Big(\int \! f^2d\mu \Big) \le
   \left(18\gamma\,C_P+24 \Big(1+\frac{M}{\varphi(8)}\Big) D \right)
   \int \! |\nabla f |^2d\mu,
   $$
   where as usual $\Phi(x)=x\varphi(x)$.
\end{theorem}

\begin{proof}
   Without loss of generality, we may consider $f\ge 0$.  Set
   $t=(\mu(f^2))^{\frac12}$.  Then
   \begin{eqnarray}
   &&  \int \Phi(f^2) d\mu-\Phi\left(\int f^2d\mu\right)
   =
   \int \Phi_{t^2}(f^2) d\mu \nonumber \\
   && \qquad\qquad =
   \int_{f^2\le 4\mu(f^2)} \Phi_{t^2}(f^2) d\mu
   +\int_{f^2>4\mu(f^2)} \Phi_{t^2}(f^2) d\mu. \label{eq:decomp}
   \end{eqnarray}
   The first term is bounded from above thanks to
   Lemma~\ref{borne-phi_t}, indeed
\begin{eqnarray*}
&&
\int_{f^2\le 4\mu(f^2)} \Phi_{t^2}(f^2) d\mu
= \int_{f\in[0,2t]} \Phi_{t^2}(f^2) \\
&\le&
9\gamma \int_{f\in[0,2t]} \Big(f-\mu(f^2)^{\frac12}\Big)^2 d\mu
\le
9\gamma \int \Big(f-\mu(f^2)^{\frac12}\Big)^2\\
&=&
18\gamma
\left(\int f^2 d\mu- \int f d \mu \left(\int f^2 d\mu\right)^{\frac12}\right)\\
&\le&
18 \gamma \left( \int f^2 d\mu- \left(\int f\,d \mu\right)^2 \right)
\le 18 \gamma C_P \int |\nabla f|^2d\mu,
\end{eqnarray*}
where we have used Cauchy-Schwartz and the Poincar\'e inequality
for $\mu$.

The second term in (\ref{eq:decomp}) is estimated as follows
\begin{eqnarray*}
\int_{f^2 > 4\mu(f^2)} \Phi_{t^2}(f^2) d\mu
& = &
\mathop\int_{f^2 > 4\mu(f^2)}
\left[f^2 \left(\varphi(f^2)-\varphi\Big(\mu( f^2)\Big) \right) \right. \\
&& \quad  \left.
-\mu(f^2) \varphi'\Big(\mu( f^2)\Big) \Big(f^2-\mu( f^2) \Big) \right]d\mu\\
&\le&
\int_{f^2 > 4\mu(f^2)} f^2
\left(\varphi(f^2)-\varphi\Big(\mu( f^2)\Big) \right) d\mu\\
&\le&
\int_{f^2 > 4\mu(f^2)} f^2
\left(\varphi\left(\frac{f^2}{\mu( f^2)}\right)+M \right) d\mu .
\end{eqnarray*}
We conclude by applying Theorem~\ref{th:crit-fsob} with $\rho=2$,
$F(x)=0$ if $x\le 4$, and $F(x)=\varphi(x)+M$ if $x>4$. Since for
$\mu(A)\le 1/4$ one has
\begin{eqnarray*}
\mu(A)F\left(\frac{2}{\mu(A)}\right)
& = &
\mu(A)\varphi\left(\frac{2}{\mu(A)}\right)
\left(1+\frac{M}{\varphi\left(\frac{2}{\mu(A)}\right)}\right)  \\
& \le &
D\left(1+\frac{M}{\varphi(8)}\right) \mathrm{Cap}_\mu(A),
\end{eqnarray*}
we obtain
$$
\int_{f^2 > 4\mu(f^2)} \Phi_{t^2}(f^2) d\mu
\le
4 (\sqrt2+1)^2 D\left(1+\frac{M}{\varphi(8)}\right) \int |\nabla f|^2d\mu.
$$
\qed
\end{proof}


\begin{remark}\label{rem:poinc}
   As already explained, the Poincar\'e constant of the measure
   $\mu$ is bounded
   above by $4 B$ where $B$ is the best constant such that every set
   $A$ with $\mu(A)\le 1/2$ verifies $\mu(A)\le B
   \mathrm{Cap}_\mu(A)$. If $\varphi(4)>0$, one has
   $$
   D:=
   \sup_{\mu(A)\le 1/2}\frac{ \mu(A)
   \varphi(2/\mu(A))}{\mathrm{Cap}_\mu(A)} \ge \varphi(4)
   \sup_{\mu(A)\le 1/2} \frac{\mu(A) }{\mathrm{Cap}_\mu(A)}
   =\varphi(4) B.
   $$
   So $C_P\le 4 D/\varphi(4)$.  In particular, if $D<+\infty$,
   then $\mu$ satisfies an additive $\phi$-Sobolev inequality.
\end{remark}

\begin{remark}
    As already mentioned, the additive $\varphi$-Sobolev inequality
   has the tensorisation property. If it is valid for a measure $\mu$
   (with second moment) then it is true for its product measures, and
   by a classical application of the Central Limit Theorem it holds
   for the Gaussian measure. For the latter it is known that the
   logarithmic Sobolev inequality, viewed as an embedding result, is
   optimal. So $\varphi$ cannot grow faster than a logarithm. Note
   that both hypothesis on $\varphi$ assumed in Theorem~\ref{th:crit-phi}
   imply that $\varphi$ is at most a logarithm.
\end{remark}

Next we present an improved criterion for measures on the real line.

\begin{theorem}\label{th:crit-phi-dim1}
   Let $\Phi$ be a continuous convex function on $[0,\infty)$, with
   $\Phi(x)=x\varphi(x)$ for $x>0$.  Assume that $\varphi$ is
   non-decreasing, concave, and $C^1$ on $(0,+\infty)$ with
   $\varphi(8)>0$.  Assume that there exist constants $\gamma,M$ such
   that for all $x,y>0$ one has
   $$
   x\varphi'(x)\le \gamma \qquad
   \mathrm{ and } \qquad \varphi(xy)\le M+\varphi(x)+\varphi(y).
   $$
   Let $\mu$ be a probability measure on
   $\dR$, with density $\rho_\mu$, and median $m$.  Let
\begin{eqnarray*}
   D_+
   &=&
   \sup_{x>m} \mu([x,+\infty))\varphi\left(\frac{2}{\mu([x,+\infty))}\right)
                     \int_m^x \frac{1}{\rho_\mu}\\
   D_-
   &=&
   \sup_{x<m} \mu((-\infty,x])\varphi\left(\frac{2}{\mu((-\infty,x])}\right)
                     \int_x^m \frac{1}{\rho_\mu}\\
   B_+
   &=&
   \sup_{x>m} \mu([x,+\infty)) \int_m^x \frac{1}{\rho_\mu}\\
   B_-
   &=&
   \sup_{x<m} \mu((-\infty,x]) \int_x^m \frac{1}{\rho_\mu},
\end{eqnarray*}
   and $B=\max(B_+,B_-)$, $D=\max(D_+,D_-)$. Then for every smooth
   function
   $$
   \int \Phi(f^2) d\mu- \Phi\left(\int f^2d\mu\right)
   \le
   \left(144 \gamma B+ 24 \Big(1+\frac{M}{\varphi(8)}\Big)D \right)
   \int {f'}^2d\mu.
   $$
\end{theorem}

\begin{proof}
   The argument is a refinement of the proof of
   Theorem~\ref{th:crit-phi}.
   We explain the points which differ. Without loss of generality we
   consider a non-negative function $f$ on $\dR$. We consider the
   associated function $g$ defined by
\begin{eqnarray*}
   g(x)
   &=&
   f(m)+\int_m^x f'(u)\mathbf1_{f'(u)>0} du \qquad \mathrm{if} \qquad x\ge m\\
   g(x)
   &=&
   f(m)+\int_m^x f'(u)\mathbf1_{f'(u)<0} du \qquad \mathrm{if} \qquad x< m.
\end{eqnarray*}
   Set $t=(\mu(g^2))^{\frac12}$.  Then Lemma~\ref{lem:phi} ensures that
\begin{eqnarray}
&&
\int \Phi(f^2)\, d\mu-\Phi\left(\int f^2d\mu\right)
\le
\int \Phi_{t^2}(f^2) d\mu \nonumber \\
& = &
\int_{f^2\le 4\mu(g^2)} \Phi_{t^2}(f^2) d\mu
+\int_{f^2>4\mu(g^2)} \Phi_{t^2}(f^2) d\mu.\label{eq:decomp-dim1}
\end{eqnarray}
For the first term, we use Lemma~\ref{borne-phi_t}
\begin{eqnarray*}
&&
\int_{f^2\le 4\mu(g^2)} \Phi_{t^2}(f^2) d\mu
= \int_{f\in[0,2t]} \Phi_{t^2}(f^2) \\
& \le &
9\gamma \int_{f\in[0,2t]} \Big(f-\mu(g^2)^{\frac12}\Big)^2 d\mu
\le
9 \gamma \int \Big(f-\mu(g^2)^{\frac12}\Big)^2d\mu\\
&\le&
18\gamma  \int (f-g)^2d\mu+18\gamma \int \Big(g-\mu(g^2)^{\frac12}\Big)^2d\mu.
\end{eqnarray*}
Next observe that
\begin{eqnarray*}
&& \int (f-g)^2d\mu \\
& = &
\int_m^{+\infty} \!
\Big(\int_m^x \! [f'-f'\mathbf1_{f'>0}]  \Big)^2d\mu(x)
+ \int_{-\infty}^m \!
\Big(\int_m^x \! [f'-f'\mathbf1_{f'<0}] \Big)^2d\mu(x)\\
& = &
\int_m^{+\infty} \left(\int_m^x f'\mathbf1_{f'\le0}  \right)^2d\mu(x)
+
\int_{-\infty}^m \left(\int_m^x f'\mathbf1_{f'\ge0}  \right)^2d\mu(x)\\
& \le &
4B_+ \int_m^{+\infty} {f'}^2 \mathbf1_{f'\le0} d\mu+
4B_- \int_{-\infty}^m {f'}^2 \mathbf1_{f'\ge0} d\mu
\end{eqnarray*}
where the last inequality relies on Hardy inequality (see
Remark~\ref{rem:hardy}). As in the proof of Theorem~\ref{th:crit-phi},
\begin{eqnarray*}
\int \Big(g-\mu(g^2)^{\frac12}\Big)^2d\mu
& \le &
2 C_P \int {g'}^2d\mu \\
& = &
2 C_P \left(\int_m^{+\infty} {f'}^2 \mathbf1_{f'>0}d\mu
+ \int_{-\infty}^m {f'}^2 \mathbf1_{f'<0}d\mu \right),
\end{eqnarray*}
and we also use the fact that the Poincar\'e constant $C_P$ of $\mu$
satisfies $C_P\le 4 B$.  Combining the previous three estimates gives
$$
\int_{f^2\le 4 \mu(g^2)}\Phi_{t^2}(f^2)\, d\mu
\le
144\gamma B \int {f'}^2 d\mu.
$$
Now we evaluate the second term in equation
(\ref{eq:decomp-dim1}): since $\Phi_t(x) \leq x(\varphi(x) - \varphi(t))$
for $x \geq t$,
\begin{eqnarray*}
\int_{f^2 > 4\mu(g^2)} \Phi_{t^2}(f^2)\, d\mu
& \le &
\int_{f^2 > 4\mu(g^2)} f^2
\left(\varphi(f^2)-\varphi\Big(\mu( g^2)\Big) \right) d\mu\\
& \le &
\int_{g^2 > 4\mu(g^2)} g^2
\left(\varphi(g^2)-\varphi\Big(\mu( g^2)\Big) \right) d\mu \\
& \le &\int_{g^2 > 4\mu(g^2)} g^2
\left(\varphi\left(\frac{g^2}{\mu( g^2)}\right)+M \right) d\mu
\end{eqnarray*}
where we have used $g\ge f\ge 0$ and the fact that $\varphi$ is
non-decreasing.  At this stage, we apply the decomposition into level
sets performed in the proof of Theorem~\ref{th:crit-fsob}, once on
$(m,+\infty)$ and once on $(-\infty,m)$. Note that the function $g$
being non-increasing before $m$ and non-decreasing after, the level
sets appearing in the proof are of the form $(-\infty,x]$, $x<m$, and
$[x,+\infty)$, $x>m$ for which the $\mu$-capacity is controlled by the
hypothesis of the theorem.
\qed
\end{proof}

The previous two theorems apply to logarithmic Sobolev inequality when
$\varphi(x)=\log(x)$, this is how Miclo and Roberto recovered the
sufficiency part of the Bobkov-G\"otze criterion. The next result
gives an application to tight versions of Rosen's inequality.

\begin{theorem}\label{th:crit-rosen}
   Let $\beta\in (0,1]$.
   Let $\mu$ be a probability measure on $\dR^n$. Assume that one of
   the following hypotheses holds:\\
   $(i)$
   There exists a constant $D$ so that every $A\subset \dR^n$ with
   $\mu(A)\le 1/2$ satisfies
   $$
   \mu(A)\log^\beta\left(1+\frac{2}{\mu(A)}\right)
   \le D \mathrm{Cap}_\mu(A).
   $$
   $(ii)$
   The dimension $n=1$, $\mu$ has density $\rho_\mu$.  Let $m$ be a
   median of $\mu$ and
\begin{eqnarray*}
D_+
&=&
\sup_{x>m}\mu([x,+\infty))\log^\beta\left(1+ \frac{2}{\mu([x,+\infty))}\right)
                     \int_m^x \frac{1}{\rho_\mu}\\
D_-
&=&
\sup_{x<m} \mu((-\infty,x])\log^\beta\left(1+\frac{2}{\mu((-\infty,x])}\right)
                     \int_x^m \frac{1}{\rho_\mu}.
\end{eqnarray*}
   Assume that $D=\max(D_+,D_-)$ is finite.

Then for every smooth $f: \dR^n \to \dR$ one has
$$
\int \!\! f^2\log^\beta\left(1+f^2\right)d\mu-\Big(\int \!\!f^2d\mu\Big)
\log^\beta\Big(1+\int\!\! f^2d\mu\Big)
\le
K D \int \! |\nabla f|^2d\mu,
$$
where one can take $K=96$ in case $(i)$  and $K=168$ in case $(ii)$.
\end{theorem}

\begin{proof}
   In view of Theorems~\ref{th:crit-phi}, \ref{th:crit-phi-dim1}
   and Remark~\ref{rem:poinc} all we have to do is to check a few
   properties of $\Phi_\beta(x)=x\phi_\beta(x)$ where
   $\varphi_\beta(x)=\log^\beta(1+x)$. We insist on the more
   significant ones.  The function $\varphi_\beta$ is increasing, and
   since $\beta\le 1$ it is also concave. From the obvious relation
   $$
   \log(1+xy) \le \log\Big((1+x)(1+y)\Big)\le \log(1+x)+\log(1+y),
   \quad x,y>0
   $$
   and the sub-additivity of $x\mapsto x^\beta$ for
   $\beta\in(0,1]$ we deduce that $\varphi_\beta(xy)\le
   \varphi_\beta(x)+\varphi_\beta(y)$.  Finally we check the
   differential properties.  Direct calculation gives
   $$
   x\varphi'_\beta(x)=\beta x\frac{\log^{\beta-1}(1+x)}{1+x} \le
   \beta\left(\frac{x}{1+x} \right)^\beta\le\beta\le 1,
   $$
   where we
   have used $(1+x)\log(1+x)\ge x$ for $x\ge 0$.  Finally,
   $\Phi_\beta$ is concave since
   $$
   \Phi_\beta''(x)=\frac{\beta\log^{\beta-2}(1+x)}{(1+x)^2}
   \left((2+x)\log(1+x)+(\beta-1) x \right)
   $$
   is non-negative due to
   $(2+x)\log(1+x)\ge (1+x)\log(1+x)\ge x$.
\qed
\end{proof}

\begin{remark}
   From the above capacity criterion it is plain that the 
Lata{\l}a-Oleszkiewicz inequality (with $T(p)=(2-p)^{\beta}$, $(\beta\in(0,1)$
implies the tight Rosen inequality. The converse is also true: to see this
starting from a tight  Rosen inequality we may obtain Poincar\'e inequality
and a defective $F$-Sobolev inequality for $F\ge 0$. This inequality 
may be tightened by Theorem~\ref{thm:tightfsob0}. Next by the results of the
latter section on homogeneous $F$-Sobolev inequality, one may obtain an
inequality between capacity and measure. 
\end{remark}


\subsection{A summary}

\begin{figure}
\psfrag{1}[]{Beckner-type inequality $(T)$ \eqref{Bec}}
\psfrag{2}[][][2][90]{$\Longleftarrow \!\!\!\!\! \Longrightarrow$ }
\psfrag{2b}{$
\begin{array}{l}
\mbox{See Theorem } \ref{th:crit-bec}.\\
\displaystyle \psi(x)=\sup_{p \in (1,2)} \frac{1-x^\frac{p-2}{p}}{T(2-p)},
\lambda = \frac 12
\end{array}
$}
\psfrag{3}[]{ $\forall A \mbox{ such that } \mu(A) \le
\lambda, \quad \mbox{Cap}_\mu(A) \ge C_\psi \mu(A) \psi(\frac{1}{\mu(A)})$}
\psfrag{4}[][][2][90]{$\Longleftarrow\!=$}
\psfrag{4b}{$
  \begin{array}{c}
  \mbox{Under assumptions}\\
  (H1) \mbox{ or } (H2) \mbox{ on } F
  \end{array}
$}
\psfrag{5}[][][2][270]{$ \Longleftarrow\!=$}
\psfrag{5b}{$
  \begin{array}{c}
  \mbox{Under assumption}\\
  (H3) \mbox{ on } F
  \end{array}
$}
\psfrag{6}[][][2][90]{$\Longleftarrow\!=\!=\!=\!=\!=\!=\!=\!=$}
\psfrag{6b}{$
  \begin{array}{c}
  \mbox{Under assumption}\\
  (H4) \mbox{ on } \phi
  \end{array}
$}
\psfrag{7}[]{Homogenous $F$-Sobolev inequality \eqref{Fsob}}
\psfrag{8}[]{Additive $\varphi$-Sobolev inequality \eqref{PHIsob}}
\psfrag{9}[][][2][270]{$\Longleftarrow\!=$}
\psfrag{9b}[]{$F=\varphi - \varphi(1)$}
\psfrag{10}[][][2][90]{$\Longleftarrow\!=$}
\psfrag{10b}{$
  \begin{array}{c}
  \mbox{Under assumption}\\
  (H5) \mbox{ on } F
  \end{array}
$}
\psfrag{11}[]{Poincar\'e inequality}
\psfrag{12}[][][2]{$\Longleftarrow\!=$}
\psfrag{12b}[]{Under (H6)}
\psfrag{12c}[]{on $\varphi$}
\psfrag{13}[][][2]{$\Longleftarrow \!=$}
\psfrag{13b}[]{(Take $p=1$)}
\psfrag{14}[]{Poincar\'e inequality}
\psfrag{15}[][][2][90]{$\Longleftarrow \!\!\!\!\! \Longrightarrow$ }
\psfrag{15b}{$
  \begin{array}{l}
  \mbox{See Remark } \ref{rem:sgcap}. \\
  \displaystyle \psi(x) \equiv 1,
  \lambda = \frac 12
  \end{array}
$}
\epsfig{file=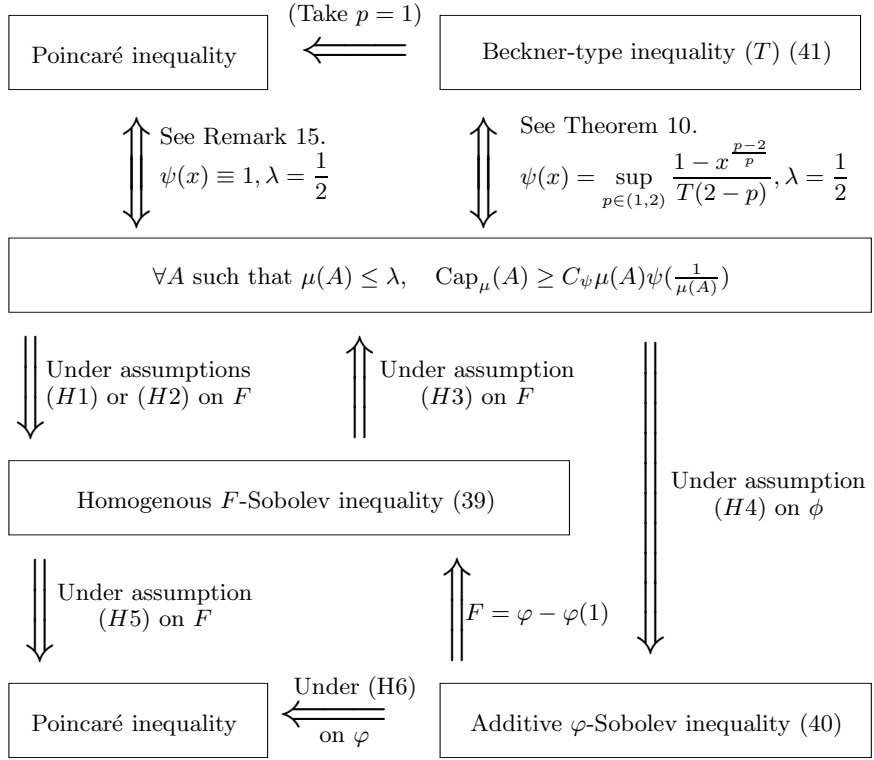,height=10cm,width=11.5cm,angle=0}
\caption{The various implications.}\label{fig:imp}
\end{figure}

In figure \ref{fig:imp} we summarize the various implications between
the inequalities studied in this section. We hope that
it will help the reader to have an overview of the picture.

First remark that thanks to Lemma \ref{lem:sup}, in figure \ref{fig:imp},
if $T: [0,1] \to \dR_+$ is
non-decreasing and $x\mapsto T(x)/x$ non-increasing, then,
$$
\frac{1}{3 T (1/ \log x)}
\leq \psi(x) = \sup_{p \in (1,2)} \frac{1-x^\frac{p-2}{p}}{T(2-p)} \leq
\frac{1}{T (1/ \log x)} \;.
$$

\begin{assumption}[H1, see Theorem \ref{th:crit-fsob}]
$F : [0,+\infty) \rightarrow \dR$ is a non-decreasing function satisfying
$F \equiv  0$  on  $[0, 2 \rho)$ for some $\rho >1$.
Finally $F(x)= \psi(x/\rho)$ for $x \ge 2 \rho$ and
$\lambda = 1/(2\rho)$.
\end{assumption}

\begin{assumption}[H2, see Theorem \ref{th:crit-fsob2}]
$F : [0,+\infty) \rightarrow \dR$ is a non-decreasing function satisfying
$F(1)= 0$ and $F$ is $\C{2}$ on $[0, 2 \rho]$. The measure
$\mu$ satisfies a Poincar\'e
inequality. Finally $F(x)= \psi(x/\rho)$  for $x \ge 0$ and
$\lambda = 1/(2\rho)$.
\end{assumption}

\begin{assumption}[H3, see Theorem \ref{th:fsobcap}]
$F : [0,+\infty) \rightarrow \dR$ is a non-decreasing function
such that there exists a constant $\gamma>4$ such that for
$x \geq 2$, $x \mapsto F(x)/x$ is non-increasing and
$F(\gamma x) \leq \gamma F(x) /4$. Then, $\psi=F$ and $\lambda = 1/2$.
\end{assumption}

\begin{assumption}[H4, see Theorem \ref{th:crit-phi}]
The function $\varphi$ is non-decrea\-sing, concave and $\C{1}$ on $(0,+ \infty)$ with
$\varphi(8)>0$. Furthermore, there exists two constants $M$ and $\gamma$ such that for any $x,y
>0$ one has
$$
x \varphi'(x) \le \gamma
\qquad \mbox{ and } \qquad
\varphi(xy) \le M + \varphi(x) + \phi(y) .
$$
The measure $\mu$ satisfies a Poincar\'e inequality. Finally $\varphi(x) = \psi(x/2)$ and $\lambda =
1/4$.
\end{assumption}

\begin{assumption}[H5, see Lemma \ref{spectralgap}]
$F : [0,+\infty) \rightarrow \dR$ is a $\C{2}$ function on a
neighborhood of $1$, $F(1)=0$ and if
$\Phi(x):=xF(x)$, $\Phi''(1)>0$.
\end{assumption}

\begin{assumption}[H6]
$\varphi : [0,+\infty) \rightarrow \dR$ is a $\C{2}$ function on a
neighborhood of $1$. Let $\Phi(x):=x\varphi(x)$.
The same proof as in Lemma~\ref{spectralgap} gives
that $\mu$ satisfies a Poincar\'e inequality
if  $\Phi''(1)>0$.
\end{assumption}


\section{Concentration property and generalized
Beckner-Lata{\l}a- Oleszkiewicz inequality.}\label{concentration}

Recall  that a probability measure
$\mu$ on $\dR^n$ satisfies a generalized Beckner inequality
if there is  a constant $C_T$ such that for any smooth  function $f$,
\begin{equation} \label{in:gbec}
\sup_{p\in(1,2)}
\frac{\int f^{2}d\mu -\left(\int |f|^{p} d\mu\right)^{\frac2p}}{T(2-p)}
\le
C_T  \int |\nabla f|^{2}d\mu.
\end{equation}
Here  $T: [0,1] \rightarrow \dR^+$ is  non decreasing, positive
on $(0,1]$ and  $T(0)=0$.
This section explores the concentration results implied by such a property.

Herbst argument, see \cite{davies-simon84,ledoux-concentration,Ane},
 derives Gaussian concentration for measures $\mu$ satisfying a log-Sobolev 
inequality along the following lines: let $h$ be  a $1$-Lipschitz function.
Applying the inequality to $\exp(\lambda h/2)$
provides the next differential
inequality for the  Laplace transform $H(\lambda)=\int \exp(\lambda h)\,d\mu$
$$
\lambda H'(\lambda) - H(\lambda) \log H(\lambda)
\le
\frac{C_{LS}}{4}\lambda^2 H(\lambda).
$$
Here $C_{LS}$ is the log-Sobolev constant.
It can be explicitly solved and gives the
subgaussian bound
$H(\lambda) \le \exp(\lambda \mu(h) + (C_{LS}/4)\lambda^2)$.
This easily yields Gaussian concentration.

On the other hand, Poincar\'e inequality only implies exponential 
concentration of Lipschitz functions. This fact goes back to Gromov and Milman
\cite{GM} (see e.g. \cite{Led99} for subsequent developments). 
In this case, the Laplace transform method works
 \cite{AS94}, but provides a relation between $H(\lambda)$ and $H(\lambda/2)$.
This approach was performed by  Lata{\l}a and Oleszkiewicz
for  their  inequality 
({\it i.e.} \eqref{in:gbec} with $T(u)=u^{2(1-\frac1\alpha)}$, $1<\alpha<2$).
 See \cite{LO00}, where optimization over $p$ is crucial. As also noted in  
\cite{Wa03b}, their argument extends as it is to general $T$. It yields

\begin{proposition} \label{prop:conc}
Let $T: [0,1] \rightarrow \dR^+$ be a non decreasing function such that
$T(0)=0$ and positive elsewhere.
Define $\theta(x)=1/T(\frac1x)$ for $x \in [1, \infty)$.
let $\mu$ be a probability measure on $\dR^n$ and
assume that there exists a constant $C_T \geq 0$ such that
for any smooth  function $f$ satisfies Inequality~\eqref{in:gbec}.
Then any $1$-Lipschitz function
$h:\dR^n \rightarrow \dR$ verifies  $\int|h|d\mu < \infty$, and \\
$(i)$ for any $t \in [0, \sqrt{T(1)}]$,
$$
\mu(\{ x: h(x) - \mu(h) \geq t\sqrt{C_T} \}) \leq e^{-\frac{t^2}{3T(1)}} ,
$$
$(ii)$ for any  $t \geq \sqrt{T(1)}$,
$$
\mu(\{ x: h(x) - \mu(h) \geq t\sqrt{C_T} \})
\leq e^{ - \sqrt2 \sup_{y\geq 1}\{t \sqrt{\theta(y)}-y\}}  .
$$
\end{proposition}

\begin{proof}
We follow the argument of \cite{LO00}. If $H(\lambda)=\mu(e^{\lambda h})$
is the Laplace transform of a $1$-Lipschitz function $h$,
Inequality  \eqref{in:gbec} for $f=\exp(\lambda h/2)$ gives
$$
H(\lambda) - H \left( \frac{p}{2}\lambda \right)^{2/p}
\le
\frac{C_T}{4} T(2-p) \lambda^2 H(\lambda) .
$$
Then, by induction, we get (see \cite{Wa03b}) for any
$\lambda < 2/\sqrt{C_T T(2-p)}$,
$$
\mu \left( e^{\lambda(h - \mu(h))} \right)
\leq
\left( 1-\frac{C_T}{4} T(2-p) \right)^{-2/(2-p)} .
$$
Chebichev inequality ensures that for any $p \in [1,2)$,
and any $\lambda < 2/\sqrt{C_T T(2-p)}$,
\begin{equation} \label{in:wang}
\mu(\{ x: h(x) - \mu(h) \geq t \sqrt{C_T} \})
\leq
e^{-\lambda t \sqrt{C_T}}
\left(1 - \frac{C_T \lambda^2}{4} T(2-p) \right)^{-\frac{2}{2-p}} .
\end{equation}
For $t < 2 \sqrt{T(1)}$, we set 
$p=1$ and $\lambda = \frac{t}{T(1) \sqrt{C_T}}$  in the latter inequality.
We get 
$$
\mu(\{ x: h(x) - \mu(h) \geq t\sqrt{C_T} \})
\leq
e^{-\frac{t^2}{T(1)}}\left(1 - \frac{t^2}{4 T(1)} \right)^{-2} .
$$
In particular, for $t<\sqrt{T(1)}$ we have
$1 - \frac{t^2}{4 T(1)} \geq e^{-\frac{t^2}{3T(1)}}$. Thus,
$$
\mu(\{ x: h(x) - \mu(h) \geq t\sqrt C_T \})
\leq
e^{-\frac{t^2}{3T(1)}} .
$$
For the second regime, choose $\lambda$  such that
$1 - \frac{C \lambda^2}{4} T(2-p) = \frac12$. It follows from \eqref{in:wang}
that for any $p \in (1,2)$
$$
\mu(\{ x: h(x) - \mu(h) \geq t\sqrt{C_T} \})
\leq
e^{-\frac{\sqrt2 t}{\sqrt{T(2-p)}} +\frac{2\ln 2}{2-p}} .
$$
Note that $2\ln 2 \leq \sqrt2$. Thus, if $y:=\frac{1}{2-p}$, we get
\begin{eqnarray*}
-\frac{\sqrt2 t}{\sqrt{T(2-p)}} +\frac{2\ln 2}{2-p}
& \leq &
- \sqrt 2 \left\{ \frac{t}{\sqrt{T(2-p)}} - \frac{1}{2-p} \right\} \\
& = &
- \sqrt 2 \{t \sqrt{\theta(y)}-y\} .
\end{eqnarray*}
One concludes the proof by optimizing in $p \in (1,2)$
or equivalently in $y \in (1, \infty)$.
\qed
\end{proof}

The next statement provides an application of the latter result 
to concentration with rate $e^{-\Phi(t)}$ for a general convex $\Phi$.
When $\Phi(t)=t^{\alpha}$, $\alpha\in(1,2)$, it reduces to the 
result by Lata{\l}a and Oleszkiewicz.

\begin{corollary}\label{coro:c}
Let $\Phi: \dR^+ \rightarrow \dR^+$ be an increasing convex function,
 with $\Phi(0)=0$.
Define $\theta (x) = \bigl(\Phi'(\Phi^{-1}(x))\bigr)^2$ for $x \in \dR^+$ and
$T(x) = 1/\theta(\frac1x)$ for $x \in \dR^+ \setminus \{0\}$, $T(0)=0$.
Here $\Phi'$ is the right derivative of $\Phi$.
Let $\mu$ be a probability measure on $\dR^n$ and
assume that there exists a constant $C_T$ such that 
it satisfies the generalized
Beckner inequality \eqref{in:gbec}.
Then, for any $1$-Lipschitz function $h:\dR^n \rightarrow \dR$,
$\int|h| < \infty$ and
for any $t \ge \sqrt{T(1)} \vee 2\Phi^{-1}(1)$,
$$
\mu(\{ x: h(x) - \mu(h) \geq t\sqrt{C_T} \})
\leq e^{ - \sqrt2\,  \Phi\left(\frac t2 \right)} .
$$
\end{corollary}

\begin{proof}
Thanks to Proposition \ref{prop:conc}, it is enough to bound from below
$\sup_{y \geq 1}\{ t \sqrt{\theta(y)}-y \}$. By assumption
$t \ge 2 \Phi^{-1}(1)$, so $\Phi(t/2) \ge 1$.
It follows that for $y=\Phi(t/2)$,
$$
\sup_{y\geq 1}\{t \sqrt{\theta(y)}-y\}
\ge
t \sqrt{\theta(\Phi(t/2))}-\Phi(t/2)
=
t \Phi'(t/2) - \Phi(t/2) .
$$
Since $\Phi$ is convex and $\Phi(0)=0$, one has   
$x \Phi'(x)\ge \Phi(x)$ for all $x\ge 0$. Hence, 
$\sup_{y\geq 1}\{t \sqrt{\theta(y)}-y\} \ge \Phi(t/2)$.
\qed
\end{proof}

Theorem~\ref{th:crit-bec} of  Section~\ref{sec:fsobolev}
 provides a criterion for a measure on the line 
to satisfy a generalized Beckner inequality.
 Under mild assumptions, and if one is not 
interested in estimating the constant, the condition may be further simplified.
\begin{proposition} \label{prop:crit}
Let $V: \dR \rightarrow \dR$ be continuous.
 Assume  that
$d\mu(x)=Z_V^{-1}e^{-V(x)} dx$ is a probability measure.
Let $T: [0,1] \rightarrow \dR^+$ be non-decreasing with $T(0)=0$ and
positive elsewhere. Assume that $x \mapsto T(x)/x$ is non-increasing.
Define $\theta(x)=1/T(1/x)$ for $x\in [1,\infty)$.
Furthermore, assume that\\
$(i)$ there exists a constant
$A>0$ such that for $|x| \geq A$, $V$ is $\C{2}$
and ${\rm sign}(x) V'(x) >0$,\\
$(ii)$ $\displaystyle \lim_{|x| \rightarrow \infty}
\frac{V''(x)}{V'(x)^2} = 0 $,\\
$(iii)$ $\displaystyle \limsup_{|x| \rightarrow \infty}
\frac{\theta(V(x) + \log V'(x) + \log Z_V)}{V'(x)^2} < \infty$.\\
Then  $\mu$
satisfies the following Beckner-type inequality:
there exists a constant $C_T \geq 0$ such that
for any smooth  function $f$,
$$
\sup_{p\in(1,2)}
\frac{\int f^{2}d\mu -\left(\int |f|^{p} d\mu\right)^{\frac2p}}{T(2-p)}
\le
C_T  \int  {f'}^{2}d\mu.
$$
\end{proposition}

\begin{proof}
 The proof is similar to the one of \cite[Proposition 15]{BR03}.
Let $m$ be a median of $\mu$. Under Hypotheses $(i)$ and $(ii)$,
when $x$ tends to $\infty$, one has (see {\it e.g.} \cite[chapter 6]{Ane})
$$
\int_m^x e^{V(t)} dt \sim \frac{e^{V(x)}}{V'(x)}
\qquad
\mbox{ and }
\qquad
\int_x^\infty e^{- V(t)} dt \sim \frac{e^{-V(x)}}{V'(x)} .
$$
Thus, for $x \ge m$,
$$
\frac{\mu([x,\infty))}{T
\left(
\frac{1}{\log (1+\frac{1}{\mu([x,\infty))})}
\right)}
\int_m^x \!\!\!  Z_V e^{V(t)} dt
\sim
Z_V \frac{\theta(V(x) + \log V'(x) + \log Z_V)}{V'(x)^2} .
$$
By Hypothesis $(iii)$, this quantity is bounded on
$[A',\infty)$ for some $A'$. Since the left hand side is 
continuous in $x\in [m,A']$, it is bounded on $(m,\infty)$.
It follows from Lemma \ref{lem:sup} that the quantity $B_+(T)$
defined in Theorem \ref{th:crit-bec} is finite.
Similarly $B_-(T)<+\infty$. We conclude with Theorem \ref{th:crit-bec}.
\qed
\end{proof}

The latter results provide
a very general condition for dimension free concentration.
Starting with an increasing convex concentration rate $\Phi:\dR^{+}\to\dR^{+}$
with $\Phi(0)=0$, we introduce the function
$T(x)=1/(\Phi'(\Phi^{-1}(x)))^{2}$. Under the additional assumption that
$\sqrt\Phi$ is concave, we know that $T(x)/x$ is non-increasing.
Therefore, under  the assumptions  of Proposition~\ref{prop:crit},
  a probability measure $d\mu(x)=Z_V^{-1}e^{-V(x)}dx$  on $\dR$ 
 satisfies the Beckner inequality with function $T$. 
By the tensorization property,  the measures $\mu^{\otimes n}$
verify the same inequality and by Corollary~\ref{coro:c}, they 
satisfy a dimension free 
concentration inequality with rate $e^{-\sqrt2\, \Phi(t/2)}$.
Note that our condition about $\sqrt\Phi$ is quite natural since,
by the Central Limit Theorem,  a dimension free
concentration inequality has at most a Gaussian rate.

The next application of our criterion  provides the best expected
 concentration rate for certain log-concave distributions.
\begin{corollary}\label{coro:cc}
Let $\Phi: \dR^{+} \rightarrow \dR^{+}$ be an increasing convex function
 with $\Phi(0)=0$ and
consider the probability measure 
 $d\mu(x)=Z_\Phi^{-1} e^{-\Phi(|x|)}\, dx$.
Assume that $\Phi$ is $\C{2}$ on $[\Phi^{-1}(1), \infty)$ and 
that $\sqrt \Phi$ is concave.

Then there exits $c>0$ such that for all $n \ge 1$, every 
$1$-Lipschitz function $h:\dR^n \rightarrow \dR$ is 
$\mu^{\otimes n}$-integrable and satisfies
$$
\mu^{\otimes n} (\{ x: h(x) - \mu^{\otimes n}(h) \geq t\sqrt{c} \})
\leq e^{ - \sqrt2\, \Phi\left(\frac t2 \right)} 
$$
provided $t \ge  2 \Phi^{-1}(1) \vee 1/(\Phi'(\Phi^{-1}(1)))$.
\end{corollary}

\begin{proof}
Set $\theta(u)=(\Phi'(\Phi^{-1}(u)))^{2}$ and $T(u)=1/\theta(1/u)$ for $u>0$.
The hypotheses on $\Phi$ ensure that $T$ is non-decreasing and $T(u)/u$ is
non-increasing. We check below that $\mu$ satisfies a 
 Beckner-type  inequality
with rate function $T$. By the above argument this implies the claimed 
concentration inequality for products. 
Let us check that $V(x)=\Phi(|x|)$ satisfies the three conditions in
Proposition \ref{prop:crit}. By symmetry it is enough to work on $\dR^{+}$.
Condition $(i)$ is obvious. Condition $(ii)$ is easily checked. Indeed since
$\sqrt \Phi$ is concave, its second derivative is non-positive when it is defined.
So for large $x$ we have 
$\Phi''/{\Phi'}^2 \leq 1/(2\Phi)$. 
So  $\lim_{+ \infty} \Phi= +\infty$ implies that  
$\lim_{+ \infty} \Phi''/{\Phi'}^2 = 0$.

Now we prove that Condition   $(iii)$ of the latter proposition is verified.
Our aim is to bound from above the quantity 
$$K(x):=\frac{\theta(\Phi(x) + \log \Phi'(x) + \log Z_\Phi )}{\Phi'(x)^2}.$$
By concavity of $\sqrt\Phi$, the function ${\Phi'}^{2}/\Phi$ is non-increasing.
Thus for $x\ge \Phi^{-1}(1)$, one has 
 $\Phi'(x)^2 \le \Phi'( \Phi^{-1}(1))^2 \Phi(x)$.
 Hence for $x$ large enough
$\log \Phi'(x) + \log Z_\Phi \le \Phi(x)$,
and  $K(x)\le \theta(2\Phi(x))/ \Phi'(x)^{2}$.

Since $\Phi$ is convex, the slope function $(\Phi(x)-\Phi(0))/x=\Phi(x)/x$ is 
non-decreasing. Comparing its values at $x$ and $2x$ shows the
inequality  $2\Phi(x)\le \Phi(2x)$.
Thus  $\theta(2\Phi(x))\le \Phi'(2x)^{2}$ and for $x$ large enough
 $K(x)\le  \Phi'(2x)^{2}/ \Phi'(x)^{2}$.
 As ${\Phi'}^{2}/\Phi$ is non-increasing we know that
 $\Phi'(2x)^2 \leq \frac{\Phi(2x)}{\Phi(x)} \Phi'(x)^2$. On the other hand, $\sqrt\Phi$
being concave, the slope function $\sqrt{\Phi(x)}/x$
is non-increasing so $\sqrt{\Phi(2x)}\le 2\sqrt{\Phi(x)}$. Finally for $x$ large
$$K(x)\le \frac{ \Phi'(2x)^{2}}{ \Phi'(x)^{2}}\le \frac{\Phi(2x)}{\Phi(x)}\le 4.$$
The proof is complete.
\qed
\end{proof}

\begin{remark}
The hypotheses of  Corollary~\ref{coro:cc} are simple but could be more
general. It is plain from Proposition~\ref{prop:crit} that we need
the convexity assumptions only for large values. The argument can be adapted
to show that the measures with potential 
 $\Phi(x)=|x|^\alpha \log(1+|x|)^\beta$ with $1<\alpha<2$
and $\beta \geq 0$ satify a dimension free concentration inequality with 
decay $e^{-C\Phi(t)}$.
\end{remark}

\begin{remark} Other concentration results for products of log-concave 
measures on the line follow from Talagrand exponential inequality, see
\cite[Theorem 2.7.1, Proposition 2.7.4]{tala95cmii}. They involve a 
different notion of enlargement depending on the log-concave density
itself. However, they imply an analogue of Corollary~\ref{coro:cc},
under the similar assumption that $\Phi(\sqrt t)$ is  subadditive.

\end{remark}


\section{Examples}\label{examples}

In this section we study fundamental examples, starting with
$|x|^{\alpha}$ Boltzmann's measures in
relation with Beckner's type inequalities.
We shall show in particular how to get dimension free
inequalities.

\subsection{ $|x|^{\alpha}$ Boltzmann's measures.}

In this subsection we are looking at the following probability measure $d\nu_\alpha^{\otimes n}(x)=
\prod_{i=1}^n Z_\alpha^{-1} e^{-2u_\alpha(x)} dx_i$ on $\dR^n$, where as in section \ref{III}, $1 <
\alpha <2$ and
\begin{equation}\label{def:ualpha}
u_{\alpha}(x) = \left\{
        \begin{array}{ll}
        |x|^{\alpha} & \textrm{ for } |x|>1 \\
        \frac{\alpha (\alpha-2)}{8}  x^4  +  \frac{\alpha(4-\alpha)}{4}  x^2
        + (1-\frac 34 \alpha + \frac 18 \alpha^2) & \textrm{ for } |x|\leq 1 .
\end{array}
\right.
\end{equation}
We will study two kind of $F$ functionals, starting from the capacity-measure
 point of view. For
each of them we give functional inequalities and derive hypercontractivity 
(or hyperboundedness)
property satisfied by the semi-group.

The first function of interest for us is
\begin{eqnarray} \label{eq:falpha}
F_\alpha: \dR^+
& \rightarrow  &
\dR \nonumber\\
x
& \mapsto &
(\log(1+x))^{2(1-\frac1\alpha)} - (\log 2)^{2(1-\frac1\alpha)} .
\end{eqnarray}
Note that it is a $\C{2}$ non-decreasing function satisfying $F_\alpha(1)=0$.
It is negative for $x<1$ and positive for $x>1$.

The second function of interest is
\begin{eqnarray} \label{eq:ftildealpha}
\widetilde F_\alpha: \dR^+ & \rightarrow  & \dR \nonumber\\
x & \mapsto &  \left\{
    \begin{array}{lcl}
    0 & \mbox{if } & x \in [0, 2 \rho] \\
    (\log(x))^{2(1-\frac1\alpha)} - (\log 2 \rho)^{2(1-\frac1\alpha)} &
    \mbox{if } & x \geq 2 \rho
    \end{array}
    \right. ,
\end{eqnarray}
where $\rho>1$ is a fixed parameter.
Note that $\widetilde F_\alpha$ is continuous
but not $\C{2}$. On the other hand, it is always non-negative.

\begin{proposition} \label{prop:Fpha}
Let $1 < \alpha <2$. Let $F_\alpha$ and $\widetilde F_\alpha$ defined in \eqref{eq:falpha} and
\eqref{eq:ftildealpha} respectively. Denote by $\nu_\alpha^{\otimes n} = \otimes_{i=1}^n
\nu_{\alpha,i}$ the product measure of $n$ copies of $d\nu_\alpha(x)=Z_\alpha^{-1} e^{-2u_\alpha(x)}
dx$.

Then, there exist two constants $C=C(\alpha)$ and
$\widetilde C= \widetilde C(\alpha,\rho)$
such that for any integer $n$,
for any smooth enough function $f: \dR^n \rightarrow \dR$,
$$
\int f^2 F_\alpha \left( \frac{f^2}{\nu_\alpha^{\otimes n} (f^2)} \right) d\nu_\alpha^{\otimes n}
\leq C \int |\nabla f|^2 d\nu_\alpha^{\otimes n} ,
$$
and
$$
\int f^2 \widetilde F_\alpha \left( \frac{f^2}{\nu_\alpha^{\otimes n} (f^2)} \right)
d\nu_\alpha^{\otimes n} \leq \widetilde C \int |\nabla f|^2 d\nu_\alpha^{\otimes n} .
$$
\end{proposition}

\begin{proof}
We start with $F_\alpha$. Fix $n=1$. Then $0$ is a median of $\nu_\alpha$. When $x$ tends to
infinity, it is easy to check that
$$
\int_0^x e^{2u_\alpha(t)} dt \sim \frac{e^{2x^\alpha}}{2\alpha x^{\alpha -1}}
\qquad
\mbox{and}
\qquad
\int_x^\infty e^{-2u_\alpha(t)} dt
\sim \frac{e^{-2x^\alpha}}{2\alpha x^{\alpha -1}} .
$$
It follows that the two constants $D_+$ and $D_-$ introduced in
Theorem \ref{th:crit-rosen} with $\beta = 2(1-\frac{1}{\alpha})$ are finite.
Then, we conclude by Theorem \ref{th:crit-rosen} that
there exists a constant $C_\alpha$ such that for every function
$f$ on $\dR$,
\begin{eqnarray}\label{tensbeta}
&& \int f^2\log^\beta \left( 1+f^2 \right) d\nu_\alpha -\left( \int  f^2 d\nu_\alpha  \right)
\log^\beta \left(1+\int f^2d\nu_\alpha \right) \\
&& \qquad \qquad \qquad \qquad \qquad \qquad \qquad \qquad \qquad \le C_\alpha \int |\nabla f|^2
d\nu_\alpha. \nonumber
\end{eqnarray}
Then, for any integer $n$, by Lemma \ref{lem:tensorisation} the latter
inequality holds for $\nu_\alpha^{\otimes n}$ in $\dR^n$. Finally,
applying the inequality to $f^2/\nu_\alpha^{\otimes n}(f^2)$ gives the
expected result.

\medskip

The case of $\widetilde F_\alpha$ is a bit more difficult. Let
$\beta=2(1-\frac{1}{\alpha})$ and $T(x)=|x|^\beta$. It is easy to
check that the hypotheses of Proposition \ref{prop:crit} are satisfied
(for $\Phi=2u_\alpha$) and thus that there exists a constant
$\widetilde C= \widetilde C(\alpha)$ such that for any function $f:
\dR \rightarrow \dR$,
$$
\sup_{p\in(1,2)} \frac{\int f^{2}d\nu_\alpha - \left(\int |f|^{p}
     d\nu_\alpha\right)^{\frac2p}}{(2-p)^\beta} \le \widetilde C \int
|\nabla f|^{2}d\nu_\alpha .
$$
Now, by tensorization property (see \cite{LO00}), the same
inequality holds for $\nu_\alpha^{\otimes n}$ with the same constant
$\widetilde C$ (independent of $n$). Thus, by Theorem
\ref{th:crit-bec} together with Lemma \ref{lem:sup} (recall that
$T(x)=|x|^\beta$), it follows that for any integer $n$, any Borel set
$A \subset \dR^n$ with $\nu_\alpha^{\otimes n} (A) \le 1/2$,
$$
\nu_\alpha^{\otimes n} (A) \left(
   \log(1+\frac{1}{\nu_\alpha^{\otimes n} (A)} ) \right)^\beta \le 2
\widetilde C \mathrm{Cap}_{\nu_\alpha^{\otimes n}}(A) .
$$
Now, for any $x \ge 2 \rho$, $\widetilde F_\alpha (\rho x) \le
(\log(1+x))^\beta$. Therefore, for any Borel set $A \subset \dR^n$
with $\nu_\alpha^{\otimes n} (A) \le 1/(2\rho)$,
$$
\nu_\alpha^{\otimes n} (A) \widetilde F_\alpha \left(
   \frac{\rho}{\nu_\alpha^{\otimes n} (A)} \right) \le 2 \widetilde C
\mathrm{Cap}_{\nu_\alpha^{\otimes n}}(A) .
$$
The expected result follows from Theorem \ref{th:crit-fsob}.  This
achieves the proof.  \qed
\end{proof}

\begin{remark} \label{rm:calpha1}
   It is not difficult to check that
   $$
   0 < \inf_{\alpha \in (1,2)} C(\alpha) < \sup_{\alpha \in (1,2)}
   C(\alpha) < +\infty .
   $$
   This means that the constant $C(\alpha)$ appearing in
   Proposition \ref{prop:Fpha} can be chosen independently of $\alpha
   \in (1,2)$. This uniformity will be useful for applications.
\end{remark}

\begin{corollary} \label{cor:fal}
   Let $1 < \alpha <2$. Let $F_\alpha$ defined in \eqref{eq:falpha}.
   Denote by $\nu_\alpha^{\otimes n} = \otimes_{i=1}^n \nu_{\alpha,i}$
   the product measure of $n$ copies of the probability measure
   $d\nu_\alpha(x)=Z_\alpha^{-1} e^{-2u_\alpha(x)} dx$. Define for any
   $q \ge 0$, any $x \ge 0$, $\tau_q^{(\alpha)}(x)=x^2 e^{q
     F_\alpha(x^2)}$.
   
   Then, there exists a universal constant $C$ such that for any
   integer $n$, any function $f:\dR^n \rightarrow \dR$ and any $t\ge
   0$,
   $$
   N_{\tau_{q(t)}^{(\alpha)}}(\PT{t} f) \le e^{\frac{C}{2} t}
   \NRM{f}_2
   $$
   where $q(t)=C t$ and $N_{\Phi}(g):= \inf \{ \lambda: \int
   \Phi(g/\lambda)d\nu_\alpha^{\otimes n} \le 1 \}$.
\end{corollary}

\begin{proof}
The result is a direct consequence of Theorem \ref{th:go},
using Proposition \ref{prop:Fpha} and Lemma \ref{lem:ex} below.
\qed
\end{proof}


\begin{lemma} \label{lem:ex}
Let $1 < \alpha <2$.
Let $F_\alpha$  defined in \eqref{eq:falpha}.
Define for any $q \ge 0$, any $x \ge 0$,
$\tau_q^{(\alpha)}(x)=x^2 e^{q F_\alpha(x^2)}$.
Then,\\
$(i)$ For any $x \ge 0$, any $q \ge 0$,
$$
(\tau_q^{(\alpha)})'' \tau_q^{(\alpha)}
\geq
\frac{5- (4 / \alpha)}{4}  {(\tau_q^{(\alpha)})'}^2
\geq
\frac 14  {(\tau_q^{(\alpha)})'}^2,
$$
$(ii)$ for any $x \ge 0$, any $q \ge 0$,
$$
\tau_q^{(\alpha)}(x) F(x^2)
\le
\tau_q^{(\alpha)}(x) F(\tau_q^{(\alpha)}(x)) + 1 .
$$
\end{lemma}

\begin{proof}
Let $\beta=2(1-\frac{1}{\alpha})$. Then $0<\beta<1$.
It is easy to check that for any $x >0$,
$$
- \frac{xF_\alpha''(x)}{F_\alpha'(x)}
=
\frac{x(1-\beta + \log(1+x))}{(1+x)\log(1+x)} \le 2- \beta .
$$
We conclude the proof of point $(i)$ applying Proposition
\ref{prop:Fcond}
(note that $2+\frac{1}{2}-\frac{5- (4 / \alpha)}{2} = 2 - \beta$).

Note that
$m_{F_\alpha} := |\min_{x \in (0,1)} xF_\alpha (x)| \le 1$.
Hence, using remark \ref{rem:Fcond} concludes the proof
of point $(ii)$.
\qed
\end{proof}

The proof of a similar result than Corollary \ref{cor:fal} for
$\widetilde F_\alpha$ is a bit more difficult due to differentiation problem at
$x=2\rho$. The result is the following:

\begin{corollary} \label{cor:tfal}
   Let $1 < \alpha <2$. Let $\widetilde F_\alpha$ defined in
   \eqref{eq:ftildealpha}. Denote by $\nu_\alpha^{ \otimes n} =
   \otimes_{i=1}^n \nu_{\alpha,i}$ the product measure of $n$ copies
   of the probability measure $d\nu_\alpha(x)=Z_\alpha^{-1}
   e^{-2u_\alpha(x)} dx$. Define for any $q \ge 0$, any $x \ge 0$,
   $\widetilde \tau_q^{(\alpha)}(x)=x^2 e^{q \widetilde
     F_\alpha(x^2)}$.
   
   Then, there exists a constant $\widetilde C= \widetilde
   C(\alpha,\rho)$ such that for any integer $n$, any function $f:
   \dR^n \rightarrow \dR$ and any $t\ge 0$,
   $$
   N_{\widetilde \tau_{q(t)}^{(\alpha)}}(\PT{t} f) \le \NRM{f}_2
   $$
   where $q(t)=\widetilde C t$ and $N_{\Phi}(g):= \inf \{ \lambda:
   \int \Phi(g / \lambda) d\nu_\alpha^{\otimes n} \le 1 \}$.
\end{corollary}

\begin{proof}
Let $g$ be a $\C{\infty}$ non-negative function with compact support in
$[-1,0]$ and such that $\int g(y) dy =1$.
For any $\varepsilon >0$ define
$g_\varepsilon (x) = \frac1\varepsilon g(\frac{x}{\varepsilon})$ and
note that
$\widetilde F_\alpha \ast g_\varepsilon (x)
:= \int \widetilde F_\alpha (x-y) g_\varepsilon (y) dy$
is a $\C{\infty}$ function.

Define for any  $\varepsilon >0$, any $q \ge 0$,
$\widetilde \tau_{q,\varepsilon}^{(\alpha)}(x)
=x^2 e^{q \widetilde F_\alpha \ast g_\varepsilon(x^2)}$.

Thanks to Lemma \ref{lem:ex2} below,
$\widetilde F_\alpha \ast g_\varepsilon$ satisfies
the hypothesis of Theorem \ref{th:go}, uniformly in $n$.
Thus, by Theorem \ref{th:go} there exists two constants
$\widetilde C= \widetilde C(\alpha,\rho)$ and
$\widetilde C'= \widetilde C'(\alpha,\rho)$ (maybe different from those one of
Lemma \ref{lem:ex2}) such that for any integer $n$, any function
$f:\dR^n \rightarrow \dR$ and any $t\ge 0$,
$$
N_{\widetilde \tau_{q(t),\varepsilon}^{(\alpha)}}(\PT{t} f)
\le
e^{\frac{1}{2}
(\widetilde F_\alpha(2\rho + \varepsilon)+\varepsilon\widetilde C')t}
\NRM{f}_2 .
$$
Then, it is easy to verify that for any function $f$, any $t$,
when $\varepsilon$ tends to $0$,
$$
N_{\widetilde \tau_{q(t),\varepsilon}^{(\alpha)}}(\PT{t} f)
\rightarrow
N_{\widetilde \tau_{q(t)}^{(\alpha)}}(\PT{t} f)
\qquad  \mbox{ and }
\qquad
e^{\frac{1}{2}
(\widetilde F_\alpha(2\rho + \varepsilon)+\varepsilon\widetilde C')t}
\rightarrow  1 .
$$
This achieves the proof.
\qed
\end{proof}

\begin{lemma} \label{lem:ex2}
Let $1 < \alpha <2$. Let $\widetilde F_\alpha$ defined in \eqref{eq:ftildealpha}. Denote by
$\nu_\alpha^{\otimes n} = \otimes_{i=1}^n \nu_{\alpha,i}$ the product measure of $n$ copies of
$d\nu_\alpha(x)=Z_\alpha^{-1} e^{-2u_\alpha(x)} dx$. Define for any $q \ge 0$, any $x \ge 0$,
$\widetilde \tau_q^{(\alpha)}(x)=x^2 e^{q \widetilde F_\alpha(x^2)}$.

Let $g$ be a $\C{\infty}$ non-negative function with compact support in
$[-1,0]$ and such that $\int g(y) dy =1$.
 Define $g_\varepsilon (x) = \frac1\varepsilon g(\frac{x}{\varepsilon})$,
and
$\widetilde F_\alpha \ast g_\varepsilon (x) :=
\int \widetilde F_\alpha (x-y) g_\varepsilon (y) dy$
for any $\varepsilon >0$,
and for any $q \ge 0$,
$\widetilde \tau_{q,\varepsilon}^{(\alpha)}(x)=
x^2 e^{q \widetilde F_\alpha \ast g_\varepsilon(x^2)}$.
Then,\\
$(i)$ for any $\varepsilon>0$ and any $q \ge 0$,
$$
(\widetilde \tau_{q,\varepsilon}^{(\alpha)})''
\widetilde \tau_{q,\varepsilon}^{(\alpha)}
\ge
\frac{3- 2(2-\alpha)/(\alpha \log(2\rho))}{4}
{(\widetilde \tau_{q,\varepsilon}^{(\alpha)})'}^2 .
$$
$(ii)$ For any $\varepsilon>0$ small enough, any $q \ge 0$, and any $x \ge 0$,
$$
\widetilde F_\alpha \ast g_\varepsilon (x^2)
\le
\widetilde F_\alpha \ast g_\varepsilon
(\widetilde \tau_{q,\varepsilon}^{(\alpha)}(x)) .
$$
$(iii)$ There exist two constants
$\widetilde C= \widetilde C(\alpha,\rho)$ and
$\widetilde C'= \widetilde C'(\alpha,\rho)$ such that for any integer $n$,
any function $f: \dR^n \rightarrow \dR$ and any $\varepsilon> 0$ small enough,
\begin{eqnarray*}
\int \! f^2 \widetilde F_\alpha \ast g_\varepsilon 
\big( \frac{f^2}{\nu_\alpha^{\otimes n} (f^2)}
\big) d \nu_\alpha^{\otimes n} 
& \le &
\widetilde C \!\int\! |\nabla f|^2 d \nu_\alpha^{\otimes n} \\ 
&&
+ (\widetilde F_\alpha(2 \rho + \varepsilon) + \varepsilon \widetilde C') 
\int f^2 d \nu_\alpha^{\otimes n} .
\end{eqnarray*}
\end{lemma}

\begin{proof}
Let  $\beta=2(1-\frac{1}{\alpha})$.

We start with  $(i)$.
The result is obviously true for $x \le 2\rho$. For $x > 2\rho$, an easy
computation gives
$$
- \frac{x \widetilde F_\alpha''(x)}{\widetilde F_\alpha'(x)}
=
\frac{1-\beta + \log x}{\log x}
\le
1+ \frac{1-\beta}{\log(2\rho)}
=1+\frac{2-\alpha}{\alpha \log 2\rho} .
$$
Thus, by Lemma \ref{lem:approx} below we get that for any $\varepsilon >0$, any
$x \ge 0$,
$$
x (\widetilde F_\alpha \ast g_\varepsilon)''(x) +
\left(1+\frac{2-\alpha}{\alpha \log 2\rho} \right)
(\widetilde F_\alpha \ast g_\varepsilon)'(x)
\ge 0 .
$$
The result follows from Proposition \ref{prop:Fcond}.

\medskip

For $(ii)$ note  that for any $\varepsilon \le 2 \rho -1$,
$\widetilde F_\alpha \ast g_\varepsilon \equiv 0$ on $[0,1]$. Thus the result
becomes obvious thanks to Remark \ref{rem:Fcond}.

\medskip

Next we deal with  $(iii)$. First  note that
$\widetilde F_\alpha \ast g_\varepsilon \equiv 0$
on $[0,2 \rho - \varepsilon]$. Then, for
$x \in [2 \rho - \varepsilon, 2 \rho]$, since
$\widetilde F_\alpha$ is non-decreasing,
\begin{eqnarray*}
\widetilde F_\alpha \ast g_\varepsilon (x)  =
\int_{\{- \varepsilon \le y \le 0 \}}
\widetilde F_\alpha(x-y)  g_\varepsilon(y) dy
 \le  \widetilde F_\alpha (2 \rho + \varepsilon) .
\end{eqnarray*}
Finally, for $x > 2 \rho$, since $\widetilde F_\alpha'$ is non-increasing,
if we set $\widetilde F_\alpha'(2 \rho^+):=
\lim_{x \rightarrow 2 \rho^+}F_\alpha'(x)$,
\begin{eqnarray*}
\widetilde F_\alpha \ast g_\varepsilon (x)
& = &
\widetilde F_\alpha(x) +
\int_{\{- \varepsilon \le y \le 0 \}}
(\widetilde F_\alpha(x-y)-\widetilde F_\alpha(x)) g_\varepsilon(y) dy  \\
& \le &
 \widetilde F_\alpha(x) +
\varepsilon \max_{\{x \le z \le x + \varepsilon\}} \widetilde F_\alpha'(z)  \\
& \le &
\widetilde F_\alpha(x) +
\varepsilon \widetilde F_\alpha'(2 \rho^+) .
\end{eqnarray*}
Hence, for any integer $n$, for any function
$f:\dR^n \rightarrow \dR$ and any $\varepsilon> 0$ small enough,
\begin{eqnarray*}
\int \!\! f^2 \widetilde F_\alpha \ast g_\varepsilon 
\Bigl( \frac{f^2}{\nu_\alpha^{\otimes n} (f^2)}
\Bigr) d \nu_\alpha^{\otimes n} 
& \le & 
\int \! \! f^2 \widetilde F_\alpha 
\Bigl(\frac{f^2}{\nu_\alpha^{\otimes n} (f^2)} \Bigr)
d \nu_\alpha^{\otimes n} \\
&&  
+ (\widetilde F_\alpha (2 \rho + \varepsilon) + 
\varepsilon \widetilde F_\alpha'(2 \rho^+)) 
\int \!\! f^2 d \nu_\alpha^{\otimes n}  .
\end{eqnarray*}
The claimed result follows from  Proposition \ref{prop:Fpha}, with 
$\widetilde C'= \widetilde F_\alpha'(2 \rho^+) =
\frac{\alpha -1}{\alpha \rho}(\log 2\rho)^\frac{\alpha -2}{\alpha}$.
\qed
\end{proof}

\begin{lemma} \label{lem:approx}
Let $F:\dR+ \rightarrow \dR^+$ be a continuous non-decreasing function
such that $F\equiv 0$ on $[0,2\rho]$, for some $\rho >1$, and $F>0$ on
$(2\rho,\infty)$.
Assume that $F$ is $\C{2}$ on $(2\rho,\infty)$ and that
$\lim_{x \rightarrow 2\rho^+} F'(x)$ and
$\lim_{x \rightarrow 2\rho^+} F''(x)$ exist.
Furthermore, assume that $F'' \leq 0$ on $(2\rho,\infty)$.

Let $g$ be a $\C{\infty}$ non-negative function with compact support in
$[-1,0]$ and such that $\int g(y) dy =1$.
Define $g_\varepsilon (x) = \frac1\varepsilon g(\frac{x}{\varepsilon})$
for any $\varepsilon >0$.

Assume that for some $\lambda>0$, $F$ satisfies for all $x \neq 2\rho$
$$
x F''(x) + \lambda F'(x) \ge 0 .
$$
Then, for any $\varepsilon>0$ small enough, any $x \geq 0$,
\begin{equation} \label{eq:approx}
x(F \ast g_\varepsilon)''(x) + \lambda (F \ast g_\varepsilon)'(x) \ge 0 .
\end{equation}
Here, $F \ast g_\varepsilon (x) := \int F(x-y) g_\varepsilon (y) dy$.
\end{lemma}

\begin{proof}
Note first that for any $\varepsilon>0$, $F \ast g_\varepsilon$ is a
$\C{\infty}$ function.
Fix $\varepsilon >0$.

If $x \in (0,2\rho -\varepsilon)$, then it is easy to check that
$(F \ast g_\varepsilon)'(x)= (F \ast g_\varepsilon)''(x)=0$.
Thus \eqref{eq:approx}
holds for any $x \in (0,2\rho -\varepsilon)$ and by continuity for any
$x \in [0,2\rho -\varepsilon)$.

Now fix $x \in (2\rho,\infty)$ and note that for any
$y \in {\rm supp}(g_\varepsilon) \subset [-\varepsilon,0]$, $x-y > 2\rho$. Thus
$F'(x-y)$ and $F''(x-y)$ are well defined. It follows that
$$
x(F \ast g_\varepsilon)''(x) + \lambda(F \ast g_\varepsilon)'(x)
=
\int [xF''(x-y) + \lambda F'(x-y)]g_\varepsilon(y) dy .
$$
Since $F'' \le 0$ and $y \le 0$, $xF''(x-y) \ge (x-y)F''(x-y)$. Hence,
the left hand side of the latter inequality is bounded below by
$$
\int [(x-y)F''(x-y) + \lambda F'(x-y)]g_\varepsilon(y) dy \ge 0
$$
just using our assumption on $F$. Thus  \eqref{eq:approx}
holds for any $x > 2\rho$ and it remains the case
$x \in [2\rho - \varepsilon, \rho]$.
By continuity, it is enough to deal with
$x \in (2\rho - \varepsilon, 2\rho)$.

Fix $x \in (2\rho - \varepsilon, 2\rho)$. Choose $h$ such that
$x+h < 2\rho$ and note
that if $x-y \leq 2\rho$, then $F(x-y)=0$. Hence,
\begin{eqnarray*}
&& \int \frac{F(x-y+h) - F(x-y)}{h} g(y) dy
= \\
&& \qquad \qquad\qquad \int_{-\varepsilon \le y < -(2\rho - x)}
\!\!\!\!\!\!\!
\frac{F(x-y+h) - F(x-y)}{h} g(y) dy  \\
&& \qquad \qquad\qquad  + \int_{-(2\rho - x)  \le y \le 0 }
\!\!\!\!\!\! \frac{F(x-y+h)}{h} g(y) dy .
\end{eqnarray*}
The second term in the latter equality is non-negative because
$F$ is non-negative. It follows by Lebesgue Theorem that
$$
(F \ast g_\varepsilon)'(x)
\ge
\int_{\{-\varepsilon \le y < -(2\rho - x)\}} F'(x-y)g(y)dy .
$$
The same holds for $(F \ast g_\varepsilon)''(x)$ because $F'$ is non-negative.
Now, as in the previous argument, by our hypothesis on $F$,
$x(F \ast g_\varepsilon)''(x) + \lambda (F \ast g_\varepsilon)'(x)$
is bounded below by
\begin{eqnarray*}
&&\int_{\{-\varepsilon \le y < -(2\rho - x)\}} \!\!\!\!\!\!\!
[xF''(x-y) +
\lambda F'(x-y)]g_\varepsilon(y) dy \\
&& \qquad \ge
\int_{\{-\varepsilon \le y < -(2\rho - x)\}} \!\!\!\!\!\!\!
[(x-y)F''(x-y) + \lambda F'(x-y)]g_\varepsilon(y) dy
\ge 0 .
\end{eqnarray*}
\qed
\end{proof}
\medskip

\subsection{A general perturbation argument.}

In Section \ref{III} we discussed a perturbation argument in order to
prove the hyperboundedness of $\PT{t}^{(\alpha)}$ the semi group
associated to $\nu_{\alpha}$. In the previous subsection we recovered
and improved these results by using the capacity-measure approach and
the Gross-Orlicz theory. We shall below show that one can also derive
the results in Proposition~\ref{prop:Fpha} by a perturbation argument
on $F_{\alpha}$-Sobolev inequalities (see \cite[section 4]{Cat03} for
a similar argument for usual log-Sobolev inequalities). The argument
can be easily generalized to others situations, but we shall not
develop a complete perturbation theory here.  \medskip

Recall that Lebesgue measure on $\dR^n$ satisfies a family of
logarithmic Sobolev inequalities i.e.  for all $\eta > 0$ and all $f$
belonging to $\dL^1(dx)\cap \dL^{\infty}(dx)$ such that $\int \, f^2
\, dx \, = \, 1$
\begin{equation}\label{lsLebesgue}
\int f^2  \log f^2 \, dx \, \leq \, 2\eta \, \int \, |\nabla f|^2 dx \, + \, \frac n2 \,
\log\left(\frac1{4\pi \eta}\right) \, ,
\end{equation}
see e.g. \cite{Dav} Theorem 2.2.3.

Set $\beta=2(1-\frac{1}{\alpha})$ which is less than 1. According to
Lemma \ref{lem:c1c2} in the next section $\log^{\beta}(1+x) -
\log^{\beta}(2) \leq \log x $ for $x\geq 1$. Since $F_{\alpha}(x)$ is
non positive for $x\leq 1$, it follows
\begin{eqnarray}\label{Lebesguebeta}
\int f^2 F_{\alpha}(f^2) dx & \leq & \int_{\{f^2\geq 1\}} f^2  \log f^2 \, dx \\ & \leq & \int f^2
\log f^2 \, dx + 1/e \, . \nonumber
\end{eqnarray}

Let $V$ be smooth and satisfying the conditions stated in Section
\ref{III}. Denote by $\nu_V$ the associated Boltzmann measure
($\nu_V(dx)=e^{-2V} dx$), and introduce $g=e^V f$ (remark that $\int
g^2 d\nu_V = 1$) . According to \eqref{lsLebesgue} and
\eqref{Lebesguebeta}, a simple calculation yields
\begin{eqnarray}\label{pert1}
\int \, g^2  F_{\alpha}\left(g^2 e^{-2V}\right)  d\nu_V & \leq &  2\eta \, \int \, |\nabla g|^2
d\nu_V + \frac n2 \, \log\left(\frac1{4\pi \eta}\right) \\ & + &  1/e + 2 \eta \, \int \, g^2
\left(\Delta V \, - \, |\nabla V|^2\right) \, d\nu_V \, . \nonumber
\end{eqnarray}

But since $\beta<1$, $(A+B)^\beta \leq A^\beta + B^\beta$ for positive
$A$ and $B$. Hence if $V\geq 0$
$$\log^\beta(1+g^2 e^{-2V}) + \log^\beta(e^{2V}) \geq
\log^\beta(e^{2V}+g^2) \geq \log^\beta(1+g^2) \, ,$$
while for $V\leq
0$ $$\log^\beta(1+g^2 e^{-2V}) + \log^\beta(e^{2|V|}) \geq
\log^\beta(1+g^2 e^{-2V}) \geq \log^\beta(1+g^2) \, .$$
It follows
\begin{eqnarray}\label{pert2}
\int \, g^2  F_{\alpha}(g^2)  d\nu_V & \leq &  2\eta \, \int \, |\nabla g|^2 d\nu_V + \frac n2 \,
\log\left(\frac1{4\pi \eta}\right) + 1/e \\ & + &  \int \, g^2 \left(\log^\beta(e^{2|V|}) +
2\eta\big(\Delta V \, - \, |\nabla V|^2\big)\right) \, d\nu_V \, . \nonumber
\end{eqnarray}
Finally introduce the convex conjugate function $H_{\alpha}$ of $x
\rightarrow x F_{\alpha}(x)$.

Using Young's inequality $$x y \leq \varepsilon xF_{\alpha}(x) +
H_{\alpha}(y/\varepsilon)$$
in \eqref{pert2} we obtain
\begin{equation}\label{pert3}
\int \, g^2  F_{\alpha}(g^2)  d\nu_V  \leq   \frac{2\eta}{1-\varepsilon} \, \int \, |\nabla g|^2
d\nu_V + c(n,\eta,\varepsilon) +
\end{equation}
$$+ \frac{1}{1-\varepsilon} \int
H_{\alpha}\left((1/\varepsilon)\Big((2|V|)^\beta + 2\eta\big(\Delta V
   \, - \, |\nabla V|^2\big)\Big)\right) \, d\nu_V \, . $$
We have
thus obtained
\begin{theorem}\label{th:perturb}
   Let $\nu_V$ be a Boltzmann measure defined for a smooth $V$ as in
   section \ref{III}. Denote by $H_{\alpha}$ the convex conjugate of
   $x \rightarrow x F_{\alpha}(x)$. Assume that $\nu_V$ satisfies the
   following two conditions \\
   $(i)$ there exist some $\xi >0$ and some $\lambda >0$ such
   that 
   $$
   \int H_{\alpha}\left(((2+\xi)|V|)^\beta + \lambda
   \big(\Delta V  -  |\nabla V|^2\big)\right)  d\nu_V <
   +\infty ,
   $$
   $(ii)$ $\nu_V$ satisfies a Poincar\'e inequality.

Then the conclusions of Proposition \ref{prop:Fpha} for $F_\alpha$ are
still true just replacing $\nu_\alpha$ by $\nu_V$. As a consequence
the conclusions of Corollary \ref{cor:fal} are also still true.

Both conditions (i) and (ii) are satisfied when $V$ satisfies
assumption OB(V) in section \ref{III} with $G(y)=c
|y|^{2(1-\frac{1}{\alpha})}$ for some $c$ and $V$ goes to infinity at
infinity.
\end{theorem}

\begin{proof}
   \eqref{pert3} and Hypothesis $(i)$ ensure that $\nu_V$ satisfies a
   defective homogeneous $F_\alpha$ Sobolev inequality. But it is
   easily seen that $F_\alpha$ fulfills the hypotheses of the
   Rothaus-Orlicz Lemma \ref{lem:Rothaus}. Hence (ii) and 
   Theorem~\ref{thm:tightfsob} allow to tight the homogeneous $F_\alpha$
   Sobolev inequality. But since 
   $$
   \log^\beta(1+g^2) \leq
   \log^\beta(1+\frac{g^2}{\int g^2}) + \log^\beta(1+\int g^2)
   $$
   \eqref{tensbeta} holds when we replace $\nu_\alpha$ by $\nu_V$.
   Hence we may use the tensorisation property.
   
   Finally $(i)$ is clearly implied by OB(V), while $(ii)$ follows from
   Remark~\ref{Kunz}. 
   \qed
\end{proof}

Again the situation is more delicate when dealing with
$\widetilde{F}_\alpha$.


\section{Isoperimetric inequalities}\label{sec:isoperimetrie}

In this section we show that the Orlicz-hypercontractivity property implies
isoperimetric inequalities. These results are more precise
than the concentration
inequalities derived in the previous section
(via the Beckner type inequalities).
Let us recall the basic  definitions.
Let $\mu$ be a Borel measure on $\dR^{n}$.
For a measurable set $A\subset \dR^{n}$ we define its $\mu$-boundary measure as
$$\mu_s(\partial A) =
\liminf_{h\to 0^{+}}\frac{\mu(A_h)-\mu(A)}{h},
$$
where $A_h=\{x\in \dR^{n},d(x,A)\le h\}=A+hB_2^{n}$
is the $h$-enlargement of $A$
in the Euclidean distance
(here $B_2^{n}=\{x\in\dR^{n};\; |x|\le 1\}$). The isoperimetric
function (or profile) of a probability measure on $\dR^{n}$ is
$$
I_\mu(a)=\inf\{\mu_s(\partial A);\; \mu(A)=a\},\quad a\in[0,1].
$$
We shall write $I_{\mu^k}$ for the isoperimetric function of the product
measure (on $\dR^{nk}$ the enlargements are for the Euclidean distance, that
is the $\ell_2$ combination of the distances on the factors).
Finally we set $I_{\mu^{\infty}} :=\inf_{k\ge 1} I_{\mu^{k}}$.

We follow Ledoux's approach of an inequality by Buser \cite{ledo94sapi}
 bounding from below the Cheeger constant  of a compact Riemannian
manifold in terms of its
spectral gap and of a lower bound on its curvature.
Ledoux also deduced a Gaussian
isoperimetric inequality from a logarithmic Sobolev inequality.
The argument was extended to the framework of Markov diffusion generators
 by Bakry and Ledoux \cite{BL96}.
Moreover
these authors obtained dimension free constants.
The following result  is a particular case of
\cite[Inequality (4.3)]{BL96}.
It allows to turn hypercontrativity  properties
into isoperimetric inequalities.

\begin{theorem}\label{th:isop-hyp}
Let $\mu$ be a probability measure on
$\dR^{n}$ with $d\mu(x)=e^{-V(x)}dx$ with $V''\ge 0$.
Let $(\PT{t})_{t \geq 0}$ be the corresponding semi-group with generator
$\Delta-\nabla V.\nabla$.
Then for every $t\ge 0$ and every smooth and bounded  function, one has
$$
\|f\|_2^{2} - \|\PT{t/2}f\|_2^{2}
\le \sqrt{2t} \|f\|_\infty  \int |\nabla f| d\mu.
$$
In particular (applying this to approximations of characteristic functions)
for any Borel set $A\subset \dR^{n}$ one has
$$
\mu(A)-\|\PT{t/2}\ind_A\|_2^{2}
\le \sqrt{2t}   \mu_s(\partial A).
$$
\end{theorem}

\begin{remark}
If one only assumes that $V''\ge -R \cdot Id$ for $R>0$ then the
statement is valid with  an additional factor
$(2t R/(1-\exp(-2tR)))^{1/2}$ on the right-hand side.
This factor is essentially a constant when  $t\le 1/R$.
\end{remark}

In order to exploit this result we need the following two lemmas.

\begin{lemma}
Let the measure $\mu$ and the semi-group $(\PT{t})_{t \geq 0}$ be as before.
Let $\tau$ be a Young  function, and assume
that for all $f \in \dL^2(\mu)$ one has $N_\tau(\PT{t}f) \le C \|f\|_2$.
Then for every Borel subset $A$ of $\dR^{n}$ one has
$\| \PT{t} \ind_A\|_2 \le C \mu(A) \tau^{-1}\left(\frac{1}{\mu(A)}\right)$,
where $\tau^{-1}$ stands for the reciprocal function of $\tau$.
\end{lemma}

\begin{proof}
Since $\PT{t}$ is symmetric for $\mu$, one gets by duality that $\PT{t}$ maps
the dual of $(\dL_\tau(\mu),N_\tau)$ into $\dL^2(\mu)$ with norm at most $C$.
So for every $A$,
$\|\PT{t} \ind_A\|_2 \le C \| \ind_A \|_{\tau^*}$. Recall that the latter
norm is
\begin{eqnarray*}
\| \ind_A\|_{\tau^*}
&=& \sup\left\{\int_A g d\mu;\; \int \tau(g) d\mu \le 1 \right\} \\
&=& \sup\left\{\int_A g d\mu;\; \int_A \tau(g) d\mu\le 1\right\}
 =\mu(A) \tau^{-1}(1/\mu(A)).
\end{eqnarray*}
Indeed  Jensen inequality yields
$\int_A \tau (g)\frac{d\mu}{\mu(A)}
 \ge \tau\left(\int_A g \frac{d\mu}{\mu(A)}\right),$
which is tight for $g=\ind_A \tau^{-1}(1/\mu(A))$.
\qed
\end{proof}

\begin{lemma}
  Let $F: \dR^+ \to \dR$
  be a non-decreasing function with $F(1)=0$,
  and continuous
  on $[1,+\infty)$. Consider
  for $q,x\ge 0$, the function $\tau_q(x)=x^{2}e^{qF(x^2)}$. Assume
  there exists constants $c_1,c_2$ such
  that for all $x\ge 1$ one has $F(x)\le c_1 \log x$ and $F(x^2)\le c_2 F(x)$.
  Then for all $q\in[0,1/c_1]$ one has
  $$
  \tau_q^{-1}(y)\le \sqrt y  e^{-\frac{q}{2c_2}F(y)}, \qquad y\ge 1.
  $$
\end{lemma}

\begin{proof}
 Set $\Theta(x)=\exp(-qF(x)/(2c_2))$. Setting $x=\tau_q(y)$, $y\ge 1$
 the claimed inequality can be rephrased as:
 $$
 x\le \tau_q(x)^{\frac12}
 \Theta(\tau_q(x))=
 x e^{\frac{q}{2}F(x^2)} e^{-\frac{q}{2c_2}F(x^{2}\exp(qF(x^2)))},
 \quad x\ge 1.
 $$
 This is equivalent to $F(x^{2}\exp(qF(x^{2})))\le c_2 F(x^{2})$. The latter
 follows from the hypotheses: for $q\le 1/c_1$, $F(x^{2}\exp(qF(x^{2})))\le
 F(x^{2+2qc_1})\le F(x^{4})\le c_2 F(x^{2})$.
 \qed
\end{proof}

\begin{theorem}\label{th:isop-hyp2}
Let $\mu$ be a probability measure on $\dR^{n}$ with $d\mu(x)=e^{-V(x)}dx$ and
$V''\ge 0$.
Assume that the corresponding semi-group $(\PT{t})_{t \geq 0}$
with generator $\Delta  -\nabla V \cdot \nabla$ satisfies for every
$t\in [0,T]$ and every function in $\dL^2(\dR^{n},\mu)$,
$$
N_{\tau_{kt}}(\PT{t} f)\le C \|f\|_2,
$$
where $k>0,C\ge1$ and for $q\ge 0,x\in \dR$,
$\tau_q(x)=x^{2}\exp(qF(x^{2}))$. Here
$F: [0,\infty)\to \dR$ is non-decreasing and satisfies $F(1)=0$,
and for $x\ge 1$, $F(x)\le c_1 \log x$, $F(x^{2})\le c_2 F(x)$.
Then if $A \subset \dR^{n}$ has small measure in the sense that
$F(1/\mu(A))\ge c_2 \log(2C^2)/\min(kT,1/c_1)$
one has the following isoperimetric inequality:
$$
\mu_s(\partial A)
\ge
\frac14 \left(\frac{k}{c_2 \log(2C^2)}\right)^{\frac12}
\mu(A)  F\left(\frac{1}{\mu(A)}\right)^{\frac12}.
$$
The symmetric inequality holds  for large sets: if
$\displaystyle F\left(\frac{1}{1-\mu(A)}\right)
\ge\frac{ c_2 \log(2C^2)}{\min(kT,1/c_1)}$, then
$$
\mu_s(\partial A)\ge \frac14 \left(\frac{k}{c_2 \log(2C^2)}\right)^{\frac12}
 (1-\mu(A))  F\left(\frac{1}{1-\mu(A)}\right)^{\frac12}.
$$
\end{theorem}

\begin{proof}
We combine the above results and choose an appropriate value of the time
parameter. If $t\le \min(2T,2/(kc_1)$ then
\begin{eqnarray*}
\mu_s(\partial A)
&\ge &
\frac{\mu(A)-\|\PT{t/2}\ind_A \|_2^{2}}{\sqrt{2t}}\\
&\ge&
\frac{\mu(A)-\left(C\mu(A)\tau_{kt/2}^{-1}\Big(\frac{1}{\mu(A)} \Big)
\right)^{2}}{\sqrt{2t}}\\
&\ge&
\mu(A) \frac{1-C^2 \exp\Big(-\frac{kt}{2c_2}
F\left(\frac{1}{\mu(A)}\right)\Big)}{\sqrt{2t}}.
\end{eqnarray*}
At this point we wish to choose $t$ so that
$\frac12=C^2 \exp\Big(-\frac{kt}{2c_2} F\big(\frac{1}{\mu(A)}\big)\Big)$.
This is compatible with the
condition $t\le \min(2T,2/(kc_1)$ provided
$F(1/\mu(A))\ge c_2 \log(2C^2)/\min(kT,1/c_1)$.
Under this condition, this value
of time yields the claimed isoperimetric inequality for small sets.
For large sets
note that applying the functional inequality of
Theorem~\ref{th:isop-hyp} to suitable
approximations of the characteristic function of $A^{c}$ gives $\sqrt{2t}
\mu_s(\partial A) \ge \mu(A^{c})-\|\PT{t/2} \ind_{A^{c}}\|_2^{2}
$, so the study of small sets apply.
\qed
\end{proof}

\begin{remark}
  Under the weaker assumption $V''\ge -R$ for $R>0$ we have similar
  results with constants depending on $R$.
\end{remark}

\begin{remark}
  Under specific assumptions on $F$ we have shown that
  $\mathrm{Cap}_\mu(A)\ge \mu(A) F(1/\mu(A))$ for all $A$
  implies continuity of the
  semigroup in the Orlicz scale $\tau_q(x)=x^{2}\exp(qF(x^{2}))$, which
  implies, at least for small sets,
  $\mu_s(\partial A)\ge K \mu(A) \sqrt{F(1/\mu(A))}$.
  Note the analogy between these relations and also the inequality
  $$
  \mu_s(\partial A)\ge \mathrm{Cap}_\mu^{(1)}(A)
  :=\inf\left\{ \int \! |\nabla f| d\mu;
  f\ge \ind_A  \, \mathrm{and}  \, \mu(f=0)\ge 1/2\right\}\!.
  $$
\end{remark}

The previous theorem provides a lower bound on the
isoperimetric profile for small
and large values of the measure only. We deal with the remaining values,
 away from 0 and 1, by means of  Cheeger's inequality.
 The dimension free version of Buser's inequality
for diffusion generator, contained in the work of
Bakry and Ledoux allows to derive
Cheeger's inequality from Poincar\'e inequality.

\begin{theorem}\label{th:cheeger}
Let $\mu$ be a probability measure on $\dR^{n}$
with $d\mu(x)=e^{-V(x)}dx$ and $V''\ge 0$.
Assume that the corresponding semi-group $(\PT{t})_{t \geq 0}$ with generator
$\Delta  -\nabla V \cdot \nabla$ satisfies the following
Poincar\'e  inequality: for all $f$
$$
\lambda \int (f-\mu(f))^{2}d\mu \le \int |\nabla f|^{2}d\mu.
$$
Then for every Borel set $A\subset \dR^{n}$ one has
$$
\mu_s(\partial A) \ge\frac{1-e^{-1}}{\sqrt2} \sqrt{\lambda} \mu(A)(1-\mu(A)).
$$
\end{theorem}

The argument is written in the setting of Riemannian
manifolds in \cite[Theorem~5.2]{ledo04sgls}.
We sketch the proof for completeness.

\begin{proof}
The spectral gap inequalities classically implies the exponential decay of
the norm of $\PT{t}$ on the space of zero mean. Therefore
\begin{eqnarray*}
\|\PT{t/2} \ind_A\|_2^{2}
&=&
\|\PT{t/2} \mu(A)\|_2^{2} + \|\PT{t/2} (\ind_A-\mu(A))\|_2^{2}\\
&\le &
\mu(A)^{2} + e^{-\lambda t} \|\ind_A-\mu(A)\|_2^{2}\\
& = &
\mu(A)^{2} + e^{-\lambda t}\mu(A)(1-\mu(A)).
\end{eqnarray*}
By Theorem~\ref{th:isop-hyp}, one has
$$
\sqrt{2t} \mu_s(\partial A) \ge (1-e^{-\lambda t}) \mu(A)(1-\mu(A)).
$$
Choosing $t=1/\lambda$ concludes the proof.
\qed
\end{proof}

Finally we apply the previous results to infinite products of exponential measures:
$m_\alpha(dx)=\exp(-|x|^{\alpha})/(2\Gamma(1+1/\alpha)) dx,\, x\in\dR$. For technical reasons, we
also consider the measures $\nu_\alpha$ defined in section \ref{examples} up to the irrelevant
constant 2. They also have a log-concave density, but more regular. The isoperimetric function of a
symmetric log-concave density on the line (with the usual metric) was calculated by Bobkov
\cite{bobk96ephs}. He showed that half-lines have minimal boundary among sets of the same measure.
Since the boundary measure of $(-\infty,t]$ is given by the density of the measure at $t$, the
isoperimetric profile is easily computed. They are readily compared to the functions
$$
L_\alpha(t)= \min(t,1-t)
\log^{1-\frac1\alpha}\left(\frac{1}{\min(t,1-t)}\right).
$$
We omit the details, some of them are written in \cite{bart01lcbe}.

\begin{lemma}
There are constants $k_1,k_2$ such that for all $\alpha\in[1,2]$, $t\in[0,1]$
one has
$$
k_1  L_\alpha(t)\le I_{m_\alpha}(t) \le  k_2  L_\alpha(t),
$$
$$
k_1  L_\alpha(t)\le I_{\nu_\alpha}(t) \le  k_2  L_\alpha(t).
$$
\end{lemma}

Our goal is to show the following infinite
dimensional isoperimetric inequality.

\begin{theorem}\label{th:iso}
There exists a constant $K>0$ such that
for all  $\alpha\in[1,2]$ and $t\in[0,1]$, one has
$$
I_{\nu_\alpha^{\infty}}(t) \ge K L_\alpha(t).
$$
\end{theorem}

Since $I_{\nu_\alpha^{\infty}}\le I_{\nu_\alpha}\le k_2 L_\alpha$,
we have, up to  a constant, the value of the
isoperimetric profile of the infinite product.

\begin{proof}
   As shown in Corollary~\ref{cor:fal} of Section~\ref{examples} the
   semi-group associated to $\nu_\alpha^{\otimes n}$ is
   Orlicz-hyperbounded. Thus we may  apply Theorem~\ref{th:isop-hyp2}
   with $F=F_{\beta(\alpha)}$ defined in \eqref{eq:falpha}, and
   $\beta(\alpha)=2(1-1/\alpha)$ and get  an isoperimetric
   inequality for small and large sets, with constants independent of
   the dimension $n$. This step requires to check a few properties of
   the function $F_{\beta(\alpha)}$. They are established in the
   following Lemma~\ref{lem:c1c2}. More precisely there are constants
   $K_1,K_2>0$ independent of $\alpha$ and $n$ such that
\begin{equation}\label{eq:isobeta}
 I_{\nu_\alpha^{\otimes n}}(t)
 \ge
 K_1 \min(t,1-t) \left[\log^{\beta(\alpha)}
 \left(1+\frac{1}{\min(t,1-t)}\right)-\log^{\beta(\alpha)}(2)\right]^{\frac12}
\end{equation}
provided
$$
\left(\log^{\beta(\alpha)}\left(1+\frac{1}{\min(t,1-t)}\right)
-\log^{\beta(\alpha)}(2)\right)^{\frac12} \ge K_2.
$$
We can prove \eqref{eq:isobeta} in the remaining range as well. Indeed, it
is plain that
$$
\sup_{x>0} \nu_\alpha([x,+\infty))
\int_0^{x}\frac{1}{\rho_{\nu_\alpha}} \le M,
$$
so that the measures $(\nu_\alpha)_{\alpha\in[1,2]}$ satisfy a
Poincar\'e inequality with a uniform constant. The latter inequality
has the tensorisation property, so the measures $\nu_\alpha^{\otimes
  n}$ also share a common Poincar\'e inequality. By
Theorem~\ref{th:cheeger}, there exists a constant $K_3>0$ such that
for all $n$, all $\alpha \in [1,2]$ and all $t\in [0,1]$
\begin{equation}\label{eq:cheeger}
 I_{\nu_\alpha^{\otimes n}}(t) \ge K_3 \min(t,1-t).
\end{equation}
Since the exponential measure
has a spectral gap, the latter argument reproves,
with a slightly worse constant, the result of \cite{BH97}.
Now assume that
$$
\log^{\beta(\alpha)}\left(1+\frac{1}{\min(t,1-t)}\right)
-\log^{\beta(\alpha)}(2) < K_2^{2},
$$
then
\begin{eqnarray*}
I_{\nu_\alpha^{\otimes n}}(t) & \ge &
\frac{K_3}{K_2} \min(t,1-t)K_2 \\
& \ge &
\frac{K_3}{K_2}
\min(t,1\!-\!t)\Bigl[\log^{\beta(\alpha)}
\bigl(1\!+\!\frac{1}{\min(t,1\!-\!t)}\bigr)
-\log^{\beta(\alpha)}(2)\Bigr]^\frac12 \!\!\!.
\end{eqnarray*}
So Inequality~\eqref{eq:isobeta} is valid for all $t$ provided one replaces
$K_1$ by $K_4:=\min(K_1,K_3/K_2)$.
Finally, the uniform Cheeger inequality \eqref{eq:cheeger}, implies that
$$
\frac{1}{K_3}I_{\nu_\alpha^{\otimes n}}(t) \ge \log^{\beta(\alpha)/2}(2) \min(t,1-t).
$$
Adding up this  relation to
\begin{eqnarray*}
 &&\frac{1}{K_4}I_{\nu_\alpha^{\otimes n}}(t)
 \ge
 \min(t,1\!-\!t) \!\left[\log^{\beta(\alpha)}
 \!\Big(1\!+\!\frac{1}{\min(t,1\!-\!t)}\Big)-
 \log^{\beta(\alpha)}(2)\right]^{\frac12}\\
 &&\qquad \ge
 \min(t,1-t) \left[\log^{\beta(\alpha)/2}
 \Big(1+\frac{1}{\min(t,1-t)}\Big)-\log^{\beta(\alpha)/2}(2)\right]
\end{eqnarray*}
yields the claimed inequality. This manipulation was important in order
to get  a non-trivial inequality when $\alpha$ tends to 1, {\it i.e.}
when $\beta(\alpha)$
tends to 0.
\qed
\end{proof}

The following technical result was used in the above proof.

\begin{lemma}\label{lem:c1c2}
   Let $\beta\in [0,1]$ then for all $x\ge 1$ one has
   \begin{equation}
      \label{eq:borne1}
      \log^{\beta}(1+x)-\log^{\beta}(2)\le \log x,
   \end{equation}
   \begin{equation} \label{eq:borne2}
       \log^{\beta}(1+x^{2})-\log^{\beta}(2)
       \le
       8\left(  \log^{\beta}(1+x)-\log^{\beta}(2)\right).
   \end{equation}
\end{lemma}

\begin{proof}
Note that \eqref{eq:borne1} is an equality for $x=1$.
It is enough to prove the
   inequality between derivatives, that is
$\beta \log^{\beta-1}(1+x)/(1+x)\le 1/x$
   for $x\ge 1$.
   If $x\ge e-1$ then $\log^{\beta-1}(1+x)\le 1$ and the inequality is obvious.
   If $x<e-1$, then $\log^{\beta}(1+x)\le 1$, therefore
    $$
\beta \frac{\log^{\beta-1}(1+x)}{(1+x)} \le\frac{1}{(1+x)\log(1+x)}
\le \frac{1}{x}.
$$

Next we address \eqref{eq:borne2}. One easily checks that for $A\ge B\ge 1$
the map $\beta>0\mapsto (A^{\beta}-1)/(B^{\beta}-1)$ is non-decreasing.
Applying this to $A=\log(1+x^{2})/\log(2)$ and
$B=\log(1+x)/\log(2)$ shows that
it is enough to prove \eqref{eq:borne2} for $\beta=1$. Let $x\ge 1$, since
$1+x^{2}\le (1+x)^{2}$ one has
\begin{eqnarray*}
   \log(1+x^{2})-\log(2)
   &\le&
   2\log(1+x)-\log(2) \\
   &=&
   2\left(\log(1+x)-\log(2)\right)
   +\log(2).
\end{eqnarray*}
If $x\ge 3$ then $\log(1+x)-\log(2)\ge \log(2)$
and the claimed inequality is proved.
For $x\in(1,3]$, we use the fundamental relation of calculus.
It provides $t_1\in(1,9)$
and $t_2\in(1,3)$ with
$$
\log(1+x^{2})-\log(2)=(x^{2}-1)\frac{1}{1+t_1}\le  2(x-1)
$$
and
$$
\log(1+x)-\log(2)=(x-1)\frac{1}{1+t_2} \ge  (x-1)/4.
$$
So the ratio is bounded from above by 8. A smarter choice than 3 would
give a better result.
\qed
\end{proof}

\begin{remark}
According to Theorem \ref{th:perturb} the conclusion of Theorem \ref{th:iso} is still true with the
same $L_\alpha$ when replacing $\nu_\alpha$ by $\nu_V$, provided $V$ is convex and the hypotheses in
Theorem \ref{th:perturb} are fulfilled.
\end{remark}

We conclude the paper with consequences of Theorem~\ref{th:iso}.
The first one is a comparison theorem. It could be stated in a more
general framework of metric probability spaces satisfying a smoothness
 assumption (see {\it e.g.} \cite{bart02lsmi}). For simplicity we write it in
the setting of Riemannian manifolds where the definition of isoperimetric
profile given in the beginning of the section applies.

\begin{theorem}\label{th:isoR}
Let $(X,d,\mu)$ be a Riemannian manifold, with the geo\-desic metric, and
a probability measure which has a density with respect to the volume.
On the product manifold we consider the geodesic distance,
which is the $\ell_2$
combination of the distances on the factors.
There exists a universal constant $K>0$ such that if for some $c>0$,
$\gamma \in [0,\frac12]$ and
all $t\in[0,1]$ one has
$$
I_\mu(t) \ge c \min(t,1-t) \log^{\gamma}\left( \frac{1}{\min(t,1-t)}\right),
$$
then for all $n\ge 1$, $t\in [0,1]$ one has
$$
I_{\mu^{\otimes n}}(t)\ge \frac{c}{K}  \min(t,1-t) \log^{\gamma}\left( \frac{1}{\min(t,1-t)}\right).
$$
\end{theorem}

\begin{remark}
   This provides a scale of infinite dimensional isoperimetric inequalities.
   Both ends of the scale where previously known.
A standard argument based
   on the central limit theorem shows that if
$\mu$ is a measure on $\dR$ with second
  moment then $\inf_n I_{\mu^{\otimes n}}$ is dominated by a multiple of the Gaussian
  isoperimetric function, which is comparable to
  $\min(t,1-t) \log^{1/2}(1/\min(t,1-t))$.
On the other hand an argument of Talagrand
  \cite{Tal91}
  shows that the weakest possible dimension free concentration result for $\mu$
  implies that it has at most exponential tails.
The isoperimetric function of the
  exponential density is $\min(t,1-t)$. So the above scale covers the whole
  range of infinite isoperimetric inequalities.
Of course finer scales could be
  obtained from our methods, with more effort.
\end{remark}

\begin{remark}
   A similar statement was proved in \cite{bart04idii} for the case when
  the distance on the product space is the $\ell_\infty$ combination of the
  distances on the factors ({\it i.e.} the maximum).
This case was much easier due
  to the product structure of balls in the product space. Also, this notion
  leads to bigger enlargement, and the scale of infinite dimension behaviour
  was larger,
  the values $\gamma\in[0,1]$ being allowed.
\end{remark}

\begin{proof}[of Theorem \ref{th:isoR}]
 The hypothesis implies that $I_\mu \ge \frac{c}{k_2} I_{\nu_\alpha}$ for
 $\alpha=1/(1-\gamma)\in[1,2]$. Theorem~10 in \cite{bart02lsmi} asserts that
 among measures having the same concave isoperimetric behaviour,
the even log-concave
 one minimizes the isoperimetric profile for the product measures, 
see also \cite{ros01ip}. So we
 have $I_{\mu^{\otimes n}}\ge   \frac{c}{k_2} I_{\nu_\alpha^{\otimes n}}$.
By the previous results
 $ I_{\nu_\alpha^{\otimes n}}\ge K L_\alpha$ and the proof is complete.
\qed
\end{proof}

 The second consequence that we wish to put forward deals  with the measures
 $dm_\alpha(x)=\exp(-|x|^{\alpha})dx/(2\Gamma(1+1/\alpha))$, $\alpha\in[1,2]$.
It shows that among sets of prescribed measure for $m_\alpha^{\otimes n}$ in $\dR^{n}$,
 coordinate half-spaces have enlargements of minimal measure, up to a universal
factor.
The result was known for $\alpha\in{1,2}$.

\begin{theorem}
   There exists a universal constant $K$ such that for every $\alpha\in[1,2]$,
   $n\ge 1$ and every Borel
   set $A\subset \dR^{n}$, if $m_\alpha^{\otimes n}(A)=m_\alpha((-\infty,t])$ then for
   $h\ge 0$,
   $$
   m_\alpha^{\otimes n}\Big(A+hB_2^{n}\Big)
   \ge
   m_\alpha\left(\left(-\infty,t+\frac{h}{K}\right]\right).
$$
\end{theorem}

\begin{proof}
   This fact is proved by integrating the inequality
   $I_{m_\alpha^{\otimes n}}\ge \frac{I_{m_\alpha}}{K}$ which provides a
   similar information about boundary measure (this corresponds to
   infinitesimal enlargements). This isoperimetric ine\-quality is a
   consequence of the fact that $I_{m_\alpha}$ is comparable to
   $I_{\nu_\alpha}$. The comparison theorem of \cite{bart02lsmi}
   implies that $I_{m_\alpha^{\otimes n}}$ is larger than a universal
   constant times $I_{\nu_\alpha^{\otimes n}}\ge K L_\alpha \ge
   \frac{K}{k_2} I_{m_\alpha}$. \qed
\end{proof}


\begin{acknowledgement}
We wish  to thank  Sergey Bobkov, 
Thierry Coulhon, Arnaud Guillin, G\'erard Kerkyacharian,
 Michel Ledoux, Yves Raynaud and Boguslaw Zegarlinski
 for fruitful discussions and remarks on the topic of this paper.
\end{acknowledgement}


\bibliographystyle{plain}


\end{document}